\documentclass[11 pt, reqno]{article}

\usepackage{amsmath,amsfonts,amssymb,amsthm}
\usepackage[colorlinks=true,hyperindex=true]{hyperref}
\usepackage{cancel,bbm}

\usepackage{mathabx} 
\usepackage{comment}
\usepackage{enumitem}
\usepackage{cite}
\usepackage{showlabels}

\usepackage{chngcntr}
\counterwithin*{equation}{section}

\usepackage[left=3cm, right=3cm, top=3 cm, bottom= 3 cm]{geometry}
\usepackage{color}
\usepackage{fancyhdr}
\usepackage{latexsym}

\usepackage{bm}

\newtheorem{ccounter}{ccounter}[section]
\newtheorem{thm}[ccounter]{Theorem}
\newtheorem{lem}[ccounter]{Lemma}
\newtheorem{cor}[ccounter]{Corollary}
\newtheorem{defn}[ccounter]{Definition}
\newtheorem{prop}[ccounter]{Proposition}
\newtheorem{ass}[ccounter]{Assumption}
\newtheorem{ex}[ccounter]{Example}

\def\bet{\begin{thm}}
\def\eet{\end{thm}}
\def\bel{\begin{lem}}
\def\eel{\end{lem}}
\def\bas{\begin{ass}}
\def\eas{\end{ass}}
\def\bec{\begin{cor}}
\def\eec{\end{cor}}
\def\bed{\begin{defn}}
\def\eed{\end{defn}}
\def\bep{\begin{prop}}
\def\eep{\end{prop}}
\def\beq{\begin{equation}}
\def\eeq{\end{equation}}
\def\proof{\noindent {\bf Proof.}\ \ }
\def\bea{\begin{equation*}}
\def\eea{\end{equation*}}
\def\tr{\mathrm{tr}}
\def\bex{\begin{ex}}
\def\eex{\end{ex}}
\def\remark{\noindent{\bf Remark. }}

\def\rr{\mathbb{R}}
\def\zz{\mathbb{Z}}
\def\cc{\mathbb{C}}
\def\1{\boldsymbol{1}}
\def\Im{\mathrm{Im}}
\def\Re{\mathrm{Re}}
\def\e{\mathrm{e}}
\def\i{\mathrm{i}}
\def\del{\partial}
\def\d{\mathrm{d}}
\def\eps{\varepsilon}
\renewcommand\leq\varleq
\renewcommand\geq\vargeq
\def\ee{\mathrm{E}}

\def\F{\mathcal{F}}
\def\O{\mathcal{O}}

\def\ee{\mathbb{E}}

\def\pp{\mathbb{P}}

\def\msc{m_{\mathrm{sc}}}
\def\rhosc{\rho_{\mathrm{sc}}}

\def\mfa{\mathfrak{m}}

\def\J{\mathcal{J}}

\def\Var{\mathrm{Var}}
\def\nn{\mathbb{N}}

\def\bw{\boldsymbol{w}}
\def\bu{\boldsymbol{u}}

\def\Cov{\mathrm{Cov}}
\def\be{\mathbf{e}}

\def\mfa{\mathfrak{a}}
\def\Oma{\Omega_{\mathfrak{a}}}
\def\tilf{\tilde{f}}
\def\ea{e_{\mathfrak{a}}}

\def\mfb{\mathfrak{b}}

\def\E{\mathcal{E}}

\def\tt{\mathbb{T}}

\def\ea{e_{\mathfrak{a}}}
\def\ui{\underline{i}}

\def\Oma{\Omega_{\mfa}}
\def\G{\mathcal{G}}
\def\M{\mathcal{M}}
\def\bz{\boldsymbol{z}}
\def\sd{\prec}
\def\Osd{\mathcal{O}_{\prec}}
\def\Gj{\mathcal{G}^{(j)}}
\def\sui{\{ \underline{i} \}}
\def\delrr{\del^{(\rr)}}
\def\LSS{LSS}
\def\Hcal{\mathcal{H}}
\begin{document}

\begin{table}
\centering

\begin{tabular}{c}

\multicolumn{1}{c}{\parbox{12cm}{\begin{center}\Large{\bf Some estimates for generalized Wigner matrix linear spectral statistics}\end{center}}}\\
\\
\end{tabular}
\begin{tabular}{  c   }
Benjamin Landon
 \\
   \\  
 \small{University of Toronto} 
 \\
 \small{Department of Mathematics} 
  \\
 \small{\texttt{blandon@math.toronto.edu}} \\ 
\end{tabular}
\\
\begin{tabular}{c}
\multicolumn{1}{c}{\today}\\
\\
\end{tabular}

\begin{tabular}{p{15 cm}}
\small{{\bf Abstract:}  We consider the characteristic function of linear spectral statistics of generalized Wigner matrices. We provide an expansion of the characteristic function with error $\O ( N^{-1})$ around its limiting Gaussian form, and identify sub-leading non-Gaussian corrections of size $\O (N^{-1/2})$. Prior expansions with this error rate held only for Wigner matrices; only a weaker error rate was available for more general matrix ensembles. We provide some applications.}
\end{tabular}
\end{table}

\tableofcontents

\section{Introduction}

This work is concerned with the behavior of eigenvalues of generalized Wigner matrices. Generalized Wigner matrices $H$ are $N \times N$ self-adjoint random matrices with independent centered entries, such that the matrix of variances $S_{ij} := \ee[ H_{ij}^2]$ is doubly stochastic. Wigner's semicircle law states that as $N \to \infty$, the eigenvalues of $H$ (in the scaling such that the entries $S_{ij}$ are order $\O(N^{-1})$) satisfy,
\beq
\lim_{N \to \infty} \frac{1}{N} \sum_{i=1}^N \delta_{\lambda_i} (E) = \frac{ \sqrt{(4-E^2)_+}}{2 \pi} \d E =: \rhosc(E) \d E
\eeq
where the convergence is, e.g., weakly almost surely. It is a classical result that the fluctuations around the semicircle law are Gaussian in that for sufficiently nice $f : \rr \to \rr$ we have that the random varaible
\beq
\LSS(f) := \sum_{i=1}^N f ( \lambda_i ) - N \int f(E) \rhosc (E) \d E 
\eeq
converges weakly to a Gaussian as $N \to \infty$ with a certain mean and variance $E (f)$ and $V (f)$ that will be defined below. This convergence can be quantified in several ways. In the present work, will be interested in deriving an expansion of the characteristic function of the form,
\beq \label{eqn:intro-expansion}
\varphi (\lambda) := \ee[ \e^{ \i \lambda \LSS (f) } ] = \e^{ - \lambda^2 V (f)/2 + \i \lambda \E (f) + N^{-1/2} P_f ( \lambda ) } + \O ( N^{-1+\eps} ),
\eeq
for any $\eps >0$. 
This is the content of our main result, Theorem \ref{thm:main}  below. Above, $N^{-1/2} P_f ( \lambda)$ should  be viewed as a sub-leading correction of size $\O ( N^{-1/2})$; without such a term, the above error rate would be $\O ( N^{-1/2})$. In our result $P_f ( \lambda)$ explicit in that it depends only on the first few moments of the matrix entries of $H$ and the coefficients of $f$ in its expansion in Chebyshev polynomials. It contains both a non-Gaussian cubic correction $ \i \lambda^3$ and a subleading correction to the variance $V(f)$.

Our motivation in deriving the expansion \eqref{eqn:intro-expansion} lies in our prior work with Sosoe \cite{landon2022almost}. In that work we derived almost-optimal regularity conditions on the function $f$ so that the Gaussian convergence holds. A technical input into that work was the expansion \eqref{eqn:intro-expansion} for Wigner matrices (which we established in \cite{landon2022almost}). By obtaining \eqref{eqn:intro-expansion} for generalized Wigner matrices, in a forthcoming revision to \cite{landon2022almost} we will extend the result to this larger class of matrices.

In Section \ref{sec:max} we give a modest application of the CLT for generalized Wigner matrices. We extend one of the results of the recent work \cite{BLZ} of Bourgade, Lopatto and Zeitouni. In that work they computed the leading order of the maximum of Wigner matrix characteristic polynomials; we extend this, as well as their optimal rigidity result, to the generalized Wigner class. The bulk of the method of that work holds for generalized Wigner matrices; there is only a minor input they require for LSS, and we show how to provide that input for generalized Wigner matrices.

\subsection{Background}

The CLT for LSS is a classical result of random matrix theory, and a vast literature has developed.  We only review a few works and refer the reader to the references in these papers for a more complete bibliography. See, e.g., the works  \cite{lytova2009central,khorunzhy19951,guionnet2002large,khorunzhy1996asymptotic,bai2009clt,
bai2005convergence,lodhia2015mesoscopic}.

Direct analysis of the characteristic function was initiated by Shcherbina \cite{shcherbina2011central} who applied Stein's method to prove that it converged point-wise to that of a Gaussian. In \cite{meso} the author, together with Sosoe, adapted this approach together with the local semicircle law to prove the convergence of LSS even on mesoscopic scales (i.e., to functions $f$ of the form $f(x) = F(N^{\alpha} (x-E))$ for some $0 < \alpha < 1$). The approach of the present paper is also based on Stein's method.

The two closest precedents to our work are \cite{bao2023quantitative} and \cite{li2021fluctuations}.  Bao and He \cite{bao2023quantitative} studied the rate of convergence of LSS to their limiting Gaussian distributions in Kolmogorov-Smirnov distance; they find optimal rates of convergence in this metric by establishing upper and lower bounds. In terms of the characteristic function, they remove the contribution of the diagonal part of $H$ to $\LSS(f)$, 
 which ensures an error of $\O ( N^{-1})$ but with no $P_f ( \lambda)$ term on the RHS. However, their analysis is restricted to the smaller class of Wigner matrices (i.e., when the matrix $S$ is constant) and their result does not seem to imply an expansion for $\varphi ( \lambda)$ directly (i.e., for the full matrix) with error rate $\O ( N^{-1})$. However, their analysis also allows them to control the characteristic function for very large $\lambda \approx N$. 

Secondly, for generalized Wigner matrices, an estimate similar to \eqref{eqn:intro-expansion} but with a polynomial error rate $\O (N^{-c})$ was derived by Li and Xu in \cite{li2021fluctuations} and also for test functions on mesoscopic scales (with an error that degenerates as the scale shortens). Other prior work on CLTs for matrices that do not have a constant variance profile include \cite{adhikari2021linear,chatterjee2009fluctuations,diaz2020global} as well as the works \cite{az,li2013central,jana2016fluctuations,shcherbina2015fluctuations,erdHos2015altshuler1,
erdHos2015altshuler2} on band matrices.

The main difficulty in the approach based on Stein's method to computing the characteristic function in the generalized Wigner case is that one needs to estimate the leading order contributions of terms like,
\beq \label{eqn:intro-tij}
\sum_{i,j} G_{ij}(z) G_{ji}(w) S_{ij}
\eeq
where $G(z) = (H-z)^{-1}$ is the resolvent of $H$. In the Wigner case $S_{ij} = N^{-1}$, this can be written as $\frac{1}{N}\tr \frac{ G (z) - G (w) }{z-w}$, which allows one to use the local semicircle law to estimate this quantity. However for general $S$, no such simplification is available. In \cite{li2021fluctuations}, Li and Xu derived a self consistent equation for this quantity, based on the methods of \cite{erdHos2013delocalization}. They achieved an estimate of $\O ( N^{-1/4})$ for this quantity (around its leading order deterministic contribution).\footnote{We ignore the dependence of the error rate on the imaginary parts of $z$ and $w$ for the present discussion.} However, as can be seen from the Wigner setting, the optimal error should be $\O (N^{-1})$. 

In order to improve on this estimate, we do not attempt to estimate \eqref{eqn:intro-tij} with high probability, and instead simply try to determine its expectation (or rather the expectation of $\e^{ \i \lambda \LSS ( f) }$ times this quantity), as this is what arises in Stein's method. We treat this quantity by means of the cumulant expansion.  Our treatment has its origins in the work \cite{lee2018local} of Lee-Schnelli, and adapted to proving CLTs via Stein's method in \cite{meso}.  

If one applies the cumulant expansion to an expression involving \eqref{eqn:intro-tij}, one is led to expressions involving more resolvent entries. One can then iterate this process, resulting in quantities with even more resolvent entries. If one was interested in obtaining the CLT only for global LSS, then only a few repeated cumulant expansions would be necessary (say, three). However, we will also derive the CLT for mesoscopic LSS; in order to do so, one must expand to an arbitrarily high order\footnote{The order of the expansion will be large but fixed as $N \to \infty$} adapted to the mesoscopic scale of the function. In order to do this, we use a graphical notation to keep track of which kinds of terms arise. Such graphical expansions have appeared in many works, and we are partially inspired by \cite{erdHos2013delocalization,alex2014isotropic}.

We first develop two general estimates for certain classes of graphs. The first are loops, which consist  of (the expectation of) products of terms of the form $G_{i_1 i_2 } (z_1) G_{i_2 i_3 } (z_2) \dots G_{i_k i_1} (z_k)$ for distinct $i_1, \dots, i_k$. The second are lines, of the form $G_{i_1 i_2 } (z_1) G_{i_2 i_3 } (z_2) \dots G_{i_k i_{k+1}} (z_k)$, for distinct $i_1, \dots, i_{k+1}$. By means of an iterated cumulant expansion, we obtain an estimate for these terms that is one order better than the naive estimate one would obtain from applying the entry-wise local semicircle law (see \eqref{eqn:entry-wise} below). These results are Propositions \ref{prop:main-loop-estimate} and \ref{prop:main-line-estimate} below. At each step of the iteration, we find we are able to truncate the terms in the cumulant expansion arising from cumulants of order three and higher; due to this, the expansion ``closes'' in the sense that expanding loops only leads to expressions with loops (with more resolvent entries), and expanding lines just leads to larger lines. 

After carrying out these expansions, we then proceed via Stein's method as in \cite{meso} and apply our estimates to treat the various terms that arise.  We identify the sub-leading terms $P_f (\lambda)$ and incorporate them into our expansion.

\subsection{Definition of model}

We now define our random matrix model.

\bed
A  generalized Wigner matrix is an $N \times N$ matrix $H$ s.t. the entries $\{ H_{ij} \}_{i \leq j} $ are independent centered random variables. If we denote,
\beq
S_{ij} := \ee[ |H_{ij}|^2]
\eeq
then we assume that,
\beq
\frac{1}{CN} \leq S_{ij} \leq \frac{C}{N}
\eeq
for some $C>0$ and that $\sum_i S_{ij} = 1$ for all $j$. For every $p$ we assume that there exists $C_p >0$ so that
\beq
\ee[ | \sqrt{N} H_{ij} |^p ] \leq C_p
\eeq
for all $i, j$ and $N$. We say that the generalized Wigner matrix is real symmetric if $H_{ij} \in \rr$ for all $i, j$. We say that the generalized Wigner matrix is complex Hermitian if the real and imaginary parts of the off-diagonal entries are independent and that $\ee[ H_{ij}^2 ] =0$ for all $i \neq j$.
\eed

We will denote by $\| f \|_{p, w}$ the following weighted $L^p$ norm of $f$,
\beq
\| f \|_{p, w} := \left( \int_{\rr} |f(x)|^p \frac{ \d x }{ \sqrt{ |4-x^2|}} \d x \right)^{1/p}
\eeq
\bed
We will say that a sequence of test functions $f= f_N \in C^2 ( \rr )$ is admissible if it is compactly supported in $[-5, 5]$ and the following holds. First, for all $\eps >0$ we have that,
\beq
\| f'\|_{1, w} \leq N^{\eps} 
\eeq
for $N$ large enough. Secondly, 
\beq
\| f \|_{1, w} \leq C, \qquad \| f''\|_{1, w} \leq C N^{1- \mfb}
\eeq
for some $1> \mfb >0$ and $C>0$.
\eed

Define,
\beq
\LSS (f) := \tr f (H) - N \int f(x) \rhosc (x) \d x.
\eeq

\subsection{Main result} 
We now introduce the quantities that arise in our result on the characteristic function of the random variable $\LSS (f)$. We will denote by  $\kappa_k (X)$ denote the $k$th cumulant of a random variable $X$.
\bed[Chebyshev polynomials] \label{def:cheby} We will define the $n$th Chebyshev polynomial $T_n(x)$ of the first kind by,
\beq \label{eqn:cheby-def}
T_n (2 \cos ( \theta ) ) = \cos ( n \theta ),
\eeq
for $n \geq 0$. The polynomials $T_n(x)$ form an orthogonal basis of the space $L^2 ( [-2, 2], (4-x^2)^{-1/2} \d x )$.\footnote{We note that our convention for the $T_n(x)$ differs from the usual convention, which defines the $n$th Chebyshev polynomial as $T_n(2x)$ where $T_n$ is as above (i.e., defined wrt a measure on $(-1, 1)$ instead of $(-2, 2)$).} We denote the $n$th Chebyshev coefficient of a function $f$ by,
\beq
t_n (f) := \frac{2}{ \pi} \int_{-2}^2 T_n(x) f(x) \frac{ \d x }{ \sqrt{4 -x^2}}.
\eeq
\eed 
\bed[Functionals in the CLT]
For $\beta =1, 2$ define the variance functional,
\beq \label{eqn:intro-V-def}
V_{ \beta} (f) := \frac{1}{2 \beta} \sum_{j=1}^\infty j t_j (f)^2 \tr (S^j) - \frac{2-\beta}{4} t_1(f)^2 \tr (S) + \frac{ \hat{s}_{4, \beta}}{2} t_2 (f)^2 +  \frac{ \hat{s}_3}{2 } t_2 (f) t_1 (f) 
\eeq
where $t_j (f)$ is the $j$th coefficient of $f$ in the basis of Chebyshev polynomials and
\beq
\hat{s}_{4, 1} := \sum_{i, j} \kappa_4 ( H_{ij} ), \qquad \hat{s}_{4, 2} := \sum_{i \neq j } \kappa_4 ( \Re[H_{ij} ] ) + \kappa_4 ( \Im[ H_{ij} ] ) + \frac{1}{2} \sum_i \kappa_4 (H_{ii} ) 
\eeq
and
\beq
\hat{s}_3 := \sum_i \kappa_3 (H_{ii} ).
\eeq
Denote also,
\begin{align}
&\E_\beta (f) := \hat{s}_{4, \beta} \frac{1}{ 2 \pi} \int_{-2}^2 f(x) \frac{ x^4 - 4 x^2 +2 }{ \sqrt{4-x^2}} \d x + \hat{s}_3 \frac{1}{ 8 \pi} \int_{-2}^2 f(x) \frac{ x^3 - x^2 - 2 x +4 }{ \sqrt{4-x^2}} \d x \notag\\
&+ \1_{ \{ \beta=1 \}} \bigg\{ \tr (S) \frac{1}{2 \pi} \int_{-2}^2 f(x) \frac{ 2 -x^2}{ \sqrt{4-x^2}} + \frac{ f(2) + f(-2)}{4} \notag\\
&+ \frac{1}{ \pi} \int_{-2}^2 \frac{ f(x)}{\sqrt{4-x^2}} ( \Re[ \msc(x + \i 0)^2 \tr (S (1- \msc(x+ \i 0)^2 S )^{-1} ) ] ) \d x \bigg\} 
\end{align}
as well as,
\beq
B (f) := \frac{ \hat{s}_3}{8} t_1(f)^3. 
\eeq
\eed 
The following gives elementary properties of the functionals appearing above. In particular, it guarantees that the variance function is non-negative. It is proven in Section \ref{sec:funct-proof}.
\bep \label{prop:funct}
If $f$ is an admissible test function then,
\beq
0 \leq V_\beta (f) \leq C \| f'\|_{1, w}^2 \log (N)^2
\eeq
for some $C>0$, and $| \E_\beta (f) | \leq C( \| f \|_{1, w}  + |f(2)| + |f(-2)|)$. Moreover, $|B (f) | \leq C N^{-1/2} \| f\|_{1, w}^3$.  
\eep

The following is our main result. It is proven in Section \ref{sec:main-proof}. 
\bet \label{thm:main}
Let $H$ be a real symmetric or complex Hermitian generalized Wigner matrix. If there is a $c>0$ so that $V_\beta(f) \geq c$ then for any $\eps >0$ we have that for all $|\lambda| \leq N^{-\eps} \sqrt{ N \|f'' \|_{1}^{-1} }$
that
\begin{align}
\ee\left[ \exp ( \i \lambda \LSS(f) ) \right] &= \e^{ - \lambda^2 V_\beta (f)/2 + \i \lambda^3 B(f)/3 + \i \lambda \E_\beta (f) } \notag\\
+&   \O \left( \| f''\|_{1,w}N^{-1} ( 1 + |\lambda| ) + N^{-1} |\lambda|^2 \right).
\end{align}
Otherwise, 
\begin{align}
\ee\left[ \exp ( \i \lambda \LSS(f) ) \right] &= \e^{ - \lambda^2 V_\beta (f)/2 + \i \lambda^3 B(f)/3 + \i \lambda \E_\beta (f) } \notag\\
&+ \O ( N^{-1} ( 1 + |\lambda| )^4 + (1 + |\lambda| )^3 N^{-1} ( 1 + \| f''\|_{1, w} ) ).
\end{align}
We have also that,
\beq
\ee[ \LSS (f) ] = \E_\beta (f) +  \O ( N^{-1} ( 1 + \| f'' \|_{1, w} ) )
\eeq
\eet

\noindent{\bf Acknowledgements.} The author is supported by an NSERC Discovery grant and a Connaught New Researcher award.

\section{Preliminaries}

\subsection{Notation}

For $ z\in \cc$ we will use the notation $z = E \pm \i \eta$, with $E \in \rr$ and $\eta >0$. 
We introduce the control parameters,
\beq
\Psi(z) := \sqrt{ \frac{\Im[ \msc(z)] }{N \eta} } + \frac{1}{N \eta}, \qquad \Psi_1 (z) := \frac{1}{ \sqrt{N \eta}}
\eeq
For an admissible test function $f$ we will also define the control parameter,
\beq
\Phi = \Phi (f) := N^{-1/2} (1+\| f''\|_1)^{1/2}, \qquad \Phi_w = \Phi_w (f) := N^{1/2} (1 + \| f''\|_{1, w} )^{1/2}.
\eeq
Note that $\Psi(z), \Psi_1(z) \geq c N^{-1/2}$ for some $c>0$. 

For a generalized Wigner matrix $H$ define the Green's function or resolvent by $G(z) := (H-z)^{-1}$ and the empirical Stieltjes transform by,
$
m_N (z) := \frac{1}{N} \sum_i G_{ii} (z).
$

For $A < B \in \rr$ we set $[\![ A, B ]\!] := \{ n \in \zz : A \leq n \leq B \}$.  The notation $A_N(i) \asymp B_N (i)$ for two $N$-dependent functions of an index $i \in \mathcal{I}$ (with $\mathcal{I}$ an abstract index set) means that there is a constant $C>0$ so that $C^{-1} A_N(i) \leq B_N (i) \leq C A_N (i)$ for all $N$ sufficiently large and all $i$. We will denote the inner product on $\cc^N$ by $\langle u, v \rangle = \sum_{i=1}^N \bar{u}_i v_i$ and the dot product by $u \cdot v = \sum_{i=1}^n u_i v_i$ for $u, v \in \cc^N$.

In order to lighten the notation, if we prove some estimate for an expression involving resolvent entries, we will often relabel the indices to just use $\{1, 2, \dots, k\}$. E.g., instead of considering $G_{i_1, i_2} \dots G_{i_{k-1} i_k}$ we consider $G_{12} G_{23} \dots G_{k-1, k}$, et cetera. Clearly this is no loss of generality, and we will not comment on this further in the paper.


\subsubsection{Stochastic domination and overwhelming probability}

We now introduce the notion of high probability events that we will use in our paper. 

\bed[Overwhelming probability] We say that an $N$-dependent family of events $\{ \Xi_u : u \in U \}$ (with $U = U^{(N)}$ an abstract index set) holds with overwhelming probability if for any large $D> 0$ we have that
\beq
\sup_{  u \in U} \pp[ \Xi_u^c ] \leq N^{-D}
\eeq
for sufficiently large $N \geq N_0 (D)$.
\eed

The following notion of stochastic domination will be useful. See \cite[Section 6.3]{erdHos2017dynamical} for standard properties of stochastic domination.
 \bed[Stochastic domination] Let $X_{u}$ and $Y_u$, for $u \in U$ be two families of nonnegative random variables and $U = U^{(N)}$ is a possibly $N$-dependent index set. We say that $X$ is stochastically dominated by $Y$ if for all $\eps >0$ and $D >0$ it holds that
 \beq
\sup_{u \in U} \pp\left[ X_u > N^{\eps} Y \right] \leq N^{-D}
 \eeq
 for all $N \geq N_0 ( \eps, D)$ sufficiently large. If $X$ is stochastically dominated by $Y$ we write $X \sd Y$. If $X$ is complex we write $X = \Osd (Y)$ if $|X| \sd Y$. Similarly, $X = Z + \Osd ( Y)$ means $|X-Z| \sd Y$.
 \eed

\subsubsection{Quasi-analytic extension}

We will have use for the following notion of a quasi-analytic extension of a test function $f$ to the complex plane.
\bed[Quasi-analytic extension] For $f : \rr \to \rr$ that is $C^2$, we define its quasi-analytic extension, $\tilf (z)$ by,
\beq
\tilf (x + \i y ) := \chi (y) ( f(x) + \i y f'(x) )
\eeq
where $\chi (y)$ is a smooth, symmetric, non-negative function satisfying $\chi (y) = 1$ for $|y| \leq 1$ and $\chi(y) = 0$ for $|y| \geq 2$. The choice of $\chi$ will be fixed throughout the paper.

We note that,
\beq \label{eqn:quasi-derivative}
\del_{\bar{z}} \tilf (z) = \frac{ \i y \chi (y) f''(x) + \i ( f(x) + i y f'(x) ) \chi'(y) }{2} .
\eeq
Here, $\del_{\bar{z}} = \frac{1}{2} (\del_x + \i \del_y) $ and $\del_{z} = \frac{1}{2} ( \del_x - \i \del_y )$ are the usual Wirtinger derivatives.
\eed

\subsection{Local laws}

We have the following local laws.
\bet
Let $\eps >0$ and fix $t_i \in \cc$ satisfying $|t_i| \leq C$ for all $i$, for some fixed $C>0$.  Uniformly for all $z$ satisfying $\Im[z] \geq N^{\eps-1}$ and $|z| \leq 10$ we have, 
\beq \label{eqn:fluct-av}
\left| \frac{1}{N} \sum_i t_i (G_{ii} (z ) - \msc (z) ) \right| \sd \frac{1}{N \eta} .
\eeq
In particular,
\beq \label{eqn:local-law}
\left| m_N (z) - \msc(z) \right| \sd \frac{1}{N \eta} .
\eeq
 For any fixed unit vectors $v, w \in \cc^N$ we have uniformly in all $z$ as above,
\beq \label{eqn:iso}
\left| v^* G (z) w - v^* w \msc (z) \right| \sd \Psi(z) .
\eeq
In particular,
\beq
\left| G_{ij} (z) - \delta_{ij} \msc (z) \right| \sd \Psi (z) \label{eqn:entry-wise}
\eeq
\eet
\proof The estimate \eqref{eqn:fluct-av} follows from, e.g., \cite[Theorem 2.6]{alt2020correlated} near the edge and \cite[Corollary 1.8]{ajanki2017universality} near the bulk. The estimate \eqref{eqn:iso} follows from \cite[Theorem 2.12]{alex2014isotropic}. \qed

\subsection{Regular matrices}

We say that a matrix $M$ is $\eps$-regular if the following hold for all $z$ satisfying $\Im[z] \geq N^{\eps-1}$ and $|z| \leq \eps^{-1}$:
\beq \label{eqn:eps-regular-def-1}
\left| \frac{1}{N} \tr (M -z)^{-1} - \msc (z) \right| \leq \frac{ N^{\eps}}{N \eta}
\eeq
and
\beq \label{eqn:eps-regular-def-2}
| (M - z )_{ij} - \delta_{ij} \msc (z) | \leq N^{\eps} \Psi(z) .
\eeq
We have the following. The proof is straightforward and deferred to Appendix \ref{a:eps-regular}.
\bel \label{lem:eps-regular}
Let $H$ be a generalized Wigner matrix for $\theta \in \cc$ and indices $a, b$ define the matrix $M^{(\theta)}$ by,
\beq
(M^{( \theta)})_{ij} = \begin{cases} \theta, & i=a,j=b \\
\bar{\theta}, & i=b, j=a\\
H_{ij}, & \mathrm{otherwise} \end{cases}
\eeq
If $a=b$ further assume $\theta \in \rr$. For any sufficiently small $\delta >0$ there is an event that holds with overwhelming probability on which every matrix $M^{(\theta)}$ with $|\theta| \leq N^{\delta-1/2}$ is $\eps$-regular. 
\eel

\subsection{Cumulant expansion}

We have the following, \cite[Lemma 3.2]{lee2018local}.

\bel \label{lem:cumu-exp}
Fix $\ell \in \nn$ and $F \in C^{ \ell+1} ( \rr; \cc)$. Let $Y$ be a centered random variable with finite $\ell+2$th moment. Then,
\beq \label{eqn:cumu-exp}
\ee[ Y F(Y) ] = \sum_{r=1}^{\ell} \frac{ \kappa^{(r+1)}(Y)}{r!} \ee[ F^{(r)} (Y)] + \ee[ \Xi]
\eeq
where the $r$th cumulant of $Y$ is denoted by $\kappa^{(r)}(Y)$, and the error term $\ee[ \Xi ]$ satisfies,
\beq
| \ee[ \Xi] | \leq C \left(\ee[ | Y|^{ \ell+2} \sup_{ |t| \leq Q } |F^{(\ell+1)} (t) | ] + \ee[ |Y|^{\ell+2} \1_{ \{ |Y|>Q \}} \sup_{t \in \rr} | F^{(\ell+1)} ( t ) | ] \right) .
\eeq
for any $Q >0$. 
\eel

We make the following convention for the  terms in the expansion on the RHS of \eqref{eqn:cumu-exp}. 
\bed [$n$th order terms] When we apply a cumulant expansion, we will refer to the term with the cumulant $\kappa^{(r)}$ of order $r$ as the ``$r$th order term(s).''

\eed

\subsection{Variance matrix and properties of semicircle law}

\bel \label{lem:spectral-gap}
The matrix of variances $S_{ij} := \frac{ s_{ij}}{N}$ has the following properties. Let $\be$ be the constant vector with entries equal to $1$. Then $S = A + N^{-1} \be \be^*$ with,
\beq
\| A \| \leq 1 - c
\eeq
for some $c>0$. 
\eel
\proof The constant vector is an eigenvector of $S$ with eigenvalue $1$. The result then follows from \cite[Lemma 7.4.2]{erdos2019matrix}. \qed

The following is proven in Appendix \ref{a:semi-calc}. 
\bel \label{lem:msc} For any $C>0$ there is a $c>0$ so that, for all $|z| \leq C$
\beq \label{eqn:msc-upper}
c \leq |\msc(z) | \leq 1
\eeq
For any $z$ and $w$ satisfying $|z| + |w| \leq C$ we have that,
\beq \label{eqn:msc-difference}
\left| \frac{ \msc(z) - \msc(w) }{ z-w} \right| \leq \frac{C_1}{ |\Im[z] | + | \Im[w] | }
\eeq
for some $C_1 >0$. Moreover,
\beq \label{eqn:msc-dif-identity}
\frac{ \msc(z) - \msc(w)}{z-w} = \frac{ \msc(z) \msc(w)}{1 - \msc(z) \msc(w)}
\eeq
\eel

\subsection{Helffer-Sjostrand formula}

For a $C^2$ function $f$ with quasi-analytic extension $\tilf (x + \i y)$ the HS formula is,
\begin{align}
 & \tr f(H) - \ee[ \tr f(H) ] = \frac{1}{ \pi} \int_{\cc^2} ( \del_{\bar{z}} \tilf (z) ) \d z \d \bar{z} \notag \\
= & \frac{1}{ 2 \pi} \int_{\rr^2} \left( \i y \chi (y) f''(x) + \i f' (x) \chi' (y) \right) N ( \ee[ m_N (x+ \i y) - \ee[ m_N ( x + \i y ) ] ) \d x \d y.
\end{align}
Define,
\beq
\Oma := \{ z : |\Im[z]|  \geq N^{\mfa-1}, | \Re[z] | \leq 6 \} .
\eeq
A straightforward argument using \eqref{eqn:local-law} shows that (see Appendix \ref{a:HS-est})
\begin{align} \label{eqn:HS-est}
 & \tr f(H) - \ee[ \tr f(H) ] \notag\\
 =&  \frac{1}{ 2 \pi} \int_{\Oma} \left( \i y \chi (y) f''(x) + \i( f(x) + \i y  f' (x) )\chi' (y) \right) N (  m_N (z) - \ee[ m_N ( z ) ] ) \d x \d y \notag\\
+ & \Osd (N^{\mfa-1} \| f''\|_1 ).
\end{align}
Note that since we assume that $\chi$ is symmetric, the integral on the RHS is real.

The following is straightforward, and the proof is deferred to Appendix \ref{a:H-est}. 
\bel \label{lem:H-est}
We have the following. Assume that $f'$ is supported in $(-5, 5)$. Suppose that $H(z)$ is a holomorphic function on $\cc \backslash \rr$. Let $S_H (z)$ denote the function,
\beq
S_H (x + \i y ) := \sup_{ |w - (x+ \i y ) | \leq |y|/2} |H(w)|.
\eeq
Then,
\begin{align} \label{eqn:general-H-est}
 & \left|  \frac{1}{ 2 \pi} \int_{\Oma} \left( \i y \chi (y) f''(x) + \i ( f(x) + \i y f' (x)) \chi' (y) \right)  H(z) \d x \d y \right| \notag\\
\leq & C(\| f \|_1+   \| f' \|_1) \left( \sup_{0.5 \leq |\Im[z] | \leq 2, |z| \leq 10 } |H(z)| \right)  \notag\\
+ & C \|f''\|_1 \left( \sup_{ |x| \leq 10 } \int_{N^{\mfa-1}  < |y| < \| f''\|_1^{-1} } |y S_H(x+ \i y) | \d y \right) \notag\\
+ & C \| f'\|_1 \left( \sup_{ |x| \leq 10 } \int_{ \| f'' \|_1^{-1}  < | y| < 2} |S_H (x + \i y ) | \d y \right)
\end{align}
In particular, if,
\beq \label{eqn:H-upper}
|H(x + \i y ) | \leq \sum_{l=1}^K \frac{ A_l}{ |y|^{s_l} }
\eeq
for some $s_l \in [1, 2]$ and $A_l >0$ we have  for admissible $f$ that,
\begin{align} \label{eqn:H-est}
 & \left|  \frac{1}{ 2 \pi} \int_{\Oma} \left( \i y \chi (y) f''(x) + \i f' (x) \chi' (y) \right)  H(z) \d x \d y \right| \notag\\
 \sd & \left(  \sum_{l=1}^K A_l +  \sum_{l=1}^k A_l \| f''\|_1^{s_l-1} \right)
\end{align}
Note also that if \eqref{eqn:H-upper} holds then,
\beq \label{eqn:del-z-H}
| \del_z^{k} H (x + \i y ) | \leq C_k \sum_{l=1}^K \frac{ A_l}{ |y|^{s_l+k} }
\eeq
for some constant $C_k >0$. 
\eel

For future reference we record here Green's theorem in complex notation which states that,
\beq
\label{eqn:green-2}
\int_{ \Omega} \del_{\bar{z}} F (z) \d z \d \bar{z} = - \frac{ \i}{2} \int_{\del \Omega} F(z) \d z
\eeq
for any $C^2$ function $F$ and sufficiently regular domain $\Omega$.

\subsection{Characteristic function}

We define,
\beq
\ea ( \lambda) := \exp \left(  \frac{\i \lambda }{ 2 \pi} \int_{\Oma} \left( \i y \chi (y) f''(x) + \i f' (x) \chi' (y) \right) N (  m_N (z) - \ee[ m_N ( z ) ] ) \d x \d y  \right)
\eeq
Due to \eqref{eqn:HS-est} we have that,
\beq
\left| \ee[ \ea ( \lambda) ] - \ee[ \exp ( \i \lambda ( \tr f(H) - \ee[ \tr f(H) ] ) ) ] \right| \prec \frac{ |\lambda| N^{\mfa} \| f''\|_1}{N} .
\eeq
Since $\chi (y)$ is symmetric, the integral in the definition of $\ea$ is real and so
$
| \ea (\lambda) | \leq 1.
$ 
We will use this estimate throughout the paper without further comment.

\subsection{Estimates for derivatives of characteristic function}

\bel \label{lem:delea}
Let $f$ be admissible. We have that for $i \neq j$,
\beq \label{eqn:deleaij}
\del_{ij} \ea = - 2 \frac{ \i \lambda \ea}{  \pi} \int_{\Oma} ( \del_{\bar{z}} \tilf (z) ) \del_z G_{ji} (z) \d z \d \bar{z} = \Osd ( |\lambda| \Phi ),
\eeq
and
\beq \label{eqn:deleaij2}
\del_{ij}^2 \ea = \Osd ( |\lambda| \Phi + |\lambda|^2 \Phi^2) + \lambda c_f \ea
\eeq
for a constant $c_f $ satisfying $c_f = \Osd (1)$. 
\eel
\proof The equality in \eqref{eqn:deleaij} follows from direct calculation using \eqref{eqn:delabgij} and the fact that
\beq
\sum_k G_{ik} (z) G_{kj} (z) = \del_z G_{ij} (z).
\eeq
 The inequality follows from the fact that $\del_{z} G_{ij} (z) = \Osd ( N^{-1/2} |\Im[z]|^{-3/2} )$ (due to \eqref{eqn:iso} and \eqref{eqn:del-z-H}) and \eqref{eqn:H-est}. For the second estimate we have by direct calculation and an argument similar to the proof of \eqref{eqn:deleaij}, 
\beq \label{eqn:delea2}
\del_{ij}^2 \ea = \Osd ( |\lambda|^2 \Phi^2)  +2 \frac{ \i \lambda \ea}{\pi} \int_{\Oma} ( \del_{\bar{z}} \tilf (z) ) \del_z ( G_{ii} G_{jj} + G_{ij}^2 ) \d \bar{z} \d z.
\eeq
We now conclude \eqref{eqn:deleaij2} using that $G_{ii} G_{jj} + G_{ij}^2 = \msc(z)^2 + \Osd ( \Psi(z) )$, and \eqref{eqn:H-est}. \qed

\bel \label{lem:del-ea-general} Let $f$ be admissible. 
Fix $n \in \nn$. Uniformly in  $i \neq j$ we have,
\beq \label{eqn:del-gen-1}
| \del_{ij}^n \ea ( \lambda ) | \sd  ((1+| \lambda |)^{\lfloor n/2 \rfloor}(1 + \Phi (1 + |\lambda| ) )  + (1+| \lambda|)^n \Phi^{n} ) \sd (1 + |\lambda| )^{ \lfloor n/2 \rfloor }
\eeq
with the second inequality holding if $|\lambda| \Phi \leq 1$. If $i= j$ then, 
\beq \label{eqn:del-gen-2}
| \del_{ii}^n \ea ( \lambda) | \sd (1+| \lambda|)^{n}.
\eeq
If $H$ is replaced by an $\eps$-regular matrix, then for any $\delta >0$, the above estimates hold with an additional factor of $N^{\delta}$ on the right hand sides, as long as $\eps>0$ is sufficiently small, depending on $\delta$ and $n$. 
\eel
\proof We begin with \eqref{eqn:del-gen-1}. The derivative $\del_{ij}^n \ea(\lambda)$ is a sum of terms of the form,
\beq\label{eqn:del-integral-1}
\ea ( \lambda) \prod_{r=1}^m \frac{ \i \lambda}{ \pi} \int_{\Oma} ( \del_{\bar{z}} \tilf (z) ) \del_{ij}^{n_r} N m_N (z) \d z \d \bar{z} ,
\eeq
where $n_i \in \nn$ s.t. $\sum_{i=1}^m n_i = n$. Similarly to the proof of Lemma \ref{lem:delea} we have that,
\beq \label{eqn:del-integral}
 \int_{\Oma} ( \del_{\bar{z}} \tilf (z) ) \del_{ij}^{n_r} N m_N (z) \d z \d \bar{z} = \Osd\left( \Phi + \1_{ \{ n_r >1 \} } \right).
\eeq
If $m \leq \lfloor n/2 \rfloor$ then using \eqref{eqn:del-integral} we can simply bound each integral in  \eqref{eqn:del-integral-1} by $\Osd (1)$ yielding an overall estimate of $\Osd ( (1 + |\lambda| )^{ \lfloor n/2 \rfloor} )$. 

On the other hand, if $m > \lfloor n/2 \rfloor$, then $|\{ r : n_r =1 \} | \geq 2m-n$, and so we can bound \eqref{eqn:del-integral-1} by $\Osd( (1+ |\lambda| )^m \Phi^{2m-n} )$. We then maximize this upper bound over choices of $m \in [\![ \lfloor n/2 \rfloor +1, n] \!]$, yielding \eqref{eqn:del-gen-1}.  

The estimate \eqref{eqn:del-gen-2} simply follows by bounding the integral on the LHS of \eqref{eqn:del-integral} by $\Osd (1)$. Finally, the statements for $\eps$-regular matrices follows simply from the fact that the above estimates depended only on estimates for the resolvent entries. \qed

\subsection{Graphical notation}

In order to develop estimates for quantities like $\ee[ \ea G_{ij} (z) ]$ we will be forced to consider more complicated quantities such as
\beq \label{eqn:example-1}
\frac{ ( \i \lambda )^2}{ \pi^2} \int_{\Oma^2} ( \del_{u_1} \del_{u_2} \ee[ G_{12} (u_1) G_{23} (u_1) G_{34} (w_1) G_{41} (u_2) ] ) \left( \prod_{i=1}^2 ( \del_{ \bar{u}_i} \tilf (u_i ) ) \d u_i \d \bar{u_i} \right) .
\eeq
In order to do so, we develop a graphical notation to assist with the evaluation of monomials in resolvent entries. Note that in \eqref{eqn:example-1}, there are two variables of integration $u_1$ and $u_2$, and a fixed variable $w_1$. In order to estimate such terms we will further apply cumulant expansions, which will result in more variables of integration $u_i$ arising. 

Due to this, we will need to distinguish between the the fixed variable $w_1$ and the two integration variables $u_1$ and $u_2$ in \eqref{eqn:example-1} (and similarly for other terms with more integration variables and fixed variables). The following definition of encoding monomials of resolvent entries by graphs is natural, except for the fact that the introduction of the two vectors $\bu$ and $\bw$ appears unmotivated at first. The motivation is due to the fact, illustrated by \eqref{eqn:example-1}, that we need to keep track of which variables are integrated over (the $\bu$) and which are fixed (the $\bw$). 

\bed [Hypergraphs]
A hypergraph $\G$ is a graph of vertices $V$ and directed edges $E$ (multiple edges and self-loops are allowed), together with the following additional data. First,  an index set $\ui \in [\![ 1, N ]\!]^{V}$. Secondly, a set of formal variables $\{ u_i\}_{i=1}^m$ for some $m \in \nn$. Finally, two vectors $\bu, \bw \in \cc^{E}$ s.t. for all $ e\in E$, exactly one of $u_e$ or $w_e$ is non-zero. Moreover, we assume that each non-zero entry of $\bu$ equals one of the $\{ u_i \}_{i=1}^m$, and furthermore that each $u_i$ appears at least once in $\bu$. 


Define $\bz = \bu + \bw$. Note that $\bz \in \cc^{E}$ and each entry contains either a variable denoted $w_e$ or $u_e$. 
Given a hypergraph $\G$, the monomial associated to $\G$ is,
\beq
\M ( \G ) := \prod_{e = (a, b) \in E } \left( G_{i_a i_b} (z_e ) - \delta_{i_a,i_b} \msc(z_e) \right) .
\eeq
Note that we can regard this, for fixed $\bw$, as a function over the variables $\{ u_i \}_{i=1}^{m}$.

We associate to a hypergraph $\G$ two control parameters,
\beq
\tilde{ \Psi} ( \G ) := \prod_{ e \in E } \Psi ( z_e), \qquad \tilde{\Psi}_1 ( \G ) := \prod_{e \in E} \Psi_1 ( z_e) .
\eeq
Note that $\M ( \G) =\Osd ( \tilde{ \Psi} ( \G))$ when we regard all of the entries $z_e$  of $\bz$ as lying in the domain  $\Oma$. 

The evaluation of a hypergraph $\G$ is,
\beq \label{eqn:T-def}
T ( \G) := \frac{ ( \i \lambda)^m}{ \pi^m} \int_{\Oma^m} \left( \del_{u_1} \del_{u_2} \dots \del_{u_m } \ee[ \ea \M ( \G ) ]  \right) \left( \prod_{i=1}^m ( \del_{ \bar{u}_i} \tilf (u_i ) ) \d u_i \d \bar{u_i} \right)
\eeq
Here, we are viewing $\ee[ \ea \M ( \G ) ]$ as a function of the $m$ variables $\{ u_i \}_{i=1}^m$ in order to define the derivatives wrt the $u_i$ and the integration over the $u_i$.

We remark also that we can regard $T( \G)$ as a function of the non-zero entries of the vector $\bw$. 

\eed
\noindent{\bf Example.} Consider the case that $(V, E)$ is a cycle on four edges and vertices. Assume that $\ui = (1, 2,3, 4)$, $\bu = (u_1, u_1, 0, u_2)$ and $\bw = (0, 0, w_1, 0)$ (so that implicitly, the formal variables are $\{ u_i\}_{i=1}^2$). Then,
\beq
\M (\G) = G_{12} (u_1) G_{23} (u_1) G_{34} (w_1) G_{41} (u_2)
\eeq
and
\beq
T ( \G ) = \frac{ ( \i \lambda )^2}{ \pi^2} \int_{\Oma^2} ( \del_{u_1} \del_{u_2} \ee[ G_{12} (u_1) G_{23} (u_1) G_{34} (w_1) G_{41} (u_2) ] ) \left( \prod_{i=1}^2 ( \del_{ \bar{u}_i} \tilf (u_i ) ) \d u_i \d \bar{u_i} \right)
\eeq
I.e., this is \eqref{eqn:example-1} above. 
\qed

\bed [Loops and lines]
A loop is a graph with vertices $V = \{ 1, \dots, K \}$ and edges $E= \{ (1, 2), (2, 3), \dots, (K-1, K), (K, 1) \}$. In the case that $K=1$, a loop has only the edges $\{ (1, 1)\}$.  A line is a graph with vertices $V= \{ 1, \dots, K \}$ (with $K \geq 2$) and edges $E= \{ (1, 2), (2, 3), \dots, (K-1, K) \}$.
\eed

For a multi-index $\ui$, we will fix the notational convention,
\beq
\{ \ui \} := \{ i_v : v \in V \} ,
\eeq
i.e., $\{ \ui \}$ is the set comprised of the entries of $\ui$. 

\subsection{Real symmetric vs complex Hermitian matrices}

We will write in detail our proofs only for the case of real symmetric generalized Wigner matrices. The proofs for the complex Hermitian case are almost the same. We give the details of how the arguments can be modified to treat the complex Hermitian case in Appendix \ref{a:cplx}. 

\subsection{Organization of remainder of paper}

In Sections \ref{sec:loops} and \ref{sec:line} we develop two general estimates for graphs involving disjoint loops with possibly one line. These estimates are the backbone for what follows. In Section \ref{sec:prelim} we develop an intermediate expansion for the $\ee[ \ea G_{ii} (z)]$ in terms of higher order terms to be estimated. These terms are of the form $\ee[ \ea G_{ij} (z) G_{ji} (w)]$ and $\ee[ \ea G_{ii} (z) G_{jj} (w) ]$. These terms are then estimated in Sections \ref{sec:g12g21} and \ref{sec:g11g22}. We finally prove our main results in Section \ref{sec:char}, along with proving various properties of the deterministic coefficients that arise. 

In Section \ref{sec:max} we extend some of the results in \cite{BLZ} to generalized Wigner matrices.

\section{Loops estimate} \label{sec:loops}

Throughout this section we will assume that $f$ is an admissible function. We will develop an estimate for the following class of hypergraphs. 

\bed \label{def:loop-admissible}
Fix $z \in \cc$. We will say that a hypergraph $\G$ is $z$-loop admissible if 
\begin{enumerate}[label={\normalfont{(\arabic*)}}]
\item Every entry of $\ui$ is unique
\item The graph $(V, E)$ underlying $\G$ is a union of disjoint loops.
\item Exactly one edge $e \in \G$ satisfies $w_e = z$, and the loop containing $e$ is of length at least $2$. 
\end{enumerate}
\eed

For any hypergraph $\G$ let $E_w := \{ e \in E : \bw_e \neq 0, \bu_e = 0 \}$. For later use, define,
\beq
 \hat{ \Psi}_1 ( \G) := \prod_{e \in E_w} \Psi_1 (z_e) .
\eeq

The main result of this section is the following estimate. The proof appears in Section \ref{sec:main-loops-proof} below.

\bep \label{prop:main-loop-estimate}
Let $\G$ be $z$-loop admissible, and let $L$ be the length of the loop containing $z$. Assume that $\bu = 0$ and that all entries of $\bw$ lie in $\Oma$, and that $| \lambda \Phi| \leq N^{-c}$ for some $c>0$. If $L \geq 3$, then,
\beq
|T (\G) | \prec \Psi_1 (z) \tilde{ \Psi}_1 ( \G)
\eeq
\eep
\remark The point of the above estimate is that the naive size of $T (\G)$, using \eqref{eqn:entry-wise}, is $\Osd ( \tilde{\Psi}_1 ( \G ) )$, and so we gain a small factor $\Psi_1 (z)$ over this estimate, due to cancellation when taking the expectation. The restriction $L \geq 3$ is optimal, as $G_{12}(z) G_{21} (\bar{z} )  = |G_{12} (z)|^2$ and so there is no improvement. \qed

\subsection{Self-consistent equation for loop graphs} \label{sec:loop-self}

This section is devoted to proving a self-consistent equation for expectations of monomials of resolvent entries, $\ee[ \ea \M ( \G)]$. We will assume throughout this section that $\G$ is a $z$-loop admissible hypergraph and that all $z_e$ satisfy $z_e \in \Oma$ (and recall that $f$ is assumed to be an admissible function). By relabelling indices we may assume that,
\beq
\M ( \G) = G_{12} (z) M_1 G_{a1} (w) M_2
\eeq
where $M_1 = G_{23} G_{34} \dots G_{a-1, a}$, and $M_2$ is a collection of loops whose indices are disjoint from the loop $G_{12} M_1 G_{a1}$. We also relabel the graph so that the vertices of the loop containing $z$ are $\{ 1, 2, \dots, a \}$; i.e., the vertex labels are the same as the index labels. If the loop containing $z$ is of length $2$, then $M_1 =1$ and $G_{a1} (w) = G_{21} (w)$. Note that it is possible that $M_2 = 1$. When the context is clear we will denote $\tilde{\Psi} = \tilde{\Psi} ( \G)$.  
Let,
\beq \label{eqn:tau-def}
\tt := \{ \ui \} \backslash \{1  \}.
\eeq

We will eventually apply a cumulant expansion to $\ee[ \ea \M ( \G ) ]$. In preparation, let us derive the following estimates. 
\bel \label{lem:del-loops-1}
For any $k \geq 1$, the quantity,
\beq \label{eqn:del-loops-1}
\del_{1j}^k ( G_{j2} (z) M_1 G_{a1} (w) M_2 ) 
\eeq
obeys the following estimates:
\begin{enumerate}[label=\normalfont(\roman*)]
\item \label{it:gen-der-1} If $j=1$, it is $ \Osd ( \tilde{ \Psi } )$. 
\item \label{it:gen-der-2} If $j=2$ and $2=a$, then it is $\Osd ( \tilde{ \Psi}  \Psi^{-1} (z) \Psi^{-1} (w) )$.
\item  \label{it:gen-der-3} If $j = a$ and $2 \neq a$, then it is $\Osd ( \tilde{\Psi} \Psi^{-1} (w) )$. 
\item \label{it:gen-der-4} If $j \in \sui$ but $j \neq 1, 2, a$ then it is $\Osd ( \tilde{\Psi} )$ as long as $j$ is not the index of a self-loop $G_{jj} - \msc$. If it is the index of a self-loop $G_{jj}(v) - \msc (v)$, then any derivative of order at least $2$ is $\Osd ( \tilde{\Psi} \Psi(v)^{-1})$, but the first derivative when $k=1$ is still $\Osd ( \tilde{ \Psi} )$.
\item \label{it:gen-der-5} If $j \notin \sui$, then it is $\Osd ( \tilde{ \Psi} )$
\end{enumerate}
These estimates also hold for any $\eps$-regular matrix with an additional factor of $N^{\delta}$ on the RHS, for any $\delta >0$ as long as $\eps >0$ is sufficiently small depending on $\delta$ and $k$. 
\eel
\proof In general, $\del_{1j}^k ( G_{ab} (v) - \delta_{ab} \msc (v) )$ is a sum of products of $k+1$ resolvent entries, and the $2(k+1)$ indices of these $k+1$ resolvent entries are comprised of one copy each of $a$ and $b$, and $k$ copies each of $1$ and $j$. Moreover, the indices $a$ and $b$ do not appear in the same resolvent entry. In particular, if  either $a \notin \{ 1, j \}$ or $b \notin \{1, j\}$, then this derivative is $\Osd ( \Psi (v) )$. From this observation, we see that the  derivatives of $M_2$ and $M_1$ can be estimated by the same product of factors of $\Psi (z_e)$ as can the original quantities $M_2$ and $M_1$, except for the following case. If there is a diagonal entry $G_{jj} (u)  - \msc (u)$ present in $M_2$, then the second and higher derivatives $\del_{1j}^k (G_{jj} (v)  - \msc (v) )$ may contain a term that is a product of diagonal resolvent entries and so is no longer $\Osd ( \Psi (v))$. These observations quickly prove items \ref{it:gen-der-4} and \ref{it:gen-der-5}. The remaining estimates of the lemma then just lists the various cases in which either the entry $G_{j2}$ or $G_{a1}$ has the indices $(1, j)$ or when $j=2$ and $G_{j2} (z) = G_{22}(z)$ is no longer $\O ( \Psi (z) )$.

The statement about $\eps$-regular matrices is simply due to the fact that all of the estimates follow from estimates on the resolvent entries, which the $\eps$-regular matrices are assumed to obey (with the fixed $N^{\eps}$ loss). \qed

In preparation for applying the cumulant expansion we derive the following estimate for the error terms.
\bel \label{lem:loops-fourth}
Let $\delta >0$ and assume $|\lambda \Phi | \leq 1$. For any $\eps$-regular matrix and sufficiently small $\eps >0$ we have that the quantity,
\beq \label{eqn:loops-fourth}
N^{-2} |\del_{1j}^3 (\ea G_{12} (z) M_1 G_{a2} (w) M_2) | 
\eeq
is bounded by $N^{\delta-1} ( 1+ |\lambda| ) \tilde{\Psi}$ if $j \in \sui$ and $N^{-2+\delta} ( 1 + |\lambda| ) \tilde{ \Psi}$ if $j \notin \sui$. 
\eel
\proof For the case $j=1$, we see that from \eqref{eqn:del-gen-2} the derivatives hitting $\ea$ contribute at most a factor of $(1+ |\lambda|)^3$ and from  Lemma \ref{lem:del-loops-1}\ref{it:gen-der-1} that the resolvent entries and their derivatives contribute a factor of $\tilde{ \Psi}$, yielding an upper bound of $N^{-2+\delta} ( 1 + |\lambda| )^3 \tilde{ \Psi}$ for this term.  Similarly, the case $j=2$ is bounded by $N^{\delta-2} ( 1 + |\lambda| ) \tilde{\Psi} \Psi(z)^{-1} \Psi(w)^{-1} \lesssim N^{\delta-1} ( 1 + |\lambda| ) \tilde{\Psi}$ by applying \eqref{eqn:del-gen-1} and  Lemma \ref{lem:del-loops-1}\ref{it:gen-der-2}. 

If $j \in \sui$ but $j \neq 1, 2$, then derivatives hitting resolvent entries can  remove at most one factor of $\Psi$ in the control parameter $\tilde{ \Psi}$, by Lemma \ref{lem:del-loops-1}\ref{it:gen-der-3},\ref{it:gen-der-4}. Combining this with \eqref{eqn:del-gen-1} we see that these terms can be bounded by $(1+ |\lambda| ) N^{-3/2+\delta} \tilde{\Psi}$. 

Finally, if $j \notin \sui$, then the resolvent entries and their derivatives always contribute a factor of $\tilde{ \Psi}$ by Lemma \ref{lem:del-loops-1}\ref{it:gen-der-5} , and so we conclude by applying again \eqref{eqn:del-gen-1}. \qed

By applying the above with the cumulant expansion, Lemma \ref{lem:cumu-exp}, we derive the following.
\bel \label{lem:loops-expand-1}
We have, for $|\lambda \Phi | \leq 1$,
\begin{align}
& z \ee[ \ea G_{12} (z) M_1 G_{a1} (w) M_2 ] = \sum_{j} \frac{s_{1j}}{N} \ee[ \del_{j1} ( \ea G_{j2} (z) M_1 G_{a1} (w) M_2 ) ] \notag\\
+ & \sum_j \frac{ s_{1j}^{(3)}}{2 N^{3/2} } \ee[ \del_{1j}^2 ( \ea G_{j2} (z) M_1 G_{a1} (w) M_2 ) ]
+ \Osd ( N^{-1}(1+ |\lambda| ) \tilde{\Psi} ) \label{eqn:loops-expand-1}
\end{align}
\eel
\proof From the fact that $(H-z) G = \1$ we have,
\beq
z \ee[ G_{12} (z) M_1 G_{a1} (w) M_2 ] = \sum_j \ee[ H_{1j} G_{j2} (z) M_1 G_{a1} (w) M_2 ] .
\eeq
The result now follows from Lemma \ref{lem:cumu-exp}, and using Lemma \ref{lem:loops-fourth} to estimate the terms at fourth order, and the fact that $|H_{ij} | \sd N^{-1/2}$. \qed

We now estimate the third order terms. 

\bel \label{lem:loops-third}
Assume $| \lambda \Phi | \leq 1$. We have for $j \in \sui$ that,
\beq \label{eqn:loops-third-1}
\left| \frac{ 1}{N^{3/2}} \ee[ \del_{1j}^2 ( \ea G_{j2} M_1 G_{a1} M_2 ) ] \right| \sd  (1+ |\lambda| ) \tilde{\Psi} N^{-1} + \delta_{a,2} \delta_{j, 2} N^{-1/2} \tilde{\Psi} 
\eeq
For other terms we have that,
\begin{align}
& N^{-3/2} \sum_{j \notin \sui } \left( \left| \ee[ ( \del_{1j} \ea ) \del_{1j} ( G_{j2} M_1 G_{a1} M_2 ) ] \right| + \left|  \ee[ ( \ea ) \del_{1j}^2 ( G_{j2} M_1 G_{a1} M_2 ) ] \right| \right) = \Osd ( \tilde{\Psi} N^{-1/2} ). \label{eqn:loops-third-2}
\end{align}
and
\beq \label{eqn:loops-third-3}
\sum_{j \notin \sui} \frac{ s_{1j}^{(3)}}{N^{3/2}} \ee[ (\del_{1j}^2 \ea ) G_{j2} M_1 G_{a1} M_2 ] = \Osd ( N^{-1/2} \tilde{ \Psi} ) .
\eeq
\eel
\proof We begin with \eqref{eqn:loops-third-1}, which will follow in a similar fashion to the first part of the proof of Lemma \ref{lem:loops-fourth}. For $j=1$, we see from Lemma \ref{lem:del-loops-1}\ref{it:gen-der-1} and \eqref{eqn:del-gen-2} that this term is $\Osd ( (1 + |\lambda| )^2 \tilde{\Psi} N^{-3/2} ) = \Osd ( N^{-1} ( 1 + |\lambda| ) \tilde{ \Psi} )$.  For $j \in \sui$ but $j \neq 1$, we consider first the case that both derivatives hit $\ea$. Then the resolvent entries contribute $\Osd ( \tilde{\Psi} \Psi(z)^{-1} )$ (the factor $\Psi(z)^{-1}$ coming in the case that $j=2$) and so this term can be bounded by $\Osd ( N^{-3/2} ( 1 + |\lambda| ) \tilde{ \Psi} \Psi(z)^{-1} ) = \Osd ( N^{-1}(1+ |\lambda| ) \tilde{ \Psi} )$, using also \eqref{eqn:del-gen-1}. Here and below we use also that $\Psi (z) \geq N^{-1/2}$. 

In the case that $j \in \sui$, $j \neq 1$ and either no or one derivatives hit $\ea$, we have  the estimate $\del_{1j} \ea = \Osd ( |\lambda| \Phi ) = \Osd (1)$ from \eqref{eqn:deleaij},  as well as the facts that derivatives hitting the resolvent entries can remove at most two factors of $\Psi(v)$ from $\tilde{\Psi}$ if $j=2$ and $a=2$ and one factor otherwise, by items \ref{it:gen-der-2}, \ref{it:gen-der-3} and \ref{it:gen-der-4} of Lemma \ref{lem:del-loops-1}.  These estimates combine to yield \eqref{eqn:loops-third-1}.

We now turn to \eqref{eqn:loops-third-2}. For both types of terms on the LHS, the resolvent entries and their derivatives are $\Osd ( \tilde{ \Psi} ) $ by  Lemma \ref{lem:del-loops-1}\ref{it:gen-der-5}. We also have $\ea = \Osd (1)$ and $\del_{1j} \ea = \Osd (1)$ by \eqref{eqn:deleaij} and so \eqref{eqn:loops-third-2} follows. 

For \eqref{eqn:loops-third-3} we have from \eqref{eqn:deleaij2} that,
\begin{align}
& \sum_{j \notin \ui} \frac{ s_{1j}^{(3)}}{N^{3/2}} \ee[ (\del_{1j}^2 \ea ) G_{j2} (z) M_1 G_{a1} (w) M_2 ] \notag\\
=  & \i \lambda c_f   \sum_{j \notin \ui} \frac{ s_{1j}^{(3)}}{N^{3/2}} \ee[ \ea G_{j2} (z) M_1 G_{a1} (w) M_2 ] + \Osd (N^{-1/2} \tilde{ \Psi} ) .
\end{align}
Now from \eqref{eqn:iso} we have that,
\beq
 \sum_{j \notin \ui} \frac{ s_{1j}^{(3)}}{N^{3/2}} G_{j2} (z) = \Osd (N^{-1} \Psi (z) )
\eeq
and so \eqref{eqn:loops-third-3} follows. \qed

We now turn to the second order terms.  We start with,
\bel \label{lem:loops-second-1}
We have,
\begin{align}
& \sum_j \frac{ s_{1j}}{N} \ee[ \ea ( \lambda)  M_1 G_{a1} (w) M_2 \del_{j1} G_{j2} (z) ] \notag \\
= &  -\msc (z)  \ee[ \ea M_1 G_{a1}(w) M_2  G_{12} (z) ] + \Osd ( \Psi(z) \tilde{\Psi} )
\end{align}
\eel
\proof By direct computation of the derivative we have,
\begin{align}
& \sum_j \frac{s_{1j}}{N} \ee[ \ea M_1 G_{a1} M_2 \del_{j1} G_{j2} (z) ]  \notag \\
= &- \sum_j \frac{s_{1j}}{N} \ee[ \ea M_1 G_{a1} M_2  G_{jj} (z) G_{12} (z) ] 
-  \sum_{j \neq 1} \frac{s_{1j}}{N} \ee[ \ea M_1 G_{a1} M_2  G_{1j}(z) G_{j2} (z) ] \notag \\
= & -\msc (z)  \ee[ \ea M_1 G_{a1} M_2  G_{12} (z) ] + \Osd ( \Psi(z) \tilde{\Psi} )
\end{align}
where in the last line we applied the entry-wise law \eqref{eqn:entry-wise}. \qed

We now define two operators on hypergraphs $\G $:
\bed \label{def:D1}
Let $\G = (V, E, \ui, \bu, \bw)$ be a hypergraph with integration variables $\{ u_i\}_{i=1}^m$. Given a vertex $v \in V$, an edge $e_1 \in E$ and  index $j$, we define a new hypergraph $D_1 ( \G, v, j, e_1) = (V', E', \ui', \bu', \bw')$ with the same integration variables $\{ u_i \}_{i=1}^m$ as follows. We add an additional vertex $v'$ to its vertex set $V' := V \cup \{ v' \}$. We add a new self-edge, $E' := E \cup \{ (v', v') \}$. The new indices $\ui'$ are $i'_l = i_l$ for $l \in V \backslash \{ v \}$, and $i'_v = j$ and $i'_{v'} = i_v$. The new spectral parameters are $\bw'_e = \bw_e$ and $\bu'_e = \bu_e$ for $e \in E$ and $\bw'_{(v', v')} = \bw_{e_1}$ and $\bu'_{(v', v')} = \bu_{e_1}$. 

\eed

\bed \label{def:D2} Let $\G = (V, E, \ui, \bu, \bw)$ be a hypergraph with integration variables $\{ u_i\}_{i=1}^m$. Given an edge $e \in E$, an index $j$ and another edge $e_1 \in E$, we define a new hypergraph $D_2 ( \G, e, j, e_1) = (V', E', \ui', \bu', \bw')$ with the same integration variables $\{ u_i \}_{i=1}^m$ as follows. Let us write $e= (a, b)$.  We add an additional vertex $v'$ to its vertex set $V' := V \cup \{ v'\}$. The new edges are the same as the old edges, except we remove $e = (a, b)$ and add edges $e_2 = (a, v')$ and $e_3 = (v', b)$. We set $i'_v = i_v$ for $v \in V$ and $i'_{v'} = j$. For all edges $e' \in E \backslash\{ e \}$ we set $\bw'_{e'} = \bw_{e'}$ and $\bu'_{e'} = \bu_{e'}$. Finally, we set $\bw'_{(a, v')} = \bw_{e_1}$, $\bu'_{(a, v')} = \bu_{e_1}$ and $\bw'_{(v', b)} = \bw_e$ , $\bu'_{(v', b)} = \bu_e$.
\eed

%
%

In our setting, the above two operators arise when $\del_{j1} $ hits $G_{a1} (w)$ in Lemma \ref{lem:loops-second-2} below. The operator $D_1$ arises from the $G_{aj}(w) G_{11} (w)$ term in $\del_{j1} G_{a1} = - G_{aj} G_{11} - G_{a1} G_{1j}$. This results in the addition of a self loop. The operator $D_2$ arises from the $G_{a1} G_{1j}$ term. This results in the loop containing $z$ increasing in length by $1$.

\bel \label{lem:loops-second-2}
We have,
\begin{align} \label{eqn:loops-second-2}
& \sum_j \frac{ s_{1j}}{N} \ee[ \ea M_1 M_2 G_{j2} (z) \del_{j1} (G_{a1} (w) ) ] \notag \\
= & - \msc (w) \sum_{j \notin \tt } \frac{s_{1j}}{N} \ee[ \ea G_{j2} (z) M_1 G_{aj} (w) M_2 ] -  \1_{ \{ a =2 \} } \frac{ s_{12} \msc(z) \msc(w)^2}{N} \ee[ \ea M_2] \notag \\
& - \sum_{ j \notin \sui } \frac{s_{1j}}{N} \ee[ \ea ( \M (D_1 (\G, 1, j, (a, 1)))  + \M ( D_2 ( \G, (1, 2), j, (a, 1) ) ) )] + \Osd ( N^{-1/2} \tilde{ \Psi} ).
\end{align}
\eel
\proof By direct computation we have,
 \begin{align} \label{eqn:loops-second-2b}
  & \sum_j \frac{s_{1j}}{N} \ee[ \ea M_1  M_2  G_{j2} (z)  \del_{j1} G_{a1} (w)] \notag \\
 = & - \sum_{j } \frac{ s_{1j}}{N} \ee[ \ea M_1 M_2 G_{j2} (z) G_{11} (w) G_{aj} (w) ]  -  \sum_{j \neq 1 } \frac{s_{1j}}{N} \ee[ M_1 M_2 G_{j2} (z) G_{a1} (w) G_{j1} (w) ] 
 \end{align}
 We now write,
 \begin{align}
 &\sum_{j } \frac{ s_{1j}}{N} \ee[ \ea M_1 M_2 G_{j2} (z) G_{11} (w) G_{aj} (w) ] 
 =  \sum_{j \neq 2 } \frac{s_{1j}}{N} \ee[ \ea M_1 M_2 G_{j2} (z) G_{11} (w) G_{aj} (w) ] \notag\\
 + &\1_{ \{ a = 2 \} }  \frac{s_{12}}{N} \ee[ \ea  M_2 G_{22} (z) G_{11} (w) G_{22} (w) ] + \Osd ( N^{-1/2} \tilde{ \Psi } )\notag \\
 = & \msc(w) \sum_{j \notin \tt } \frac{s_{1j}}{N} \ee[ \ea G_{j2} (z) M_1 G_{aj} (w) M_2 ] \notag \\
 + & \sum_{ j \notin \sui } \frac{s_{1j}}{N}  \ee[ \ea G_{j2} (z) M_1 G_{aj} (w) M_2 (G_{11} (w) - \msc (w) ) ] \notag \\
 + & \1_{ \{ a = 2 \} }  \frac{s_{12}}{N} \ee[ \ea  M_2 G_{22} (z) G_{11} (w) G_{22} (w) ] + \Osd ( N^{-1/2} \tilde{ \Psi } ) \label{eqn:loops-second-2a}
 \end{align}
 In the first estimate, we estimated the $j=2$ term by $\Osd (N^{-1} \tilde{ \Psi} \Psi(z)^{-1}  ) = \Osd ( N^{-1/2} \tilde{ \Psi} )$ if $a \neq 2$. In the second estimate we wrote $G_{11} (w) = \msc(w) + (G_{11} (w) - \msc (w))$, and also estimated the contribution of the terms $j \in \sui$ and $j \neq 1,2$  by $\Osd ( N^{-1/2} \tilde{ \Psi} )$.  The remaining terms on the second last line are of the form $\M( D_1 ( \G, 1, j, (a, 1)))$. Moreover, the terms on the third last  line of \eqref{eqn:loops-second-2a} contribute the first term on the RHS of \eqref{eqn:loops-second-2}. 
 
 When $a=2$ we can write the last term in \eqref{eqn:loops-second-2a} as,
 \beq
 \1_{ \{ a = 2 \} }  \frac{s_{12}}{N} \ee[ \ea  M_2 G_{22} (z) G_{11} (w) G_{22} (w) ]  = \1_{ \{ a =2 \} } \frac{ s_{12}}{N} \msc(z) \msc(w)^2 \ee[ \ea M_2] + \O ( N^{-1/2} \tilde{ \Psi} ) .
 \eeq
 We also have,
 \begin{align}
 & \sum_{j \neq 1 } \frac{s_{1j}}{N} \ee[ M_1 M_2 G_{j2} (z) G_{a1} (w) G_{j1} (w) ]  \notag \\
 = & \sum_{ j \notin \sui } \frac{ s_{1j}}{N} \ee[  G_{1j} (w) G_{j2} (z) M_1 G_{a1} (w) M_2 ]  + \Osd ( N^{-1/2} \tilde{\Psi} )
 \end{align}
 which yields the claim, as these terms are of the form $\M( D_2 ( \G, (1, 2), j, (a, 1)))$. \qed

 We introduce two more operations on hypergraphs $D_3$ and $D_4$.
\bed  \label{def:D3D4} Let $\G = (V, E, \ui, \bu, \bw)$ be a hypergraph with integration variables $\{ u_i\}_{i=1}^m$. Given two edges $e_1 = (a, b)$ and $e_2 = (c, d)$ in $E$, and an index $j$ we define new hypergraphs $D_3 ( \G, e_1, e_2, j)$ and $D_4 ( \G, e_1, e_2, j)$ as follows. Both have the same integration variables $\{ u_i \}_{i=1}^m$ as before. In both cases we add a new vertex $v'$ so that $V' = V \cup \{ v'\}$ and delete edges $e_1$ and $e_2$. In both cases we add edge $(v', b)$ and set $i'_{v'} = j$, and set $\bu'_{( v',b)} = \bu_{e_1}$ and $\bw'_{( v',b)} = \bw_{e_1}$.

In the case of $D_3$, we add edges $(c, v')$, $( a,d)$ and set $\bu'_{(c, v')} = \bu'_{(a,d)} = \bu_{e_2}$, $\bw'_{(c, v')} = \bw'_{(a,d)} = \bw_{e_2}$. In the case of $D_4$ we add edges $(c, a)$ and $(d, v')$ and set $\bu'_{(c, a)} = \bu'_{(d, v')} = \bu_{e_2}$, $\bw'_{(c, a)} = \bw'_{(d, v')} = \bw_{e_2}$

\eed

 \bel \label{lem:loops-second-3}
 We have,
 \begin{align} \label{eqn:loops-second-3}
 & \sum_j \frac{s_{1j}}{N} \ee[ \ea G_{j2}(z) G_{a1} (w) \del_{j1} (M_1 M_2 ) ] +  \Osd ( N^{-1/2} \tilde{\Psi} )  \notag \\
 = & - \sum_{j \notin \sui} \frac{s_{1j}}{N} \sum_{e \in E \backslash \{ (1, 2), (a, 1) \} } \ee[ \ea ( \M(D_3 ( \G, (1, 2), e, j )) +\M( D_4 ( \G, (1, 2), e, j) )) ] .
 \end{align}
 \eel
 \proof First, since no resolvent entry in either $M_1$ or $M_2$ contains the index $1$, we can estimate the terms for $j \in \sui$ by $\Osd ( N^{-1/2} \tilde{\Psi} )$. By inspection, the other derivatives generate the terms with $D_3$ and $D_4$. \qed
 
 We introduce a final operation on hypergraphs.
 \bed \label{def:D5}
 Let $\G = (V, E, \ui, \bu, \bw)$ be a hypergraph with integration variables $\{ u_i\}_{i=1}^m$. Given an edge $e \in E$, an index $j$ and a spectral parameter $u'$ we define a new hypergraph $D_5 ( \G, e, j, u') = (V', E', \ui', \bu', \bw')$ as follows. We add an additional integration variable $u_{m+1} = u'$ so that the integration variables are $\{ u_i \}_{i=1}^{m+1}$. Let us write $e= (a, b)$.  We add an additional vertex $v'$ to its vertex set $V' := V \cup \{ v'\}$. The new edges are the same as the old edges, except we remove $e = (a, b)$ and add edges $e_1 = (a, v')$ and $e_2 = (v', b)$. We set $i'_v = i_v$ for $v \in V$ and $i'_{v'} = j$. For all edges $e' \in E \backslash\{ e \}$ we set $\bw'_{e'} = \bw_{e'}$ and $\bu'_{e'} = \bu_{e'}$. Finally, we set $\bw'_{(a, v')} = 0$, $\bu'_{(a, v')} = u_{m+1}$ and $\bw'_{(v', b)} = \bw_e$ , $\bu'_{(v', b)} = \bu_e$.
 \eed
 
 \bel \label{lem:loops-second-4}
 We have for $|\lambda \Phi | \leq 1$ that,
 \begin{align} \label{eqn:loops-second-4}
 & \sum_j \frac{ s_{1j}}{N} \ee[ ( \del_{j1} \ea ) G_{j2} (z) G_{a1} (w) M_1 M_2 ] \notag \\
 = & -\sum_{j \notin \sui} \frac{ 2 \i \lambda}{ \pi} \frac{ s_{1j}}{N} \int_{\Oma} \del_{\bar{u}} \tilf (u) \del_u \ee[ \ea \M(D_5 ( \G , (1, 2) , j, u) )] \d u \d \bar{u}
 +  \Osd ( N^{-1/2} \tilde{ \Psi} ).
 \end{align}
 \eel
 \proof The $j=1$ term can be estimated by $\Osd ( N^{-1} \tilde{ \Psi} ( 1 + | \lambda| ))$ by applying \eqref{eqn:del-gen-2}. The terms for $j \in \sui$ with $j \neq 1$ can be estimated by $\Osd (N^{-1/2}  |\lambda| \Phi \tilde{ \Psi})$ using \eqref{eqn:deleaij}. Since for $j \neq 1$,
 \beq
 \del_{1j} \ea = - \frac{ 2 \i \lambda}{ \pi} \ea \int_{\Oma} \del_{\bar{u}} \tilf (u) \del_u G_{1j} (u) \d u
 \eeq
 the claim follows. \qed

Define the vector,
\beq
v_j := \ee[ G_{j2} (z) M_1 G_{aj} (w) M_2] \1_{ \{ j \notin \tt\} }
\eeq
and the hypergraph $\Gj$ to be the same as $\G$ except we replace the vertex $1$ with a vertex $j$ with $\ui_{j} = j$.  The above expansions allow us to deduce the following self-consistent equation for the vector $v_j$. Recall the definition of $\tt$ in \eqref{eqn:tau-def}. 
\bep \label{prop:loops-expand}
We have for $j \notin \tt$,
\begin{align}
 & v_j - \msc(z) \msc(w) \sum_{k} \frac{s_{jk}}{N} v_k = \1_{ \{ a = 2 \} } \frac{ \msc(z)^2 \msc (w)^2 s_{12}}{N} \ee[ \ea M_2] \notag\\
+ & \msc(z) \sum_{ k \notin \tt \cup \{ j \} } \frac{s_{jk}}{N} \ee[ \ea ( \M (D_1 (\Gj, j, k, (a, j)))  + \M ( D_2 ( \Gj, (j, 2), k, (a, j) ) ) )] \notag\\
+ & \msc(z) \sum_{ k \notin \tt \cup \{ j \} } \frac{s_{jk}}{N} \sum_{e \in E \backslash \{ (j, 2), (a, j) \} } \ee[ \ea ( \M(D_3 ( \Gj, (j, 2), e, k )) +\M( D_4 ( \Gj, (j, 2), e, k) )) ] \notag\\
+ & \msc(z) \sum_{ k \notin \tt \cup \{ j \} } \frac{ 2 \i \lambda}{ \pi} \frac{ s_{jk}}{N} \int_{\Oma} \del_{\bar{u}} \tilf (u) \del_u \ee[ \ea \M(D_5 ( \Gj , (j, 2) , k, u) )] \d u \d \bar{u}
 +  \Osd ( \Psi(z) \tilde{ \Psi} ). \label{eqn:loops-expand-2}
\end{align}
\eep
\proof For notational simplicity we do only the case $j=1$. This then follows from Lemma \ref{lem:loops-expand-1}, Lemma \ref{lem:loops-third} to estimate the third order terms on the RHS of \eqref{eqn:loops-expand-1}, and then Lemmas \ref{lem:loops-second-1}, \ref{lem:loops-second-2}, \ref{lem:loops-second-3} and \ref{lem:loops-second-4} to compute the second order terms. Note that we also use the identity $\msc(z) +z = -1 / \msc(z)$. \qed

Define the vector $Y_j$ by,
\beq
Y_j := \1_{ \{a = 2 \}} \1_{ \{ j \notin \tt \} } \frac{ \msc(z)^2 \msc(w)^2 s_{12}}{N} \ee[ \ea M_2].
\eeq
For $j \notin \tt$, define $X_j$ to be the sum of the three summations over $k$ on the last three lines of \eqref{eqn:loops-expand-2}, and $X_j = 0$ for $j \in \tt$. Let $\be$ be the vector of all ones and define the matrices $A$ and $B$ by,
\beq \label{eqn:loops-AB-def}
S = N^{-1} \be \be^T +A, \qquad B_{ij}  = S_{ij} \1_{ \{ i \in \tt \} }.
\eeq
We can rewrite \eqref{eqn:loops-expand-2} as the vector equation,
\beq \label{eqn:loops-expand-3}
(1 - \msc(z) \msc(w) (A-B) ) v = \msc(z) \msc(w) N^{-1} \be \be^T v + X + Y + \eps
\eeq
where $\eps$ is a vector satisfying $\| \eps \|_\infty = \Osd( \Psi(z) \tilde{ \Psi } )$. Note that the case when $j \in \tt$ is trivial because $[(S-B)v]_j=0$ for such $j$.  With the above equation as input we can prove the following.
\bel \label{lem:loops-expand-2}
For any $j \notin \tt$ we have that,
\begin{align} 
v_j &= \sum_{k} C_{jk} (Y_k + X_k +Z_k) + \Osd ( \Psi(z) \tilde{ \Psi} ) \notag\\
 +& \Osd \left(  \frac{ \Psi(z) + \Psi(w)}{ N ( |\Im[z] | + | \Im[w] | ) } \tilde{ \Psi} \Psi(z)^{-1} \Psi(w)^{-1} \right). \label{eqn:loops-expand-4}
\end{align}
where $C_{jk}$ are coefficients satisfying $|C_{jk} | \leq C(\delta_{jk} + N^{-1} )$, and $Z_k$ is the constant vector with entries,
\beq
Z_k = \msc(z) \msc(w) \frac{\delta_{a2}}{N} \ee[ \ea M_2] \left( \frac{ \msc(z) - \msc(w) }{z-w} - \msc(w) \msc(w) \right).
\eeq
Moreover, the $C_{jk}$ are holomorphic in the variables $z_e$ as they vary over $\cc \backslash \rr$. 
\eel
\proof We invert the vector equation \eqref{eqn:loops-expand-3} by moving $(1 - \msc(z) \msc(w) (A-B)$ to the RHS. It follows that the coefficients will be $C_{jk} := (1 - \msc(z) \msc(w) (A-B) )^{-1}_{jk}$ which satisfy the desired estimates due to Lemma \ref{lem:inverse} (the hypotheses of which are satisfied due to Lemma \ref{lem:spectral-gap} and \eqref{eqn:msc-upper}). We compute,
\begin{align}
& \frac{1}{N} \be^T v = \frac{1}{N} \sum_{ j \notin \tt } \ee[ G_{aj} (w) G_{j2} (z) M_1 M_2 \ea ] \notag\\
= & \frac{1}{N} \sum_{j}  \ee[ G_{aj} (w) G_{j2} (z) M_1 M_2 \ea ] - \1_{ \{ a =2 \}} \frac{ \msc(z) \msc(w)}{N} \ee[ M_2 \ea ] + \Osd ( N^{-1/2} \tilde{ \Psi} ) \notag\\
= & \frac{\delta_{a2}}{N} \ee[ \ea M_2] \left( \frac{ \msc(z) - \msc(w) }{z-w} - \msc(w) \msc(w) \right) + \Osd ( N^{-1/2} \tilde{ \Psi} ) \notag\\
+& \Osd \left(  \frac{ \Psi(z) + \Psi(w)}{ N ( |\Im[z] | + | \Im[w] | ) } \tilde{ \Psi} \Psi(z)^{-1} \Psi(w)^{-1} \right).
\end{align}
where in the last line we applied Lemma \ref{lem:G-dif}. In the second line we used \eqref{eqn:entry-wise} in the case that $a=2$ to replace the $j=2$ term with the second term on the RHS.  The claim now follows from this and \eqref{eqn:loops-expand-3}. \qed

\subsection{Iterative expansion} \label{sec:loop-iterative}

The goal of this section is to prove Proposition \ref{prop:main-loop-estimate} via an iterative expansion. Our iterative expansion will repeatedly use Lemma \ref{lem:loops-expand-2} to estimate $T( \G)$ in terms of contributions from a family of hypergraphs.  We now define this family in Definition \ref{def:loops-family}. We first require some notation.

Let $\G$ be a $z$-loop admissible hypergraph.  Denote the edge whose spectral parameter is $z$ by $e_z = (v_1, v_2)$ and assume $i_{v_1} = a$. Let $e_1 = (v_3, v_1)$ be the edge leading to $e_z$.  Let $\Gj$ be the same hypergraph as $\G$ but with $i_v = j$. Let $\tt := \{ \ui \} \backslash \{ a \}$. 

\bed \label{def:loops-family} For a $z$-loop admissible hypergraph $\G$ as above, we let $\hat{\G}_z$ be the hypergraph obtained by removing the loop containing $z$ from $\G$ and let $L_z ( \G)$ be the length of the loop containing $z$. We furthermore define a set of hypergraphs $F_1 ( \G)$ as follows. It is the union of the following sets. 
\begin{enumerate}[label=\normalfont(\roman*)]
\item The hypergraphs $D_1 ( \Gj, v_1, k, e_1 )$ for all $j, k$ such that $j \notin \tt$ and $k \notin \tt \cup \{ j \}$. 
\item The hypergraphs $D_2 (\Gj, e_z, k, e_1 )$ for all $j, k$ such that $j \notin \tt $ and $k \notin \tt \cup \{ j \}$. 
\item The hypergraphs $D_3 ( \Gj, e_z, e, k)$ and $D_4 ( \Gj, e_z, e, k)$ where $e$ varies over $e \in E \backslash \{ e_z, e_1\}$ and the indices $j, k$ vary over $j \notin \tt$ and $k \notin \tt \cup \{ j \}$.
\item The hypergraphs $D_5 ( \Gj, e_z, k, u)$ where $j \notin \tt$ and $k \notin \tt \cup \{ j \}$.
\end{enumerate}
\eed

Let now $\G_0$ be a $z$-loop admissible hypergraph s.t. $L_z ( \G_0) \geq 3$ and $\bu = 0$. Assume that all entries of $\bw$ line in $\Oma$ for some $\mfa >0$. We define the following sets inductively. Let $\F_0 := \{ \G_0 \}$ and for $n \geq 1$ set,
\beq
\F_n := \bigcup_{ \Hcal \in \F_{n-1} } F_1 ( \Hcal ).
\eeq
By inspection, one finds the following:
\bel \label{lem:iterative-properties} \begin{enumerate}[label=\normalfont(\roman*)]
\item  For any $n \geq 1$, all hypergraphs in $\F_n$ contain $n$ more edges than $\G_0$.
\item Every graph in any $\F_n$ is $z$-loop admissible. 
\item For every $\Hcal \in \F_n$ we have $\hat{ \Psi}_1 ( \Hcal ) \leq \tilde{ \Psi}_1 ( \G_0 ) = \hat{\Psi}_1 ( \G_0)$. 
\item If $\Hcal_1 \in F_1 ( \Hcal)$ for some $\Hcal \in \F_n$, then
\beq
\tilde{\Psi}_1 ( \Hcal_1) = \Psi_1 ( v) \tilde{ \Psi}_1 ( \Hcal)
\eeq
for some $v$ that equals some spectral parameter of $\Hcal$. 
\end{enumerate}
\eel

The next lemma translates Lemma \ref{lem:loops-expand-2} into a form useable for our iteration argument.
\bel \label{lem:loops-expand-3}
Let $n \geq 0$. Uniformly in $\Hcal \in \F_n$ we have,
\begin{align}
T ( \Hcal  ) &= \1_{ \{ L_z ( \Hcal ) =2 \} } a_1 T (  \hat{ \Hcal}_z  )  + \sum_{ \Hcal_1 \in F_1 ( \Hcal ) } a_{\Hcal_1} T ( \Hcal_1 )\notag  \\
&+ \Osd ( \Psi_1 (z)  \hat{\Psi}_1 ( \Hcal ) ) \label{eqn:loops-expand-5}
\end{align} 
for some coefficients obeying 
\beq
\sum_{ \Hcal_1 \in F_1 ( \Hcal ) } \left| a_{ \Hcal_1 } \right| \leq C_n, \qquad |a_1| \leq \frac{C_n}{N | \Im[z] | } .
\eeq
\eel
\proof This follows almost directly from applying Lemma \ref{lem:loops-expand-2} to the expectation in the integrand of the definition of $T ( \Hcal )$ in \eqref{eqn:T-def}. In particular, the terms with $Y_k$ and $Z_k$ on the RHS of \eqref{eqn:loops-expand-4} correspond to the first term on the RHS of \eqref{eqn:loops-expand-5}, and $X_k$ terms correspond to to the sum on the RHS. The only thing to verify is that the error term on the RHS of \eqref{eqn:loops-expand-4} can be integrated in the $u_i$ to obtain the error $\O ( \Psi_1 (z) \hat{ \Psi}_1 (  \Hcal ) )$.  First, note that the error terms on the RHS of \eqref{eqn:loops-expand-4} are $\Osd ( \Psi_1 (z) \tilde{ \Psi}_1 ( \Hcal ) )$ and are holomorphic functions in the spectral parameters of $\Hcal$ (as the error term is a difference of holomorphic functions). Let us assume that the hypergraph $\Hcal$ contains $K$ independent integration variables $u_1, \dots, u_K$ (note that any single $u_i$ may be repeated amongst the non-zero entries of $\bu$).   We have that,
\beq
\tilde{ \Psi}_1 ( \Hcal) \leq \hat{ \Psi}_1 ( \Hcal ) \prod_{j=1}^K \Psi_1 ( u_j ).
\eeq
Therefore, the error on the RHS of \eqref{eqn:loops-expand-4} integrates to at most,
\beq
\Osd \left( ( | \lambda | \Phi)^K \Psi_1 (z) \hat{ \Psi}_1 ( \Hcal ) \right)
\eeq
by repeatedly applying \eqref{eqn:H-est}. \qed

We single out the case where $L_z ( \Hcal) =2 $ for special treatment, as a naive estimate of the size of these terms would be too large.

\bel \label{lem:loop-2}
Let $\Hcal \in \F_n$ for some $n \geq 0$. Assume that $L_z ( \Hcal ) = 2$, and let $\J =  \hat{\Hcal}_z$. Then, with $a_1$ as in \eqref{eqn:loops-expand-5},
\beq
| a_1 T ( \J ) | = \Osd( \Psi_1 (z) \tilde{ \Psi}_1 ( \G_0 ) ).
\eeq
\eel
\proof First note that $n \geq 1 $ because by assumption $L_z (\G_0) \geq 3$. Let $\hat{z}$ be the spectral parameter of the other edge in this loop. It either equals $w_i$ for some $w_i$ in the vector $\bw$ of $\Hcal$ or $u_i$ for some integration variable. In either case, we claim that there is an edge $e \in \J$ s.t. $z_e = \hat{z}$.

Consider the sequence of parents of $\Hcal$, i.e., the sequence of graphs $\hat{J}_i$ s.t. $\hat{J}_0 = \G$, $\hat{J}_i \in F_{1} ( \hat{J}_{i-1} )$  and $\hat{J}_n = \Hcal$. Each of these graphs arise from one of the actions $D_1, D_2, D_3, D_4$ or $D_5$ as described above.  Let $k$ be the smallest integer s.t. the graphs $\{ \hat{J}_l \}_{k \leq l \leq n }$ all have the property that the loop containing $z$ is of size $2$ and the other spectral parameter is $\hat{z}$. Then $1 \leq k \leq n$. Consider how $\hat{J}_k$ was constructed from $\hat{J}_{k-1}$. Since actions $D_1$, $D_2$ and $D_5$ either leave the $z$-loop unchanged or increase its length by $1$, it cannot have arisen from either of  these operations. So it must have come from either $D_3$ or $D_4$. But this action creates an additional edge elsewhere in the graph with the same spectral parameter of the edge that gets connected to the edge containing $z$. Therefore our claim that there is an edge $e \in \J$ with $z_e = \hat{z}$ is proven.

Now, using this and Lemma \ref{lem:iterative-properties} we see that whether or not $\hat{z}$ is an integration variable, the parameter $\tilde{ \Psi}_1 ( \J )$ has at least one copy of $\Psi_1 (u_i)$ for each integration variable and at least one copy of $\Psi (v)$ for all $v$ that are not equal to $z$ in the original graph $\G_0$ (recall that we assume that $\G_0$ has no integration variables). Therefore, if there are $K$ independent variables of integration in $\Hcal$ we have that,
\beq
| a_1 T ( \J ) | = \Osd \left( \frac{1}{N |\Im[z]| } ( | \lambda| \Phi )^K \tilde{ \Psi}_1 ( \G ) \Psi_1(z)^{-1} \right) = \Osd( \Psi_1(z) \tilde{ \Psi}_1 ( \G ) ).
\eeq
This completes the proof. \qed

The following allows us to truncate our iterative expansion at large $n$. 
\bel \label{lem:loops-expand-4}
Let $\Hcal \in \F_n$. Let $\min_{ z \in \G} N \Im[z] \geq N^{4\alpha}$ and $| \lambda| \Phi \leq N^{-4\alpha}$ for some $\alpha >0$ satisfying $\alpha < \mfa /10$. Then,
\beq
| T ( \Hcal ) | \leq  N^{ - n \alpha }.
\eeq
for all $N$ sufficiently large, depending only on $n$.
\eel
\proof Suppose that there are $K$ independent variables of integration in $\Hcal$. Note that $n \geq K$. Then we have,
\beq
\M( \Hcal ) = \Osd ( N^{ - (n-K) 2\alpha}  \prod_{i=1}^K \Psi_1 (u_i ) )
\eeq
and so the claim follows after repeatedly applying \eqref{eqn:H-est}. \qed

\subsubsection{Proof of Proposition \ref{prop:main-loop-estimate}} \label{sec:main-loops-proof}

For any $n$, we have by iterating \eqref{eqn:loops-expand-5} and applying Lemma \ref{lem:loop-2} that,
\begin{align}
| T ( \G) | \sd \Psi_1 (z) \tilde{\Psi}_1 ( \G ) + \max_{ \Hcal \in \F_n } | T ( \Hcal ) |.
\end{align}
For large enough $n$, the second term on the RHS is less than $\Psi_1(z) \tilde{ \Psi}_1 ( \G )$ by Lemma \ref{lem:loops-expand-4} and so the claim follows. \qed

\subsection{Better estimate in $L=2$ case} \label{sec:short-loops}

The estimate of Proposition \ref{prop:main-loop-estimate} applied in the case when the loop containing the spectral parameter $z$ was of length at least $3$. In this section we establish an estimate when the loop is of length $2$ and the monomial $M_2$ is of a specific form. 

Specifically, we consider terms of the type, $\ee[ \ea G_{12} (z) G_{21} (w) X_n ]$,
where
\beq
X_n := \prod_{k=3}^{2+n}  (G_{kk} (v_k) - \msc (v_k ) ).
\eeq
and $n \geq 1$ with $X_0 = 1$.  Define now,
\beq
\tilde{\Psi} := \Psi(z) \Psi(w) \prod_{k=3}^{2+m} \Psi(v_k), \qquad \tilde{\Psi}_1 := \Psi_1(z) \Psi_1(w) \prod_{k=3}^{2+m} \Psi_1(v_k),
\eeq
We let $\sui$ be the set of indices appearing in $X_n$ and $G_{12} (z) G_{21} (w)$ and once again let $\tt := \sui \backslash \{ 1 \}$. 
\bel \label{lem:weakgij-expand}
We have for $|\lambda \Phi | \leq 1$,
\begin{align}
& z \ee[ \ea G_{12} (z) G_{21} (w) X_n ] = \sum_{j} \frac{ s_{1j}}{N} \ee[ \del_{j1} (G_{j2} (z) G_{21} (w) X_n \ea ) ] 
+  \Osd ( N^{-1/2} \tilde{ \Psi} )
\end{align}
\eel
\proof This follows immediately from \eqref{eqn:loops-expand-1} and Lemma \ref{lem:loops-third}. \qed

The following estimates the second order terms in the above expansion.
\bel \label{lem:weakgij-second}
We have for $j  \in \sui$ and $|\lambda \Phi | \leq 1$,
\beq \label{eqn:weakgij-second-1} 
N^{-1}  \ee[ \del_{j1} ( \ea G_{j2} (z) G_{21}(w)  X_n ) ] = - \frac{ \msc(z) \msc(w)^2}{N} \ee[ X_n] \delta_{j=2} +\Osd ( N^{-1/2} \tilde{ \Psi } )
\eeq
We have for $j \notin \sui$
\beq \label{eqn:weakgij-second-2} 
| \ee[ G_{j2} (z) G_{21} (w) X_n \del_{j1} \ea  ] | \sd |\lambda \Phi | \tilde{\Psi}_1 \min\{ \Psi_1(z), \Psi_1(w), \Phi \}
\eeq
and
\beq \label{eqn:weakgij-second-3}
| \ee[ \ea G_{j2} (z) G_{21} (w) \del_{j1} X_n] | \sd \tilde{ \Psi}_1 \min \{ \Psi_1 (w),  \Psi_1 (z) \} .
\eeq
Finally,
\begin{align}\label{eqn:weakgij-second-4}
& \sum_{j \notin \sui} \frac{ s_{1j}}{N} \ee[ \ea X_n \del_{j1} (G_{j2} (z) G_{21} (w) ) ] = - \msc(z) \ee[ G_{12} (z) G_{21} (w) X_n \ea ]  \notag\\
&- \msc(w) \sum_{j \notin \tt} \frac{ s_{1j}}{N} \ee[ G_{j2} (z) G_{2j} (w) X_n \ea ]  - \sum_{j \notin \ui} \frac{s_{1j}}{N} \ee[ \ea G_{j2} (z) G_{2j} (w) X_n (G_{11} (w) - \msc (w) ) \ea ] \notag\\
& - \sum_{j \notin \sui}  \frac{s_{1j}}{N} \ee[ G_{12} (z) G_{21} (w) (G_{jj} (z) - \msc (z) ) X_n \ea ] + \Osd ( \tilde{\Psi}_1 \min\{ \Psi_1 (z) ,\Psi_1 (w) \} )
\end{align}
\eel
\proof The estimate \eqref{eqn:weakgij-second-1} follows easily using \eqref{eqn:del-gen-1}, \eqref{eqn:del-gen-2} and \eqref{eqn:entry-wise} and direct computation. For \eqref{eqn:weakgij-second-2} we have
\beq
\ee[ G_{j2} (z) G_{21} (w) X_n \del_{1j} \ea ] = \frac{ -2 \i \lambda}{\pi} \int_{ \Oma} ( \del_{\bar{u}}  \tilf (u ) ) \del_u \ee[ G_{j1} (u) G_{j2} (z) G_{21} (w) X_3 ] \d u 
\eeq
The expectation in the integral is $\Osd ( \tilde{ \Psi}_1 \Psi_1(u) \min \{ \Psi_1(z), \Psi_1 (w), \Psi_1(u) \} )$ by Proposition \ref{prop:main-loop-estimate}, and so \eqref{eqn:weakgij-second-2}  follows by applying \eqref{eqn:H-est}. 

Since $\del_{j1} (G_{kk} (v_k ) - \msc (v_k ) ) = -2 G_{k1} (v_k) G_{kj} (v_k)$, the estimate \eqref{eqn:weakgij-second-3} follows from Proposition \ref{prop:main-loop-estimate}. The estimate \eqref{eqn:weakgij-second-4} follows in a similar manner, by applying Proposition \ref{prop:main-loop-estimate}. Note that we added back the term with $j=1$ in the first sum on the RHS of \eqref{eqn:weakgij-second-4} at a cost of $\Osd (N^{-1} \tilde{ \Psi}_1 )$.  \qed

Recalling the definition of $\tt$ above, define
\beq
v_j := \ee[ G_{j2} (z) G_{2j} (w) \ea X_n] \1_{ \{ j \notin \tt \} } .
\eeq
Then, Lemmas \ref{lem:weakgij-expand} and \ref{lem:weakgij-second} imply, (recall $A$ and $B$ defined in \eqref{eqn:loops-AB-def}),
\beq \label{eqn:weakgij-self-consist}
(1 - \msc (z) \msc (w) (A-B) )v = \msc(z) \msc(w) N^{-1} \be \be^T v + U + V + \eps
\eeq
where $\| \eps \|_\infty \leq \tilde{\Psi}_1 \min \{ \Psi_1 (z), \Psi_1 (w) \}$ and
\beq
U_j = \1_{ \{ j \notin \tt \}} \frac{ \msc(z)^2 \msc(w)^2 s_{j2}}{N} \ee[ \ea X_n]
\eeq
and
\begin{align}
V_j &= \1_{ \{ j \notin \tt\} } \msc(z) \sum_{k \notin \tt \cup\{ j\} } \frac{ s_{jk}}{N} \ee[ \ea (G_{k2}(z) G_{2k} (w) X_n (G_{kk} (w) - \msc (w) ) ]  \notag\\
&+\1_{ \{ j \notin \tt\} } \msc(z) \sum_{k \notin \tt \cup \{ j \}} \frac{ s_{jk}}{N} \ee[ \ea (G_{j2}(z) G_{2j} (w) X_n (G_{kk} (z) - \msc (z) ) ]
\end{align}

We therefore may obtain the following.
\bel
For $|\lambda \Phi | \leq 1$ we have,
\begin{align} \label{eqn:weakgij-iterative}
\| v \|_\infty & \sd \left( \left| \frac{ \msc(z) - \msc(w) }{ N (z-w) } \right| + N^{-1} \right) | \ee[ \ea X_n] |  
 + \| V \|_\infty + \tilde{\Psi}_1 \min \{ \Psi_1 (z), \Psi_1 (w) \}.
\end{align} 
\eel
\proof  By applying Lemma \ref{lem:G-dif} (similarly to the proof of Lemma \ref{lem:loops-expand-2})  one sees that
\beq
\left| N^{-1} \be^T v \right| \sd  \left( \left| \frac{ \msc(z) - \msc(w) }{ N (z-w) } \right| + N^{-1} \right) | \ee[ \ea X_n] | + \min\{ \Psi_1(z),\Psi_1 (w) \} \tilde{ \Psi}_1 .
\eeq
The claim now follows from \eqref{eqn:weakgij-self-consist} and Lemma \ref{lem:inverse} (again, similar to the proof of Lemma \ref{lem:loops-expand-2}). \qed

Iterating the above, we can quickly conclude the following.
\bel \label{lem:weak-g12g21g33}
We have for $| \lambda \Phi | \leq 1$, and $u, v, w \in \Oma$,
\begin{align} \label{eqn:weak-g12g21g33}
 & \left| \ee[ \ea (G_{12} (z) G_{21} (w) ) (G_{33} (v) - \msc (v) ) ] \right| \notag\\
 \leq&  \min\{ \Psi^2_1 (z), \Psi^2_1 ( w) \}  \max_{j \in [\![1, N]\!]} | \ee[ \ea (G_{jj} (v) - \msc (v) ) ] | +  \Psi_1 (v) \Psi_1 (w) \Psi_1(z) \min \{ \Psi_1 (z), \Psi_1 (w) \} 
\end{align} 
\eel
\proof Consider a vector $v$ consisting of expectations of the form $\ee[ \ea G_{j2} (z) G_{2j} (w) X_n ]$ where $X_n$ is a product of $n$ terms of the form $G_{kk} (v_k ) - \msc (v_k)$. Then the vector $V$ on the RHS of \eqref{eqn:weakgij-iterative} is also of this form, except it has an additional term of the form $G_{kk} (v ) - \msc (v)$ where $v$ is either $z$ or $w$. We therefore can iterate \eqref{eqn:weakgij-iterative} repeatedly to the term $\| V \|_\infty$ that appears on the RHS. The first iteration starts with the term on the LHS of \eqref{eqn:weak-g12g21g33}. For this first iteration, the first term on the RHS of \eqref{eqn:weakgij-iterative} is bounded by the first term on the RHS of \eqref{eqn:weak-g12g21g33}. For the subsequent iterations, the first term on the RHS of \eqref{eqn:weakgij-iterative} can be absorbed into the last term on the RHS of \eqref{eqn:weak-g12g21g33}. At all steps of the iteration the last error term on the RHS of \eqref{eqn:weakgij-iterative} can be absorbed into the last term on the RHS of \eqref{eqn:weak-g12g21g33}. The iteration may be terminated in a fashion similar to the proof of Proposition \ref{prop:main-loop-estimate}. \qed


\section{Line estimate} \label{sec:line}

In Section \ref{sec:loops} we obtained estimates for expectations of resolvent entries whose associated graphs was a product of loops. In this section we obtain an estimate for hypergraphs containing a line. We again assume throughout that $f$ is an admissible function. The overall strategy is roughly similar in that we will iterate the cumulant expansion and inspect which terms arise.

\bed \label{def:line-admissible}
We will say that a hypergraph is line admissible if 
\begin{enumerate}[label=\normalfont(\roman*)]
\item Every entry of $\ui$ is unique
\item It is a disjoint union of loops and exactly one line.
\end{enumerate}
\eed
We will assume $\M ( \G)  = G_{12} (z) M_1 M_2$ where $M_1 = G_{23} G_{34} \dots G_{K-1, K}$ and $M_2$ is a disjoint union of loops. The parameter $z$ will not play a particularly special role here, unlike the argument in Section \ref{sec:loops}. Our goal is to prove the following. The proof appears in Section \ref{sec:main-line-proof}. 
\bep \label{prop:main-line-estimate}
Let $\G$ be a line admissible hypergraph with $\bu = 0$ and all $z_e \in \Oma$, and assume $|\lambda \Phi| \leq N^{-c}$ for some $c>0$. Then we have,
\beq
| \ee[ \ea \M ( \G) ] | \sd N^{-1/2} \tilde{ \Psi}.
\eeq
\eep

We first develop some estimates in preparation for applying the cumulant expansion.
\bel \label{lem:line-der}
The term,
\beq \label{eqn:line-der}
\del_{j1}^k (G_{j2} (z) M_1 M_2 )
\eeq
obeys the following estimates:
\begin{enumerate}[label=\normalfont(\roman*)]
\item If $j=1$ it is $\Osd ( \tilde{ \Psi} )$. \label{it:line-der-1}
\item If $j=2$ it is $\Osd ( \tilde{ \Psi} \Psi(z)^{-1} )$. \label{it:line-der-2}
\item If $ j \in \sui$ but $j\neq 1, 2$, then it is $\Osd ( \tilde{ \Psi} )$ as long as $j$ is not the index of a self-loop $(G_{jj}(v) - \msc (v))$. If $j$ is the index of a self-loop $(G_{jj} (v) - \msc (v))$ then any derivative of order at least $2$ is $\Osd ( \tilde{ \Psi} \Psi(v)^{-1} )$ but the first derivative is still $\Osd ( \tilde{ \Psi} )$ \label{it:line-der-3}
\item If $j \notin \sui$, then it is $\Osd ( \tilde{ \Psi} )$ \label{it:line-der-4}
\end{enumerate}
If $\delta >0$ and $M$ is $\eps$-regular, than the above esimates hold with an $N^{\delta}$ factor as long as $\eps$ is sufficiently small.
\eel
\proof The proof is almost identical to Lemma \ref{lem:del-loops-1} and so is omitted. \qed
\bel \label{lem:line-fourth}
We have that for $| \lambda \Phi| \leq 1$ that,
\beq \label{eqn:line-fourth}
N^{-2} | \del_{j1}^{3} ( \ea G_{j2} M_1 M_2 \ea ) |
\eeq
is $\Osd ( (1 + |\lambda| )^2 \tilde{ \Psi} / N^{3/2} )$ for $j \in \sui$ and $\Osd (( 1+ |\lambda| ) N^{-2} \tilde{ \Psi} )$ otherwise. The same estimates hold with an additional factor of $N^{\delta}$ if $M$ is $\eps$-regular, for $\eps >0$ sufficiently small. 
\eel
\proof When $j=1$, the derivatives hitting $\ea$ contribute at most $(1 + |\lambda| )^3$ by \eqref{eqn:del-gen-2}, and the derivatives hitting resolvent entries can be estimated using  Lemma \ref{lem:line-der}\ref{it:line-der-1}, finding $\Osd ( (1 + |\lambda| )^3 \tilde{ \Psi} N^{-2} )$. For $j=2$, we use \eqref{eqn:del-gen-1} and  Lemma \ref{lem:line-der}\ref{it:line-der-2} to find an estimate of $\Osd ( (1 + | \lambda| ) \tilde{ \Psi} / N^{3/2} )$. For $j \in \sui$ and $j \neq 1, 2$ we use \eqref{eqn:del-gen-1} and  Lemma \ref{lem:line-der}\ref{it:line-der-3}  to find an estimate of $\Osd ( (1+ |\lambda| ) \tilde{ \Psi} / N^{3/2} )$. 

Finally, the terms when $j \notin \sui$ are estimated using \eqref{eqn:del-gen-1} and   Lemma \ref{lem:line-der}\ref{it:line-der-4}. \qed

From this and the cumulant expansion, Lemma \ref{lem:cumu-exp} we find the following. The proof is the same as Lemma \ref{lem:loops-expand-1} and so omitted.
\bel \label{lem:line-expand-1}
We have for $| \lambda \Phi | \leq 1$,
\begin{align}
& z \ee[ \ea  G_{12} (z) M_1 M_2 ] = \sum_j\frac{s_{1j}}{N} \ee[ \del_{j1} (G_{j2} (z) M_1 M_2 \ea ) ] \notag \\
+ & \sum_j \frac{ s_{1j}^{(3)}}{2 N^{3/2}} \ee[ \del_{1j}^2 (G_{j2} (z) M_1 M_2 \ea ) ] + \Osd ( (1+ |\lambda| ) N^{-1} \tilde{ \Psi} ). \label{eqn:line-expansion-1}
\end{align}
\eel

We now begin estimating the third order terms on the RHS of \eqref{eqn:line-expansion-1}. 

\bel \label{lem:line-third}
The following hold for $| \lambda \Phi | \leq 1$. We have for $j \in \sui$ that,
\begin{align} \label{eqn:line-third-1}
\left| N^{-3/2} \ee[ \del_{j1}^2 ( \ea G_{j2} M_1 M_2 ) ] \right| \sd (1 + |\lambda| ) \tilde{ \Psi} N^{-1}.
\end{align}
For the other terms we have,
\begin{align}
& \sum_{ j \notin \sui } \frac{1}{N^{3/2}} | \ee[ ( \del_{j1} \ea ) \del_{j1} (G_{j2} M_1 M_2 ) ] | \sd |\lambda \Phi | N^{-1/2} \tilde{ \Psi} \label{eqn:line-third-2} \\
& \sum_{j \notin \sui} \frac{1}{N^{3/2}} | \ee[ \ea \del_{j1}^2 ( G_{j2} M_1 M_2 ) ] | \sd N^{-1/2} \tilde{ \Psi} \label{eqn:line-third-3} \\
& \sum_{j \notin \sui } \frac{ s_{1j}^{(3)}}{N^{3/2}} \ee[ ( \del_{j1}^2 \ea ) G_{j2} M_1 M_2 ] = \Osd( | \lambda | \Phi N^{-1/2} \tilde{ \Psi} ) \label{eqn:line-third-4}
\end{align}
\eel
\proof We start by proving \eqref{eqn:line-third-1}. For $j=1$ we use \eqref{eqn:del-gen-2} and   Lemma \ref{lem:line-der}\ref{it:line-der-1} to find an estimate of $\Osd( (1 + | \lambda| )^2 N^{-3/2} \tilde{ \Psi} )$. When $j=2$ we use \eqref{eqn:del-gen-1} and  Lemma \ref{lem:line-der}\ref{it:line-der-2} to find an estimate of $\Osd ( (1 + | \lambda| ) \tilde{ \Psi} N^{-1} )$. For other $j \in \sui$ we use \eqref{eqn:del-gen-1} and  Lemma \ref{lem:line-der}\ref{it:line-der-3} to find an estimate of $\Osd ( (1 + |\lambda| ) \tilde{ \Psi} N^{-1} )$. This completes the proof of \eqref{eqn:line-third-1}.

The estimates  \eqref{eqn:line-third-2} and \eqref{eqn:line-third-3} follow immediately from \eqref{eqn:deleaij2} and  Lemma \ref{lem:line-der}\ref{it:line-der-4}. For \eqref{eqn:line-third-4} we use \eqref{eqn:deleaij2} to conclude,
\begin{align}
& \sum_{j \notin \sui } \frac{ s_{1j}^{(3)}}{N^{3/2}} \ee[ ( \del_{1j}^2 \ea ) G_{j2} M_1 M_2 ] 
= \Osd( \tilde{ \Psi} |\lambda \Phi | N^{-1/2} ) + \i \lambda c_f \sum_{j \notin \sui} \frac{ s_{1j}^{(3)}}{N^{3/2}} \ee[ \ea G_{j2} M_1 M_2 ] \notag  \\
= & \Osd ( |\lambda \Phi | N^{-1/2} \tilde{ \Psi} ) + \Osd( N^{-1} |\lambda| \tilde{ \Psi } )
\end{align}
with the second line following from the isotropic local law \eqref{eqn:iso}. This completes the proof. \qed


We now turn to computing the second order terms on the RHS of \eqref{eqn:line-expansion-1}. The following is fairly straightforward.
\bel \label{lem:line-second-1}
We have for $j \in \sui$ and $j \neq 2$,
\beq \label{eqn:line-second-1}
\frac{1}{N} | \ee[ \del_{j1} ( \ea G_{j2} M_1 M_2 ) ] | \sd (1 + |\lambda| )N^{-1} \tilde{\Psi}
\eeq
and for $j =2$,
\beq
\frac{1}{N} | \ee[ \del_{j1} ( \ea G_{j2} M_1 M_2 ) ] | \sd  N^{-1/2} \tilde{ \Psi}
\eeq
\eel
\proof These estimates follow from using \eqref{eqn:del-gen-1}, \eqref{eqn:del-gen-2} and items \ref{it:line-der-1}, \ref{it:line-der-2} and \ref{it:line-der-3} of Lemma \ref{lem:line-der} as appropriate. \qed

We recall here Definition \ref{def:D2} of the operation $D_2$ on hypergraphs. We define here also another operation $D_6$ on hypergraphs:
\bed
Given a hypergraph $\G = (V, E, \ui, \bu, \bw)$ with integration variables $\{ u_i \}_{i=1}^m$, an index $j$ and an edge $e \in E$ we define a new hypergraph $D_6 ( \G, j, e)$ with the same integration variables $\{ u_i \}_{i=1}^m$ as follows. We add a new vertex $v'$ and a self-edge $(v', v')$. We set $i'_{v'} = j$ and set $\bu'_{(v', v') } = \bu_e$, $\bw'_{(v', v')} = \bw_e$. The other variables are all the same as the original graph.
\eed

\bel \label{lem:line-second-2}
We have,
\begin{align}
 & \sum_{ j \notin \sui } \frac{ s_{1j}}{N} \ee[ ( \del_{j1} G_{j2} ) M_1 M_2 \ea ] = - \msc (z) \ee[ G_{12} (z) M_1 M_2 \ea ] + \Osd ( N^{-1} \tilde{ \Psi} ) \notag \\
- & \sum_{j \notin \sui } \frac{s_{1j}}{N} \ee[ \ea ( \M ( D_6 ( \G, j, (1, 2) ) ) + \M ( D_2 ( \G, (1, 2), j, (1, 2) ) ) ) ] . \label{eqn:line-second-2}
\end{align}
\eel
\proof Since $\del_{j1} G_{j2} = - G_{jj} G_{12} - G_{1j} G_{j2}$ we have,
\begin{align}
&\sum_{ j \notin \sui } \frac{ s_{1j}}{N} \ee[ ( \del_{j1} G_{j2} ) M_1 M_2 \ea ] = - \sum_{j \notin \sui} \frac{s_{1j}}{N} \msc (z) \ee[ G_{12} (z) M_1 M_2 \ea ] \notag\\
- & \sum_{j \notin \sui} \frac{s_{1j}}{N} \ee[ G_{12} M_1 M_2 (G_{jj} (z) - \msc (z) ) \ea ] - \sum_{j \notin \sui} \frac{s_{1j}}{N} \ee[ G_{1j} (z) G_{j2} (z) M_1 M_2 \ea ] \label{eqn:line-second-2a}
\end{align}
The first sum on the last line has terms of the form $\ee[ \ea \M (D_6 (\G, j, (1, 2) )]$. 
The second sum on the last line has terms  of the form $\ee[ \ea \M ( D_2 ( \G, (1, 2), j, (1, 2) )]$. The claim follows from this and the fact that $\sum_j \frac{s_{1j}}{N} = 1 + \O ( N^{-1} )$. \qed

We recall the definitions of the operations on hypergraphs $D_3, D_4, D_5$ given in Definitions \ref{def:D3D4} and \ref{def:D5}. 
\bel \label{lem:line-second-3}
\begin{align}
 & -\sum_{j \notin \sui } \frac{ s_{1j}}{N} \ee[ \ea G_{j2} \del_{j1} (M_1 M_2 ) ] \notag \\
= & \sum_{j \notin \sui}  \sum_{ e \in E \backslash \{ (1, 2) \} } \frac{ s_{1j}}{N} \ee[ \ea(  \M (D_3 ( \G, (1, 2), e, j) ) + \M (D_4 ( \G , (1, 2) , e, j ) ) ) ]  \label{eqn:line-second-3}
\end{align}
and 
\begin{align}
 & \sum_{j \notin \sui } \frac{ s_{1j}}{N} \ee[ G_{j2} M_1 M_2 \del_{j1} \ea ] \notag\\
 = & - \frac{ 2\i \lambda}{\pi} \sum_{j \notin \sui} \frac{ s_{1j}}{N} \int_{\Oma} \d u ( \del_{ \bar{u}} \tilf ( u ) ) \del_u \ee[ \ea \M(D_5 ( \G, (1, 2), j, u )) ]  \d u \d \bar{ u} \label{eqn:line-second-4}
\end{align}
\eel
\proof This is direct computation and so omitted. \qed

We collect the above estimates into the following expansion for line admissible hypergraphs.
\bep \label{prop:line-expand} Let $\G$ be a line admissible hypergraph as above, and assume  $| \lambda \Phi | \leq 1 $. We have,
\begin{align}
& \ee[ \ea G_{12} (z) M_1 M_2 ] = \msc(z)\sum_{j \notin \sui } \frac{s_{1j}}{N} \ee[ \ea ( \M ( D_6 ( \G, j, (1, 2) ) ) + \M ( D_2 ( \G, (1, 2), j, (1, 2) ) ) ) ] \notag\\
& + \msc(z) \sum_{j \notin \sui}  \sum_{ e \in E \backslash \{ (1, 2) \} } \frac{ s_{1j}}{N} \ee[ \ea(  \M (D_3 ( \G, (1, 2), e, j) ) + \M (D_4 ( \G , (1, 2) , e, j ) ) ) ]   \notag\\
& + \frac{2 \i \lambda \msc(z)}{\pi} \sum_{j \notin \sui} \frac{ s_{1j}}{N} \int_{\Oma} \d u ( \del_{ \bar{u}} \tilf ( u ) ) \del_u \ee[ \ea \M(D_5 ( \G, (1, 2), j, u )) ]  \d u \d \bar{ u} + \Osd ( N^{-1/2} \tilde{ \Psi} ) \label{eqn:line-exx}
\end{align}
\eep
\proof This follows from \eqref{eqn:line-expansion-1} and Lemmas \ref{lem:line-third}, \ref{lem:line-second-1}, \ref{lem:line-second-2} and \ref{lem:line-second-3}. \qed

\subsection{Iterative expansion} \label{sec:line-iterative}

We now carry out our iterative expansion for line-admissible hypergraphs. We begin by fixing some notation. Let $\G$ be a line-admissible hypergraph, and let $e= (v_1, v_2)$ be the first edge in this line.  The following families of hypergraphs will appear in our iterative expansion. 
\bed Let $\G$ and $e= (v_1, v_2)$ be as above. We define a set of hypergraphs $F_2 ( \G)$ to be the union of the following sets of hypergraphs.
\begin{enumerate}[label=\normalfont(\roman*)]
\item The graphs $D_6 ( \G , j, e)$ and $D_2 ( \G, e, j, e)$  for $j \notin \sui$
\item The graphs $D_3 ( \G, e, e_1, j)$ and $D_4 (\G, e, e_1, j)$ for $j \notin \sui$ and $e_1 \in E \backslash \{ e \}$.
\item The graphs $D_5 ( \G, e, j, u)$ for $j \notin \sui$.
\end{enumerate}
\eed
For any hypergraph $\G$ we recall that $E_w$ is the support of the vector $\bw$ and define,
\beq
\hat{ \Psi} ( \G ) := \prod_{ e \in E_w} \Psi ( z_e ).
\eeq
Let now $\G_0$ be a line-admissible hypergraph s.t. $\bu = 0$ and let $e = (v_1, v_2)$ be the first edge in this line. Set $\F_0 := \{ \G_0 \}$ and inductively set
\beq
\F_n := \bigcup_{ \Hcal \in \F_{n-1} } F_2 ( \Hcal ).
\eeq
By inspection, we find the following.
\bel
\begin{enumerate}[label=\normalfont(\roman*)]
\item For any $n \geq 1$, all hypergraphs in $\F_n$ contain $n$ more edges than $\G_0$
\item Every graph in $\F_n$ is line admissible
\item For every $H \in \F_n$ we have $\hat{ \Psi} ( \Hcal) \leq \tilde{ \Psi} ( \G_0 ) = \hat{\Psi} ( \G_0)$
\item If $\Hcal_1 \in F_2 ( \Hcal)$ for some $\Hcal \in \F_n$, then 
\beq
\tilde{ \Psi} ( \Hcal_1 ) = \Psi (v) \tilde{ \Psi}  ( \Hcal)
\eeq
for some $v$ that is a spectral parameter of $\Hcal$.
\end{enumerate}
\eel

\bel
Let $n \geq 0$. Uniformly in $\Hcal \in \F_n$ we have,
\beq \label{eqn:line-expand-1}
T ( \Hcal) = \sum_{ \Hcal_1 \in \F_2 ( \Hcal ) } b_{ \Hcal_1} T( \Hcal_1 ) + \Osd ( \tilde{ \Psi} (\G_0) N^{-1/2} )
\eeq
for some coefficients obeying, 
\beq
\sum_{ \Hcal_1 \in F_2 ( \Hcal) } |b_{ \Hcal_1} | \leq C_n .
\eeq
Let $\alpha >0$ satisfy $\alpha < \frac{ \mfa}{10} $ as well as $\min_{ z \in \G } N | \Im[z] | \geq N^{4 \alpha}$, and $| \lambda \Phi | \leq N^{-4 \alpha}$. Then for all $\Hcal \in \F_n$ we have,
\beq \label{eqn:line-expand-2}
| T ( \Hcal ) | \leq N^{ - n \alpha} .
\eeq
\eel
\proof The expansion \eqref{eqn:line-expand-1} follows from Proposition \ref{prop:line-expand} in the exact same way as Lemma \ref{lem:loops-expand-3} follows from Lemma \ref{lem:loops-expand-2}, and so we omit the detailed proof. Similarly, \eqref{eqn:line-expand-2} follows in the same way as Lemma \ref{lem:loops-expand-4}. \qed

\subsubsection{Proof of Proposition \ref{prop:main-line-estimate}} \label{sec:main-line-proof}

By iterating \eqref{eqn:line-expand-1} we see that,
\beq
| \ee[ \ea \M ( \G ) ] | \sd \max_{ \Hcal \in \F_n} | T ( \Hcal ) | + N^{-1/2} \tilde{ \Psi} ( \G ).
\eeq
The claim now follows from taking $n$ sufficiently large and using \eqref{eqn:line-expand-2}. \qed

\subsection{A better estimate for $G_{ij}$} \label{sec:gij-better}

In this section we derive an even better estimate for $\ee[ \ea G_{ij} (z) ]$ with $i \neq j$ than was derived in Proposition \ref{prop:main-line-estimate}. We start first with the following. We set $\tt = \{ 1, 2 \}$. 
\bel For $|\lambda \Phi| \leq 1$ we have, \label{lem:gij-expand}
\begin{align}
& z \ee[ \ea G_{12} (z) ] = \sum_{ j \neq 1} \frac{s_{1j}}{N} \ee[ \del_{j1} (G_{j2} (z) \ea ) ] \notag\\
+ & \Osd (  |\lambda| \Phi N^{-1/2} \Psi (z) + N^{-1/2} \Psi(z)^2 ). \label{eqn:gij-expand}
\end{align}
\eel
\proof
 First, by combining Lemma \ref{lem:line-expand-1}, \eqref{eqn:line-third-1}, \eqref{eqn:line-third-2}, \eqref{eqn:line-third-4} and \eqref{eqn:line-second-1} we see that,
 \begin{align} \label{eqn:gij-expand-a}
&  z \ee[ \ea G_{12} (z) ] = \sum_{ j \neq 1} \frac{s_{1j}}{N} \ee[ \del_{1j} (G_{j2} (z) \ea ) ] \notag\\
+& \sum_{j \notin \tt} \frac{s_{1j}^{(3)}}{ 2 N^{3/2}} \ee[ \ea \del_{1j}^2 (G_{j2} ) ] + \Osd ( (1 + |\lambda| ) N^{-1} \Psi (z) + |\lambda| \Phi N^{-1/2} \Psi (z) ). 
 \end{align}
We now compute that for $j\notin \tt$ we have $\del_{1j}^2 G_{j2} = 2 \msc (z)^2 G_{j2} (z) + \Osd ( \Psi(z)^2)$ to see that,
\beq \label{eqn:gij-expand-b}
\sum_{j \notin \tt} \frac{s_{1j}^{(3)}}{2N^{3/2}} \ee[ \ea \del_{1j}^2 (G_{j2} ) ] =  \msc(z)^2 \sum_{j \notin \tt} \frac{s_{1j}^{(3)}}{N^{3/2}} \ee[ \ea(G_{j2} ) ] + \Osd ( N^{-1/2} \Psi(z)^2 )
\eeq
and the first term on the RHS is $\Osd ( N^{-1} \Psi (z) )$ by the istropic local law \eqref{eqn:iso}. \qed

We now turn to computing the second order terms.
\bel \label{lem:gij-second}
Assume $|\lambda \Phi | \leq 1$. We have,
\beq \label{eqn:g12-second-1}
N^{-1} | \ee[ \del_{21} (\ea G_{22} ) ] | \sd N^{-1} |\lambda \Phi| + N^{-1} \Psi (z).
\eeq
and
\beq \label{eqn:g12-second-2}
\sum_{j \notin \tt} \frac{ s_{1j}}{N} \ee[ G_{j2} (z) \del_{j1} \ea ] = \Osd ( |\lambda \Phi| N^{-1/2} \Psi (z) )
\eeq
and 
\beq \label{eqn:g12-second-3}
\sum_{j \notin \tt} \frac{s_{1j}}{N} \ee[ \ea \del_{j1} G_{j2} (z) ] = - \msc(z) \ee[ G_{12} (z) \ea ] + \Osd ( N^{-1/2} \Psi(z)^2 ).
\eeq
\eel
\proof For \eqref{eqn:g12-second-1} we have $\del_{21} G_{22} = \Osd ( \Psi (z))$ and $\del_{12} \ea = \O ( |\lambda| \Phi )$ by \eqref{eqn:del-gen-1} and so the estimate follows. For \eqref{eqn:g12-second-2} we have that,
\beq
\ee[ G_{j2} (z) \del_{j1} \ea ] = - \frac{2 \i \lambda}{\pi} \int_{\Oma} ( \del_{\bar{u}} \tilf (u) ) \del_u \ee[ \ea G_{1j} (u) G_{j2} (z) ]  \d u \d \bar{u}.
\eeq
From Proposition \ref{prop:main-line-estimate} we have $ \ee[ \ea G_{1j} (u) G_{j2} (z) ]  = \Osd ( N^{-1/2} \Psi(z) \Psi (u) )$ and so \eqref{eqn:g12-second-2} follows by applying \eqref{eqn:H-est}. For \eqref{eqn:g12-second-3} we have,
\begin{align} \label{eqn:gij-second-3a}
& -\sum_{j \notin \tt} \frac{s_{1j}}{N} \ee[ \ea \del_{j1} G_{j2} (z) ] = \sum_{j \notin \tt} \frac{s_{1j}}{N} \msc(z) \ee[ G_{12} (z) ] \notag\\
+ & \sum_{j \notin \tt} \frac{s_{1j}}{N} \ee[ \ea  G_{12} (z) (G_{jj} (z) - \msc (z) ] + \sum_{j \notin \tt} \frac{s_{1j}}{N} \ee[ G_{1j} (z) G_{j2} (z) \ea ] \notag\\
= & \msc(z) \ee[ \ea G_{12} (z) ] + \Osd ( \Psi(z)^2 N^{-1/2} )
\end{align} 
with the last line following from applying Proposition \ref{prop:main-line-estimate} twice, as well as that $\sum_{j \notin \tt} \frac{s_{1j}}{N} = 1 + \O ( N^{-1} )$. 
\qed

The above two lemmas easily imply the following, which improves upon the estimate Proposition \ref{prop:main-line-estimate} in the special case of a single resolvent entry.
\bep \label{prop:gij-est}
If there is a $c>0$ so that $| \lambda \Phi | \leq N^{-c}$ then for $i \neq j$ we have,
\beq
| \ee[ G_{ij} (z) \ea ] | \sd | \lambda \Phi | N^{-1/2} \Psi (z) + N^{-1/2} \Psi(z)^2 .
\eeq
\eep

\section{Preliminary expansion for $\ee[ \ea G_{11}]$} \label{sec:prelim}

In this section we derive an expansion for $\ee[ \ea G_{ii} (z) ]$. We will use it to derive an estimate for the expression $\ee[ \ea (G_{ii} (z) - \msc(z ) ]$ that is better than the naive estimate $\Osd ( \Psi (z))$. Secondly, the first part of our expansion will be used to derive our main estimates on the characteristic function. Throughout this section we will use the notation $z = E + \i \eta$, i.e., $\eta = | \Im[z]|$.  Throughout this section we will always assume that,
\beq
|\lambda \Phi| \leq N^{-c}
\eeq
for some $c>0$ and that $z \in \Oma$. 

 We start with the following.
 
 \bel \label{lem:G11-fifth} 
 We have that,
 \beq
 N^{-5/2} | \del_{11}^4 ( \ea G_{11} (z) ) | \sd (1 + |\lambda| )^4 N^{-5/2}.
 \eeq
 and that for $j \neq 1$,
 \beq \label{eqn:G11-fifth-der-2}
 N^{-5/2}  | \del_{j1}^4 ( \ea G_{1j} (z) ) | \sd  N^{-5/2} \Psi(z)(1 + |\lambda| )^2 + N^{-5/2} |\lambda \Phi| (1 + |\lambda| )
 \eeq
 For any $\delta >0$, the same estimates hold for any $\eps$-regular matrix with an additional factor of $N^{\delta}$ on the RHS, if $\eps$ is sufficiently small depending on $\delta$.
 \eel
 \proof The estimate for $j=1$ follows directly from \eqref{eqn:del-gen-2}. We turn to the estimate for $j \neq 1$. Note that for $k=0, 2, 4$ we have that,
 \beq
 \del_{1j}^k G_{1j} = \Osd ( \Psi (z) )
 \eeq
 by direct computation. The LHS is $\Osd (1)$ otherwise. Furthermore, for $k=1, 3$ we have 
 \beq \label{eqn:delea3}
 \del_{1j}^k \ea = \Osd ( |\lambda \Phi| (1 + |\lambda| ) ,
 \eeq
 and for $k=2, 4$ it is $\Osd ( (1 + |\lambda| )^2)$. Here, the $k=1$ case follows from \eqref{eqn:deleaij}, and the $k=2, 4$ cases from \eqref{eqn:del-gen-1}. The $k=3$ case follows from a direct computation similar to \eqref{eqn:deleaij2}.  The estimate \eqref{eqn:G11-fifth-der-2} now follows. The statement about $\eps$-regular matrices is now straightforward. \qed

 We now apply the cumulant expansion to arrive at our starting point for estimating $\ee[ \ea G_{11} (z) ]$.

\bel \label{lem:prelim-expand-1}
We have,
\begin{align} \label{eqn:prelim-expand-1}
 & z \ee[ \ea G_{11} (z) ] + \ee[ \ea] = \sum_{j} \frac{s_{1j}}{N} \ee[ \del_{1j} ( \ea G_{j1} ) ] + \sum_{j } \frac{s_{1j}^{(3)}}{2 N^{3/2}} \ee[ \del_{1j}^2 ( \ea G_{j1} ) ] \notag\\
+ &   \sum_{j } \frac{s_{1j}^{(4)}}{6 N^{2}} \ee[ \del_{1j}^3 ( \ea G_{j1} ) ] + \Osd ( (1 + |\lambda| )^4 N^{-5/2} + \Psi(z) N^{-3/2} ( 1 + |\lambda| )^2 + N^{-3/2} |\lambda \Phi| ( 1 + |\lambda| ) )
\end{align}
\eel
\proof This follows from the cumulant expansion of Lemma \ref{lem:cumu-exp} to fifth order applied to,
\beq
z \ee[ G_{11} \ea ] + \ee[ \ea] = \sum_j \ee[ H_{1j} G_{j1} \ea ],
\eeq
using the estimates of Lemma \ref{lem:G11-fifth}.  \qed

In Sections \ref{sec:prelim-fourth} and \ref{sec:prelim-third} below, we will compute the fourth and third order terms, respectively, on the RHS of \eqref{eqn:prelim-expand-1}. The results will then be summarized in Section \ref{sec:intermediate}.

\subsection{Fourth order terms} \label{sec:prelim-fourth}

We first prove the following.
\bel \label{lem:prelim-fourth-1}
We have,
\begin{align}
& \sum_j \frac{ s_{1j}^{(4)}}{N^2} \ee[ \del_{1j}^3 ( \ea G_{j1} ) ] = - 6 \sum_{j \neq 1} \frac{ s_{1j}^{(4)}}{N^2} \msc(z)^4 \ee[ \ea ] \notag\\
+ & 3 \sum_{j \neq 1 } \frac{ s_{1j}^{(4)}}{N^2} \ee[ ( \del_{1j}^2 \ea ) ( \del_{1j} G_{j1} ) ] 
- 12 \msc(z)^3 \sum_{j \neq 1} \frac{ s_{1j}^{(4)}}{N^2} \ee[ \ea (G_{11} (z) - \msc (z) ) ] \notag\\
+& \Osd ( N^{-2} \eta^{-1}  + ( 1 + |\lambda| )^3 N^{-2} + |\lambda \Phi| \Psi(z) (1 + |\lambda| ) N^{-1} ) \label{eqn:prelim-fourth-1}
\end{align}
\eel
\proof The term when $j=1$ can be estimated by $\Osd(N^{-2} (1 + |\lambda| )^3 )$ using \eqref{eqn:del-gen-2}. When $j \neq 1$ we have that by direct computation and \eqref{eqn:entry-wise}
\beq \label{eqn:prelim-fourth-1a}
\del_{j1}^3 G_{j1} = - 6 G_{11}^2 G_{jj}^2 + \O ( \Psi(z)^2)
\eeq
and so,
\begin{align} \label{eqn:prelim-fourth-1a-1}
& \frac{-1}{6}\sum_{j \neq 1}  \frac{ s_{1j}^{(4)}}{N^2} \ee[\ea \del_{1j}^3 (  G_{j1} ) ] =  \sum_{j \neq 1}  \frac{ s_{1j}^{(4)}}{N^2} \ee[\ea G_{11} (z)^2 G_{jj}(z)^2] + \Osd ( N^{-1} \Psi(z)^2 ) \notag\\
= &   \sum_{j \neq 1} \frac{ s_{1j}^{(4)}}{N^2} \left\{ \msc(z)^4\ee[\ea] + 2 \msc(z)^3   \ee[\ea (G_{11} (z) + G_{jj} (z) - 2\msc (z) ) ]\right\} + \Osd ( N^{-1} \Psi (z)^2 ) \notag\\
= & \sum_{j \neq 1 } \frac{ s_{1j}^{(4)}}{N^2} \msc(z)^4 \ee[\ea] +2 \msc(z)^3 \sum_{j \neq 1 } \frac{ s_{1j}^{(4)}}{N^2} \ee[ \ea (G_{11} (z) - \msc (z) ) ] + \Osd (N^{-2} \eta^{-1} )
\end{align}
The second line follows from writing $G_{ii} = (G_{ii} - \msc) + \msc$ and expanding, and dropping all contributions quadratic and higher in $(G_{ii} - \msc)$ (and similarly for the $G_{jj}$ term). The last line follows from fluctuation averaging \eqref{eqn:fluct-av}. This takes care of the terms when all of the $\del_{ij}$ hit $G_{j1}$. 

We estimate the other contributions using,
\beq
( \del_{1j} \ea ) ( \del_{1j}^2 G_{1j} )  = \Osd ( |\lambda| \Phi \Psi (z) ) , \qquad \del_{1j}^3 \ea = \Osd ( |\lambda \Phi| ( 1 + |\lambda| ) )
\eeq
which hold for $j \neq 1$. The first estimate follows from \eqref{eqn:deleaij} and \eqref{eqn:entry-wise}, and the second from \eqref{eqn:delea3}.  \qed

For the remaining fourth order terms we compute,
\bel \label{lem:prelim-fourth-2}
We have,
\begin{align} \label{eqn:prelim-fourth-2}
&\sum_{j \neq 1 } \frac{ s_{1j}^{(4)}}{N^2} \ee[  ( \del_{1j}^2 \ea ) ( \del_{1j} G_{j1} ) ] = \Osd ( |\lambda|N^{-2} \eta^{-1} + |\lambda|(1+ |\lambda| ) \Phi^2 N^{-1} )   \notag\\
& -  \ee[ \ea] \msc(z)^2 \sum_{j \neq 1 } \frac{ s_{1j}^{(4)}}{N^2} \frac{2 \i \lambda}{\pi} \int_{\Oma} ( \del_{\bar{u}} \tilf (u) ) \del_u \msc(u)^2 \d u \d \bar{u} \notag\\
 &-  \sum_{j \neq 1} \frac{ s_{1j}^{(4)}}{N^2} \frac{2 \i \lambda}{\pi} \int_{\Oma} ( \del_{\bar{u}} \tilf (u) ) \del_u ( \msc(z)^2 \msc(u) \ee[ \ea(  G_{11} (u) - \msc(u) ) ] ) \d u \d \bar{u} \notag\\ 
  &- \sum_{j \neq 1} \frac{ s_{1j}^{(4)}}{N^2}\frac{2 \i \lambda}{\pi} \int_{\Oma} ( \del_{\bar{u}} \tilf (u) ) \del_u (\msc(u)^2 \msc(z) \ee[ \ea(  G_{11} (z)  -\msc(z) ) ] ) \d u \d \bar{u} \notag\\ 
\end{align}
\eel
\proof By using that 
\beq \label{eqn:prelim-fourth-2a}
\del_{1j} G_{j1} =  - G_{jj} G_{11} + \Osd ( \Psi(z)^2 )
\eeq
 and 
\beq \label{eqn:prelim-fourth-2b}
\del_{1j}^2 \ea = \frac{ 2 \i \lambda}{\pi} \int_{ \Oma} ( \del_{\bar{u}} \tilf (u) ) \del_u ( G_{11} (u) G_{jj} (u) ) \d u \d \bar{u} + \Osd ( | \lambda \Phi^2| (1 + |\lambda| ) ),
\eeq
(the above following from \eqref{eqn:delea2})
we see that,
\begin{align}
 -& \sum_{j \neq 1 } \frac{ s_{1j}^{(4)}}{N^2} \ee[  ( \del_{1j}^2 \ea ) ( \del_{1j} G_{j1} ) ] = \Osd ( |\lambda| \Psi(z)^2 N^{-1} + |\lambda|(1+ |\lambda|) \Phi^2 N^{-1} ) \notag\\
 +& \sum_{j \neq 1} \frac{ s_{1j}^{(4)}}{N^2} \frac{2 \i \lambda}{\pi} \int_{\Oma} ( \del_{\bar{u}} \tilf (u) ) \del_u \ee[ \ea G_{11} (u) G_{jj} (u) G_{11} (z) G_{jj} (z) ] \d u \d \bar{u} \notag\\
 = & \ee[ \ea] \msc(z)^2 \sum_{j \neq 1 } \frac{ s_{1j}^{(4)}}{N^2} \frac{2 \i \lambda}{\pi} \int_{\Oma} ( \del_{\bar{u}} \tilf (u) ) \del_u \msc(u)^2 \d u \d \bar{u} \notag\\
 + & \sum_{j \neq 1} \frac{ s_{1j}^{(4)}}{N^2} \frac{2 \i \lambda}{\pi} \int_{\Oma} ( \del_{\bar{u}} \tilf (u) ) \del_u (\msc(z)^2 \msc(u) \ee[ \ea(  G_{11} (u) + G_{jj} (u)- 2 \msc(u) ) ] ) \d u \d \bar{u} \notag\\ 
  & \sum_{j \neq 1} \frac{ s_{1j}^{(4)}}{N^2} \frac{2 \i \lambda}{\pi}\int_{\Oma} ( \del_{\bar{u}} \tilf (u) ) \del_u ( \msc(u)^2 \msc(z) \ee[ \ea(  G_{11} (z) + G_{jj} (z)- 2 \msc(z) ) ] ) \d u \d \bar{u} \notag\\ 
 + & \Osd ( |\lambda| \Psi(z)^2 N^{-1} + |\lambda|(1+ |\lambda| ) \Phi^2 N^{-1} ) 
\end{align}
with the second estimate following from linearizing the $G_{ii}$ and $G_{jj}$ around $\msc$, similar to \eqref{eqn:prelim-fourth-1a-1}. The terms with $G_{jj}$ on the third last and second last lines can now be estimated by $\Osd ( |\lambda| \Phi^2 N^{-1} )$ and $\Osd ( N^{-1} |\lambda| \eta^{-1})$, respectively, using fluctuation averaging \eqref{eqn:fluct-av} and \eqref{eqn:H-est}. This yields the claim. \qed

\subsection{Third order terms} \label{sec:prelim-third}

We split up the estimation of the third order terms on the RHS of \eqref{eqn:prelim-expand-1} into the diagonal and off-diagonal terms.

\subsubsection{Diagonal terms}

Starting with the diagonal term we have,
\bel \label{lem:prelim-third-1}
We have, 
\begin{align}
& \frac{ s_{11}^{(3)}}{N^{3/2}} \ee[ \del_{11}^2 ( \ea G_{11} (z) ) ] =  \frac{ s_{11}^{(3)}}{N^{3/2}} \ee[ \ea]\bigg\{ \msc(z)^3 +  2 \msc(z)^2 \left( \frac{ \i \lambda}{\pi}  \int_{\Oma} \d u \d \bar{u} \del_{\bar{u}} \tilf (u) \del_u \msc(u) \right)  \notag\\
+ & \msc(z) \left( - \frac{ \i \lambda}{\pi}  \int_{\Oma} \d u \d \bar{u} \del_{\bar{u}} \tilf (u) \del_u \msc(u) \right)^2 + \msc(z)  \left( \frac{ \i \lambda}{\pi}  \int_{\Oma} \d u \d \bar{u} \del_{\bar{u}} \tilf (u) \del_u \msc(u)^2 \right)  \bigg\} \notag\\
&+ \Osd ( N^{-3/2} \Psi (z)(1 + |\lambda| )^2 + |\lambda| ( 1 + |\lambda| ) \Phi N^{-3/2} )  \label{eqn:prelim-third-4}
\end{align}
\eel
\proof Using that $\del_{11}^2 G_{11} = \msc(z)^3 + \O ( \Psi )$ we see that,
\beq
\frac{ s_{11}^{(3)}}{N^{3/2}} \ee[ \ea \del_{11}^2 G_{11} (z) ] = \msc(z)^3 \frac{ s_{11}^{(3)}}{N^{3/2}} \ee[ \ea] + \Osd (N^{-3/2} \Psi(z) ).
\eeq
Using the fact that $\del_{11}^k G_{11} = (-1)^k (G_{11} )^{k+1} = (-1)^k\msc(z)^{k+1} + \O ( \Psi)$ and that,
\begin{align}
\del_{11} \ea &= - \frac{ \i \lambda}{\pi} \ea \int_{\Oma} \d u \d \bar{u} \del_{\bar{u}} \tilf (u) \del_u \msc(u) \d u \d \bar{u} + \O ( |\lambda| \Phi ) \notag\\
\del_{11}^2 \ea &= \ea \left( - \frac{ \i \lambda}{\pi}  \int_{\Oma} \d u \d \bar{u} \del_{\bar{u}} \tilf (u) \del_u \msc(u) \d  u \d \bar{u} \right)^2 \notag \\
&+ \ea \left( \frac{ \i \lambda}{\pi}  \int_{\Oma} \d u \d \bar{u} \del_{\bar{u}} \tilf (u) \del_u \msc(u)^2 \d u \d \bar{u} \right) + \O ( \Phi |\lambda| (1 + |\lambda| ) )
\end{align}
we conclude the proof. We note that the above estimates follow from computing the derivatives explicitly (similar to the first equality in \eqref{eqn:del-gen-1}), and then replacing all the appearances of $G_{ii}(u)$ by $\msc (u)$ using \eqref{eqn:entry-wise} and \eqref{eqn:H-est}. \qed

\subsubsection{Off-diagonal terms}

The next few lemmas deal with the off-diagonal third order terms on the RHS of \eqref{eqn:prelim-expand-1}. 

\bel \label{lem:prelim-third-2}
We have, 
\begin{align} \label{eqn:prelim-third-1}
& \sum_{ j \neq 1 } \frac{ s_{1j}^{(3)}}{N^{3/2}} \ee[ \ea \del_{j1}^2 G_{j1} (z)] = 6 \msc (z)\sum_{j \neq 1 } \frac{ s_{1j}^{(3)}}{N^{3/2}} \ee[ \ea ( G_{jj} -\msc(z) )G_{1j} (z) ] \notag\\
+ & \Osd ( N^{-1/2} \Psi(z)^3 + |\lambda \Phi | N^{-1} \Psi(z) )
\end{align}
\eel
\proof Since 
\beq \label{eqn:prelim-third-1a}
\del_{j1}^2 G_{j1} = 6 G_{11} G_{jj} G_{1j} + \Osd ( \Psi(z)^3)
\eeq
 we have first that,
\begin{align}
\sum_{ j \neq 1 } \frac{ s_{1j}^{(3)}}{N^{3/2}} \ee[ \ea \del_{j1}^2 G_{j1} ] = 6 \sum_{ j \neq 1 } \frac{ s_{1j}^{(3)}}{N^{3/2}} \ee[ \ea G_{11} G_{jj} G_{j1} ] + \Osd ( N^{-1/2} \Psi(z)^3 ). 
\end{align}
Now we write,
\begin{align}
 & \sum_{j \neq 1 } \frac{ s_{1j}^{(3)}}{N^{3/2}} G_{11} G_{jj} G_{1j} = \msc(z)^2\sum_{j \neq 1 } \frac{ s_{1j}^{(3)}}{N^{3/2}}G_{1j} + \msc(z) \sum_{j \neq 1 } \frac{ s_{1j}^{(3)}}{N^{3/2}} (G_{11}- \msc(z)) G_{1j} \notag\\
+ & \msc (z)\sum_{j \neq 1 } \frac{ s_{1j}^{(3)}}{N^{3/2}} ( G_{jj} -\msc(z) )G_{1j} + \Osd( N^{-1/2} \Psi(z)^3 ) \notag\\
= & \msc(z)^2\sum_{j \neq 1 } \frac{ s_{1j}^{(3)}}{N^{3/2}}G_{1j}  + \msc (z)\sum_{j \neq 1 } \frac{ s_{1j}^{(3)}}{N^{3/2}} ( G_{jj} -\msc(z) )G_{1j} + \Osd ( N^{-1/2} \Psi(z)^3 )
\end{align}
where we used the isotropic local law \eqref{eqn:iso} in the second estimate. Now by Proposition \ref{prop:gij-est} we have,
\beq
\sum_{j \neq 1 } \frac{ s_{1j}^{(3)}}{N^{3/2}} \ee[ \ea G_{1j} (z)] = \Osd (N^{-1} \Psi(z)^2 + | \lambda \Phi | N^{-1} \Psi (z) )
\eeq
which completes the proof. \qed

\bel \label{lem:prelim-third-3}
We have,
\begin{align} \label{eqn:prelim-third-2}
& \sum_{ j \neq 1 } \frac{ s_{1j}^{(3)}}{N^{3/2}} \ee[ ( \del_{1j} \ea ) ( \del_{1j} G_{1j} (z) ) ] =  -\sum_{ j \neq 1 } \frac{ \msc(z) s_{1j}^{(3)}}{N^{3/2}} \ee[ ( \del_{1j} \ea ) ( G_{jj} (z)-  \msc (z) ) ] \notag\\
+ & \Osd ( N^{-1/2} |\lambda \Phi | \Psi(z)^2 + N^{-1} |\lambda| ( 1 + |\lambda| ) \Phi^2 )
\end{align}
\eel
\proof Since $\del_{1j} G_{1j} = - G_{jj} G_{11} + \Osd ( \Psi(z)^2)$ and $\del_{1j} \ea = \Osd ( | \lambda \Phi |)$ by \eqref{eqn:deleaij} we have,
\begin{align}
& \sum_{ j \neq 1 } \frac{ s_{1j}^{(3)}}{N^{3/2}} \ee[ ( \del_{1j} \ea ) ( \del_{1j} G_{1j} ) ] = -\sum_{ j \neq 1 } \frac{ s_{1j}^{(3)}}{N^{3/2}} \ee[ ( \del_{1j} \ea ) G_{11} G_{jj} ] + \Osd (N^{-1/2} |\lambda \Phi | \Psi(z)^2 ) \notag\\
= & -\msc(z)^2 \sum_{ j \neq 1 } \frac{ s_{1j}^{(3)}}{N^{3/2}} \ee[ ( \del_{1j} \ea )  ]- \sum_{ j \neq 1 } \frac{ \msc(z) s_{1j}^{(3)}}{N^{3/2}} \ee[ ( \del_{1j} \ea ) (G_{11} + G_{jj} - 2 \msc (z) ) ] \notag\\
+ &  \Osd (N^{-1/2} |\lambda \Phi | \Psi(z)^2 ).
\end{align}
In the second estimate we linearized the $G_{ii}$ and $G_{jj}$ around $\msc$. 
We have,
\begin{align}
&  \sum_{ j \neq 1 } \frac{ s_{1j}^{(3)}}{N^{3/2}} \ee[ ( \del_{1j} \ea ) (G_{11}  -  \msc (z) ) ] \notag\\
 =&  - \frac{2 \i \lambda }{\pi} \int_{\Oma} ( \del_{\bar{u}} \tilf (u) ) \del_u \sum_{j \neq 1} \frac{ s_{1j}^{(3)}}{N^{3/2}} \ee[ \ea G_{1j} (u) (G_{11} (z) - \msc(z) ) ] \d u \d \bar{u} = \Osd ( |\lambda \Phi | N^{-1} \Psi(z) )
\end{align}
by the isotropic local law \eqref{eqn:iso} and \eqref{eqn:H-est}. Secondly, we have for $j \neq 1$,
\beq
\ee[ \del_{1j} \ea ] = \Osd( N^{-1/2} |\lambda| ( 1 + |\lambda| ) \Phi^2 )
\eeq
by applying Proposition \ref{prop:gij-est} and \eqref{eqn:H-est} to the formula in the first equality of \eqref{eqn:deleaij}. This completes the proof. \qed

\bel \label{lem:prelim-third-4}
We have,
\begin{align} \label{eqn:prelim-third-3}
& \sum_{ j \neq 1 } \frac{ s_{1j}^{(3)}}{N^{3/2}} \ee[ G_{1j} (z) \del_{1j}^2 \ea ]  \notag\\
=& \frac{2 \i \lambda}{\pi} \sum_{ j \neq 1 } \frac{ s_{1j}^{(3)}}{N^{3/2}} \int_{\Oma} \d u \d \bar{u} ( \del_{ \bar{u}} \tilf (u) ) \del_u( \msc (u) \ee[ \ea G_{1j} (z) ( G_{jj} (u) -  \msc (u) ) ] ) \notag\\
+ & \Osd ( N^{-1/2} |\lambda \Phi^2| ( 1+ |\lambda|) \Psi (z) + |\lambda| N^{-1} \Psi (z)^2 )
\end{align}
\eel
\proof We have by \eqref{eqn:prelim-fourth-2b},
\begin{align}
& \sum_{ j \neq 1 } \frac{ s_{1j}^{(3)}}{N^{3/2}} \ee[ G_{1j} (z) \del_{1j}^2 \ea ]  = \sum_{ j \neq 1 } \frac{ s_{1j}^{(3)}}{N^{3/2}}\ee[ \ea G_{1j} (z) ] \frac{2 \i \lambda}{\pi}  \int_{\Oma} \d u \d \bar{u} ( \del_{ \bar{u}} \tilf (u) ) \del_u \msc(u)^2 \notag\\
+ & \sum_{ j \neq 1 } \frac{ s_{1j}^{(3)}}{N^{3/2}} \frac{2 \i \lambda}{\pi} \int_{\Oma} \d u \d \bar{u} ( \del_{ \bar{u}} \tilf (u) ) \del_u ( \msc (u) \ee[ \ea G_{1j} (z) (G_{11} (u)  + G_{jj} (u) - 2 \msc (u) ) ] )  \notag\\
+ & \Osd ( N^{-1/2} \Phi^2 |\lambda| (1 + |\lambda| ) \Psi (z)  )
\end{align}
The term with $G_{11} (u) - \msc(u)$ on the second line contributes $\Osd ( N^{-1} |\lambda| \Phi \Psi (z) )$ using the isotropic law \eqref{eqn:iso}. We use Proposition \ref{prop:gij-est} to bound the term on the first line by $\Osd ( |\lambda| N^{-1} | \lambda \Phi| \Psi(z) + |\lambda| N^{-1} \Psi (z)^2 )$. The claim follows. \qed

\subsection{Intermediate expansion}

 \label{sec:intermediate} 
 
\bep \label{prop:intermediate-expansion}
We have that,
\begin{align}
& z  \sum_a \ee[ \ea G_{aa} (z)] + N \ee[ \ea ] = \sum_{ja} \frac{ s_{1j}}{N} \ee[ \del_{aj} ( \ea G_{ja} ) ] \notag \\
&  -  \sum_{j \neq a } \frac{ s_{aj}^{(4)}}{N^2} \msc(z)^4 \ee[ \ea ] - \ee[ \ea] \msc(z)^2 \sum_{j \neq a } \frac{ s_{aj}^{(4)}}{N^2} \frac{\i \lambda}{\pi} \int_{\Oma} ( \del_{\bar{u}} \tilf (u) ) \del_u \msc(u)^2 \d u \d \bar{u}  \notag\\
&+ \sum_a  \frac{ s_{aa}^{(3)}}{2 N^{3/2}} \ee[ \ea] \bigg\{ \msc(z)^3 +  2 \msc(z)^2 \left( \frac{ \i \lambda}{\pi} \ea \int_{\Oma} \d u \d \bar{u} \del_{\bar{u}} \tilf (u) \del_u \msc(u) \right)  \notag\\
+ & \msc(z) \left( - \frac{ \i \lambda}{\pi}  \int_{\Oma} \d u \d \bar{u} \del_{\bar{u}} \tilf (u) \del_u \msc(u) \right)^2 
+ \msc(z)  \left( \frac{ \i \lambda}{\pi}  \int_{\Oma} \d u \d \bar{u} \del_{\bar{u}} \tilf (u) \del_u \msc(u)^2 \right) \bigg\} \notag\\
+ &\Osd ( N^{-1} \eta^{-1} (1 + |\lambda| ) + ( 1 + |\lambda| )^3 N^{-1} +N^{1/2} \Psi(z)^3 + N^{1/2} |\lambda \Phi | \Psi(z)^2) \notag \\
+ &  \Osd ( N^{1/2} |\lambda \Phi^2| ( 1+ |\lambda|) \Psi (z)  ) \label{eqn:intermediate-expansion} 
\end{align}
\eep
\proof The proof starts by applying Lemma \ref{lem:prelim-expand-1} to each term on the LHS of \eqref{eqn:intermediate-expansion}  but with the index $1$ in \eqref{eqn:prelim-expand-1} replaced  by $a$. Then, one uses Lemmas \ref{lem:prelim-fourth-1}, \ref{lem:prelim-fourth-2}, \ref{lem:prelim-third-1}, \ref{lem:prelim-third-2}, \ref{lem:prelim-third-3} and \ref{lem:prelim-third-4} to deal with the terms that arise on the  RHS of \eqref{eqn:prelim-expand-1} (but again with the index $1$ replaced by $a$). The only issue is that there are some extra terms in \eqref{eqn:prelim-fourth-1}, \eqref{eqn:prelim-fourth-2}, \eqref{eqn:prelim-third-1}, \eqref{eqn:prelim-third-2} and \eqref{eqn:prelim-third-3} that we must estimate (they are simpler to estimate after taking the sum over $a$) in order to arrive at \eqref{eqn:intermediate-expansion}. Therefore, the remainder of the proof is devoted to estimating these terms.

The first extra term is the second term on the second line of \eqref{eqn:prelim-fourth-1}. We estimate it by,
\begin{align}
\sum_{j} \sum_{a \neq j } \frac{ s_{aj}^{(4)}}{N^2} \ee[ \ea ( G_{aa}(z) - \msc (z) ) ] = \Osd (N^{-1} \eta^{-1} )
\end{align}
using fluctuation averaging \eqref{eqn:fluct-av}. The next extra terms are those on the third and fourth lines of \eqref{eqn:prelim-fourth-2}. By applying \eqref{eqn:fluct-av}, one sees that the term on the third line is $\Osd ( |\lambda| \Phi^2)$ and the term on the fourth line is $\Osd (|\lambda| N^{-1} \eta^{-1} )$ (when we take the sum over $a$).

 The next extra term is the term on the RHS of the first line of \eqref{eqn:prelim-third-1}. We estimate it by,
\begin{align}
 & \sum_{j} \sum_{a \neq j } \frac{ s_{aj}^{(3)}}{N^{3/2}} \ee[ \ea (G_{jj}(z) - \msc(z) ) G_{aj} (z) ] \notag\\
=&  \sum_j \ee \left[ \ea (G_{jj} - \msc(z) ) \sum_{ j \neq a} \frac{ s_{aj}^{(3)}}{N^{3/2}} G_{aj} (z) \right] = \Osd ( \Psi(z)^2)
\end{align}
using the isotropic local law \eqref{eqn:iso} for the sum over the index $a$. By almost the same argument one finds an estimate of $\Osd ( |\lambda \Phi| \Psi (z) )$ for the terms on the RHS of the first lines of \eqref{eqn:prelim-third-2} and \eqref{eqn:prelim-third-3}. This completes the proof. \qed

As a simple corollary we find, 
\bec
We have,
\begin{align}
&z  \sum_a \ee[ \ea (G_{aa} (z) - \ee[ G_{aa} (z) ] ) ] = \sum_{ja} \frac{ s_{aj}}{N} \left( \ee[ \del_{aj} ( \ea G_{ja} (z) ) ] - \ee[ \ea] \ee[ \del_{aj} G_{ja} (z) ] \right) \notag\\
& - \ee[ \ea] \msc(z)^2 \sum_{j \neq a } \frac{ s_{aj}^{(4)}}{N^2} \frac{ \i \lambda}{\pi} \int_{\Oma} ( \del_{\bar{u}} \tilf (u) ) \del_u \msc(u)^2 \d u \d \bar{u}  \notag\\
&+ \sum_a  \frac{ s_{aa}^{(3)}}{2 N^{3/2}} \ee[ \ea] \bigg\{   2 \msc(z)^2 \left( \frac{ \i \lambda}{\pi} \ea \int_{\Oma} \d u \d \bar{u} \del_{\bar{u}} \tilf (u) \del_u \msc(u) \right)  \notag\\
+ &  \msc(z) \left( - \frac{ \i \lambda}{\pi}  \int_{\Oma} \d u \d \bar{u} \del_{\bar{u}} \tilf (u) \del_u \msc(u) \right)^2 
+  \msc(z)  \left( \frac{ \i \lambda}{\pi}  \int_{\Oma} \d u \d \bar{u} \del_{\bar{u}} \tilf (u) \del_u \msc(u)^2 \right) \bigg\}  \notag\\
+ &\Osd ( N^{-1} \eta^{-1} (1 + |\lambda| ) + ( 1 + |\lambda| )^3 N^{-1} +N^{1/2} \Psi(z)^3 + N^{1/2} |\lambda \Phi | \Psi(z)^2) \notag \\
+ &  \Osd ( N^{1/2} |\lambda \Phi^2| ( 1+ |\lambda|) \Psi (z)  )  \label{eqn:intermediate-characteristic}
\end{align}
as well as,
\begin{align} \label{eqn:intermediate-expectation}
& z   \sum_a \ee[ G_{aa} (z) ]  +N= \sum_{ja} \frac{ s_{aj}}{N} \ee[ \del_{aj} G_{ja} (z) ]  \notag\\
-&   \sum_{j \neq a } \frac{ s_{aj}^{(4)}}{N^2} \msc(z)^4 + \sum_{a} \frac{ s_{aa}^{(3)}}{2 N^{3/2}} \msc(z)^3 + \Osd ( N^{-1} \eta^{-1} + N^{1/2} \Psi(z)^3 ) .
\end{align}

\eec
\proof The first estimate follows from applying \eqref{eqn:intermediate-expansion} twice (once with $\lambda =0$) and subtracting the two results. The second just follows from taking $\lambda =0$ in \eqref{eqn:intermediate-expansion}. \qed

\subsection{Better estimate for $G_{11}$} \label{sec:g11-better}

In this section we derive an estimate for $\ee[ \ea (G_{11} (z) - \msc (z) )]$ one order better than the naive estimate $\Psi(z)$. We start with the following.
\bep \label{prop:g11-expand-1}
We have that,
\begin{align} 
& z \ee[ \ea G_{11} (z) ] + \ee[ \ea ] = \sum_{j} \frac{ s_{1j}}{N} \ee[ \del_{j1} ( \ea G_{j1} (z) ) ] \notag\\
&+ \Osd ( N^{-1} (1 + |\lambda| ) + N^{-1/2} \Psi(z)^2 + N^{-1/2} | \lambda \Phi | \Psi(z) ) \label{eqn:g11-expand-1} 
\end{align}
\eep
\proof This follows from Lemma \ref{lem:prelim-expand-1} and then using Lemmas \ref{lem:prelim-fourth-1},   \ref{lem:prelim-third-1}, \ref{lem:prelim-third-2}, \ref{lem:prelim-third-3} and \ref{lem:prelim-third-4} to estimate the terms on the RHS of \eqref{eqn:prelim-expand-1}. Specifically, the errors in the $\Osd$ on the last line of \eqref{eqn:prelim-expand-1} are estimated by $\Osd ( (1 + | \lambda| ) N^{-1} )$. Then, all terms on the RHS of \eqref{eqn:prelim-fourth-1} and \eqref{eqn:prelim-third-4} can be estimated by $\Osd ( (1 + | \lambda| ) N^{-1} )$.  All terms on the RHS of \eqref{eqn:prelim-third-1} can be estimated by $\Osd ( N^{-1/2} \Psi(z)^2)$. All terms on the RHS of \eqref{eqn:prelim-third-2} and \eqref{eqn:prelim-third-3} can be estimated by $\Osd( |\lambda \Phi| N^{-1/2} \Psi (z) )$. \qed

We now compute the second order terms on the RHS of \eqref{eqn:g11-expand-1} in the following lemma.
\bel \label{lem:g11-second-1}
We have,
\beq \label{eqn:g11-expand-2}
\sum_{j} \frac{ s_{1j}}{N} \ee[ G_{j1} \del_{j1} \ea ] = |\lambda| \Osd (   \Phi \Psi(z)( \Psi(z)+ \Phi (1 + |\lambda| ) ) + \min\{ \Psi_1(z)^2 , \Phi^2  \} )
\eeq
and
\beq \label{eqn:g11-expand-3}
\sum_j \frac{s_{1j}}{N} \ee[ \ea \del_{j1} G_{j1}(z) ] = - \msc(z) \ee[ \ea G_{11} (z) ] + \Osd ( N^{-1} \eta^{-1} )
\eeq
\eel
\proof We start with \eqref{eqn:g11-expand-2}. The $j=1$ term contributes $\Osd (  |\lambda|  N^{-1} )$ using \eqref{eqn:del-gen-2}. The other terms we may write as,
\begin{align}
\ee[ G_{1j}(z) \del_{1j} \ea ]= - \frac{2 \i \lambda}{\pi} \int_{\Oma} ( \del_{ \bar{u}} \tilf (u) ) \del_u \ee[ \ea G_{1j} (z) G_{1j} (u) ] \d  u \d \bar{u} .
\end{align}
We now claim that the estimate  \eqref{eqn:loops-expand-4}, applied in the case that $v_j = \ee[ \ea G_{1j} (z) G_{j1} (u)]$, implies, 
\beq \label{eqn:g11-expand-2a}
\left|  \ee[ \ea G_{1j} (z) G_{1j} (u) ] \right| \sd \Psi(z)\Psi(u) ( \Psi(z) + \Psi (u) + |\lambda \Phi| ) + \frac{1}{N ( | \Im[z]| + | \Im[u]| )} .
\eeq
Indeed, in this case the terms involving $X_k$ on the RHS of \eqref{eqn:loops-expand-4} are $\Osd ( \Psi(z) \Psi(u) ( \Psi(u) + | \lambda \Phi | ) )$ (i.e., they involve hypergraphs with three edges, with the additional edge either having spectral parameter $u$ or a new integration variable), the first error term on the RHS of \eqref{eqn:loops-expand-4} is $\Osd ( \Psi(z)^2 \Psi(u))$, and the error term on the second line of \eqref{eqn:loops-expand-4} is $\Osd ( N^{-1} ( |\Im[z]| + |\Im[u] | )^{-1})$. Futhermore, the terms involving $Y_k$ and $Z_k$ are $\Osd ( N^{-1} ( |\Im[z]| + | \Im[u]| )^{-1} )$ by their definition and \eqref{eqn:msc-difference}.

Now, the the estimate \eqref{eqn:g11-expand-2} follows from \eqref{eqn:g11-expand-2a} by applying \eqref{eqn:H-est}. We now turn to \eqref{eqn:g11-expand-3}.  Since $\del_{1j} G_{1j} = - G_{11} G_{jj} + \Osd ( \Psi (z)^2)$ we have,
\begin{align}
 &  \sum_{j \neq 1} \frac{ s_{1j}}{N} \ee[ \ea \del_{1j} G_{j1}(z) ] = - \msc(z) \ee[ \ea G_{11} (z)] \notag\\
 -  & \sum_{j \neq 1 } \frac{ s_{1j}}{N} \ee[ \ea G_{11} (z) (G_{jj} (z) - \msc (z) ) ] + \Osd( \Psi(z)^2) .
\end{align}
The first term on the second line is $\Osd ( N^{-1} \eta^{-1} )$ by fluctuation averaging \eqref{eqn:fluct-av}, which completes the proof. \qed

From \eqref{eqn:g11-expand-1}, \eqref{eqn:g11-expand-2} and \eqref{eqn:g11-expand-3} we immediately conclude the following.
\bep \label{prop:g11-est}
We have,
\beq
\ee[ \ea (G_{11} (z) - \msc (z) ) ] = \Osd ( N^{-1} \eta^{-1} + |\lambda| [ \Phi^2(1+ |\lambda| ) \Psi (z)  + \min \{ (N \eta)^{-1}, \Phi^2 \} ])
\eeq
\eep

Finally we have the following.

\bec \label{cor:weak-g12g21g33}
Let $z, w, v \in \Oma$. We have,
\begin{align}
&| \ee[ \ea G_{12} (z) G_{21} (w) (G_{33} (v) - \msc (v ) ) ] | \sd \Psi_1 (z) \Psi_1 (w) \Psi_1 (v) \min\{ \Psi_1 (z), \Psi_1 (w) \} \notag\\
&\frac{1}{N( | \Im[z] | + | \Im[w] | ) } \left( N^{-1}|\Im[v]| ^{-1} + |\lambda| [ \Phi^2(1+ |\lambda| ) \Psi (v)  + \min \{ (N |\Im[v]| )^{-1}, \Phi^2 \} ] \right)
\end{align}

\eec
\proof This follows immediately from  using Proposition \ref{prop:g11-est} to bound the first term on the second line of \eqref{eqn:weak-g12g21g33} of Lemma \ref{lem:weak-g12g21g33}. \qed

\section{Estimates for $G_{12} G_{21}$ } \label{sec:g12g21}

In this section we find a better an estimate for $\ee[ \ea ( \lambda) G_{12} (z) G_{21} (w ) ]$. Throughout this section we assume that $f$ is admissible and define $z = E \pm \i \eta$ and $w = x \pm \i y$, for $\eta, y >0$. We assume that $z, w \in \Oma$ and that $| \lambda \Phi| \leq N^{-c}$ for some $c>0$. We begin with a cumulant expansion.

\bel \label{lem:g12g21-expand}
We have,
\begin{align} \label{eqn:g12g21-expand-1}
&z \ee[ \ea G_{12} (z) G_{21} (w) ] = \sum_j \frac{s_{1j}}{N} \ee[ \del_{j1} (G_{j2} (z) G_{21} (w) \ea ) ] \notag\\ 
+ & \Osd ( N^{-1/2} \Psi(z) \Psi(w) (\Psi(z) + \Psi (w) ) + |\lambda \Phi | N^{-1/2} \Psi(z) \Psi (w) )
\end{align}
\eel
\proof From \eqref{eqn:loops-expand-1} and \eqref{eqn:loops-third-1} we have,
\begin{align}
&z \ee[ \ea G_{12} (z) G_{21} (w) ] = \sum_j \frac{s_{1j}}{N} \ee[ \del_{1j} (G_{j2} (z) G_{21} (w) \ea ) ] \notag\\
+& \sum_{j \neq 1 } \frac{ s_{1j}^{(3)}}{2 N^{3/2}} \ee[  \del_{1j}^2 ( \ea G_{j2} (z) G_{21} (w) ) ] + \Osd ( N^{-1} (1 + |\lambda| ) \Psi(z) \Psi(w) ) \label{eqn:g12g21-third-1b}
\end{align}
It therefore remains to estimate the sum on the second line of \eqref{eqn:g12g21-third-1b}. 
For the $j=2$ term in the sum we have, 
\beq \label{eqn:g12g21-third-1}
| N^{-3/2} \ee[ \del_{12}^2 ( \ea G_{22} (z) G_{12} (w) ) ] | \sd (1 + |\lambda| )N^{-3/2} \Psi (w) + N^{-3/2} |\lambda \Phi| + N^{-3/2} ( \Psi(z) + \Psi (w) )
\eeq
by direct computation, \eqref{eqn:deleaij}, \eqref{eqn:del-gen-1} and \eqref{eqn:entry-wise}.

By direct calculation and \eqref{eqn:entry-wise} we have for $j \neq 1, 2$,
\begin{align}
\del_{1j}^2 G_{j2} (z) = 2 \msc(z)^2 G_{j2} (z) + \Osd ( \Psi (z)^2), \quad \del_{1j}^2 G_{21} (w) = 2\msc(w)^2 G_{12} (w) + \Osd ( \Psi(w)^2 ) \notag\\
\del_{1j} G_{j2} (z) = - G_{12} (z) \msc(z) + \Osd ( \Psi(z)^2), \quad \del_{1j} G_{21} (w) = - \msc(w) G_{j2} (w) + \Osd ( \Psi(w)^2 ). \label{eqn:g12g21-third-1a}
\end{align}
From the above we find,
\begin{align}
&\sum_{j \neq 1 , 2} \frac{ s_{1j}^{(3)}}{N^{3/2}} \del_{1j}^2 (G_{j2} (z) G_{21} (w) ) = \Osd ( N^{-1/2} \Psi(z)\Psi(w) ( \Psi(z) + \Psi (w) ) \notag\\
+&  \sum_{j \neq 1 , 2} \frac{ 2 s_{1j}^{(3)}}{N^{3/2}} ( (\msc(z)^2 + \msc(w)^2)G_{j2} (z) G_{12} (w) + \msc(z) \msc(w) G_{12} (z) G_{j2} (w)  ) \notag\\
= &\Osd ( N^{-1/2} \Psi(z)\Psi(w) ( \Psi(z) + \Psi (w) )
\end{align}
with the last line following from the isotropic local law \eqref{eqn:iso}.  

Using the estimates on the second line of \eqref{eqn:g12g21-third-1a} and \eqref{eqn:deleaij} we find,
\begin{align}
\sum_{j \neq 1, 2} \frac{ s_{1j}^{(3)}}{N^{3/2}} \ee[(  \del_{1j} \ea) ( \del_{1j} (G_{j2} (z) G_{21} (w) ) ) ] = \Osd ( |\lambda \Phi | N^{-1/2} \Psi(z) \Psi (w) ) .
\end{align}
Using \eqref{eqn:deleaij2} we find,
\begin{align}
\sum_{j \neq 1, 2} \frac{ s_{1j}^{(3)}}{N^{3/2}} \ee[(  \del_{1j}^2 \ea)  (G_{j2} (z) G_{21}(w) )  ] = \Osd ( N^{-1/2} |\lambda \Phi | \Psi (z) \Psi (w) )
\end{align}
using the isotropic local law \eqref{eqn:iso}. This completes the estimation of the sum on the second line of \eqref{eqn:g12g21-third-1b} and the proof. \qed

We now estimate terms appearing on the RHS of the expansion \eqref{eqn:g12g21-expand-1}.

\bep \label{prop:g12g21-second}
We have the estimates,
\beq \label{eqn:g12g21-second-1}
\frac{1}{N} | \ee[ \del_{11} ( \ea G_{12} (z) G_{21} (w) ) ] | \sd (1 + |\lambda| )N^{-1} \Psi(z) \Psi(w),
\eeq
and 
\begin{align} \label{eqn:g12g21-second-2}
& \frac{s_{12}}{N}   \ee[ \del_{21} ( \ea G_{22} (z) G_{21} (w) ) ]  =  - \frac{s_{12} \msc(z) \msc(w)^2}{N} + \Osd( N^{-2} (\eta^{-1} + y^{-1} ) ) \notag\\
+ & \Osd( |\lambda \Phi | N^{-1} \Psi(w) ) +  |\lambda| N^{-1} \Osd (   \Phi^2 ( 1+ |\lambda| ) \Psi(z)  + \min\{ N^{-1} y^{-1} + N^{-1} \eta^{-1}, \Phi^2 \}  )
\end{align}
and for $j \neq 1, 2$,
\beq \label{eqn:g12g21-second-3}
| \ee[ G_{j2} (z) G_{21} (w) \del_{1j} \ea ] | \sd |\lambda|  \Psi_1(z) \Psi_1 (w) \Phi \min\{ \Psi_1(z),\Psi_1 (w), \Phi \} .
\eeq
Finally,
\begin{align} \label{eqn:g12g21-second-4}
&\sum_{ j \neq 1, 2} \frac{ s_{1j}}{N} \ee[ \ea \del_{j1} (G_{j2}(z) G_{21} (w) ) ] = - \msc(z) \ee[ \ea G_{12} (z) G_{21} (w) ] \notag\\
-&  \sum_{j \neq  2} \frac{ s_{1j}}{N} \msc(w) \ee[ \ea  G_{j2} (z) G_{2j} (w) ] + \Osd (( \eta + y )^{-1}  N^{-2} (y^{-1} + \eta^{-1}  )) \notag\\
+ & \frac{|\lambda| }{N(\eta+y)} \Osd (   \Phi^2 (1 + |\lambda| ) ( \Psi (z) + \Psi (w) ) + \min \{ N^{-1} (y^{-1} + \eta^{-1} ), \Phi^2 \}  )
\end{align}
\eep
\proof The estimate \eqref{eqn:g12g21-second-1} follows easily from \eqref{eqn:del-gen-2} and \eqref{eqn:entry-wise}. For the estimate \eqref{eqn:g12g21-second-2}, the term where the derivative hits $\ea$ contributes $\Osd ( |\lambda \Phi | N^{-1} \Psi (w) )$ by \eqref{eqn:deleaij}. The term where the derivative hits $G_{22} (z)$ contributes $\Osd( N^{-1} \Psi(z) \Psi(w))$. Finally,
\begin{align}
& \frac{s_{12}}{N}   \ee[  \ea G_{22} (z) \del_{21} (G_{21} (w) ) ] = -\frac{s_{12}}{N} \ee[ \ea G_{22} (z) G_{11} (w) G_{22} (w) ] + \Osd ( N^{-1} \Psi(w)^2 ) \notag\\
= & - \frac{s_{12} \msc(z) \msc(w)^2}{N} + \Osd( N^{-2} (\eta^{-1} + y^{-1} ) ) \notag \\
+ & N^{-1} |\lambda| \Osd (  \Phi^2 ( 1+ |\lambda| ) (\Psi(z) + \Psi(w)) + \min\{ N^{-1} y^{-1} + N^{-1} \eta^{-1}, \Phi^2 \}  )
\end{align}
The second estimate follows from Proposition \ref{prop:g11-est} after linearizing each of the $G_{ii}$ around $\msc$. 

Turning now to \eqref{eqn:g12g21-second-3} we have by the equality in \eqref{eqn:deleaij}, 
\beq
 \ee[ G_{j2} (z) G_{21} (w) \del_{j1} \ea ] = -\frac{ 2 \i \lambda}{\pi} \int_{\Oma} ( \del_{\bar{u}} \tilf (u) ) \del_u \ee[ \ea G_{j2} (z) G_{21} (w) G_{1j} (u ) ] \d u \d \bar{u}.
\eeq
From Proposition \ref{prop:main-loop-estimate} we have,
\beq
|  \ee[ \ea G_{j2} (z) G_{21} (w) G_{1j} (u ) ]  | \leq \Psi_1 (z) \Psi_1 (w) \Psi_1 (u ) \min \{ \Psi_1 (z), \Psi_1 (w) , \Psi_1 (u ) \} 
\eeq
and so \eqref{eqn:g12g21-second-3} follows, using \eqref{eqn:H-est}. 

Turning now to \eqref{eqn:g12g21-second-4} we have by direct computation,
\begin{align}
& \sum_{ j \neq 1, 2} \frac{ s_{1j}}{N} \ee[ \ea \del_{j1} (G_{j2}(z) G_{21} (w) ) ] = - \msc(z) \ee[ \ea G_{12} (z) G_{21} (w) ] + \Osd(N^{-1} \Psi(w) \Psi (z) ) \notag\\
-&  \sum_{j \neq 2} \frac{ s_{1j}}{N} \msc(w) \ee[ \ea  G_{j2} (z) G_{2j} (w) ]  - \sum_{j \neq 1, 2} \frac{ s_{1j}}{N} \ee[ \ea G_{12} (z) G_{21} (w) (G_{jj} (z) - \msc(z) )] \notag\\
- & \sum_{j \neq 1, 2} \frac{ s_{1j}}{N} \ee[ \ea G_{j2} (z) G_{2j} (w) (G_{11} (w) - \msc(w) )] \notag\\
- & \sum_{j \neq 1, 2} \frac{s_{1j}}{N} \ee[ \ea (G_{1j} (z) G_{j2} (z) G_{21} (w) + G_{j2}(z) G_{21} (w) G_{j1} (w) ) ]  \label{eqn:g12g21-second-4a}
\end{align}
The claimed estimates now follow from Corollary \ref{cor:weak-g12g21g33} applied to the term on the third line and second term on the second line of \eqref{eqn:g12g21-second-4a}, and Proposition \ref{prop:main-loop-estimate} applied to the terms on the last line. \qed

Lemma \ref{lem:g12g21-expand} and Proposition \ref{prop:g12g21-second} easily imply the following and so the proof is omitted.
\bec \label{cor:g12g21-self}
We have that,
\begin{align}
 & z \ee[ \ea G_{12} (z) G_{21} (w) ] = - \msc(z) \ee[ \ea G_{12} (z) G_{21} (w) ] - \sum_{j \neq  2} \frac{ s_{1j}}{N} \ee[ \ea G_{j2} (z) G_{j2} (w) ] \notag\\
- & \frac{ s_{12} \msc(z) \msc(w)^2}{N} +  \Osd ( (N y)^{-1} (N \eta)^{-1} ) \notag\\
+ & \frac{|\lambda|}{N (y + \eta)} \Osd ( \Phi^2 (1 + |\lambda| ) ( \Psi (z) + \Psi (w) ) + \min \{ N^{-1} (y^{-1} + \eta^{-1} ), \Phi^2 \}  ) \notag \\
+ & \Osd ( |\lambda| \Psi_1 (z) \Psi_1 (w) \Phi \min \{ \Psi_1 (z), \Psi_1 (w), \Phi \} )
\end{align}

\eec

We define now,
\beq
v_j := \ee[ \ea G_{j2} (z) G_{2j} (w) ] \1_{ \{ j \neq 2 \} } .
\eeq
Recall the definition of $A$ in Lemma \ref{lem:spectral-gap} and let $B$ be the matrix $B_{ij} = \frac{S_{ij}}{N} \1_{ \{ i=2 \} }$. 
Corollary \ref{cor:g12g21-self} implies that we have the equation,
\beq \label{eqn:g12g21-self-1}
(1 - \msc(z) \msc(w) (A - B) ) v = \msc(z) \msc(w) N^{-1} \be \be^* v + X + \eps
\eeq
where $X_j := s_{j2} \1_{ \{ j \neq 2 \} } \frac{ \msc(z)^2 \msc(w)^2}{ N }$ and the vector $\eps$ obeys,
\begin{align}
\| \eps \|_\infty \sd  &(N y)^{-1} (N \eta)^{-1} ) \notag\\
+ & N^{-1} (y + \eta )^{-1} (|\lambda| [(  \Phi^2 (1 + |\lambda| ) ( \Psi (z) + \Psi (w) ) + \min \{ N^{-1} (y^{-1} + \eta^{-1} ), \Phi^2 \}  ) \notag \\
+ & ( |\lambda| \Psi_1 (z) \Psi_1 (w) \Phi \min \{ \Psi_1 (z), \Psi_1 (w), \Phi \} )
\end{align}
We now compute the inner product on the RHS of \eqref{eqn:g12g21-self-1}. 
\bel \label{lem:g12g21-dif}
We have,
\begin{align}
& N^{-1} \be \be^* v = \ee[ \ea] \left( \frac{ \msc(z) - \msc(w)}{N (z-w) } - \frac{ \msc(z) \msc(w)}{N} \right) + N^{-1} (y + \eta)^{-1} \Osd( N^{-1} (y^{-1} + \eta^{-1} ) ) \notag\\
&+ N^{-1} (y + \eta)^{-1} \Osd (|\lambda| [(  \Phi^2(1+ |\lambda| ) (\Psi (z)+ \Psi(w))  + \min \{ N^{-1} ( y^{-1} + \eta^{-1} ), \Phi^2 \} ])
\end{align}
\eel
\proof We write,
\begin{align}
N^{-1} \be \be^* v = \frac{1}{N} \sum_{j} \ee[ \ea G_{j2} (z) G_{2j} (w) ] - \frac{1}{N} \ee[ \ea G_{22}(z) G_{22} (w) ].
\end{align}
For the second term we write,
\begin{align}
 & \frac{1}{N} \ee[ \ea G_{22}(z) G_{22} (w) ] = \frac{ \msc(z) \msc(w)}{N}\ee[ \ea]  \notag\\
+ & N^{-1}\Osd ( N^{-1}(y^{-1} + \eta^{-1}) + |\lambda| [ \Phi^2(1+ |\lambda| ) (\Psi (z)+ \Psi(w))  + \min \{ N^{-1} ( y^{-1} + \eta^{-1} ), \Phi^2 \} ])
\end{align}
using Proposition \ref{prop:g11-est}. For the first term we write,
\begin{align}
 & \frac{1}{N} \sum_{j} \ee[ \ea G_{j2} (z) G_{2j} (w) ] = \frac{\msc(z) - \msc(w)}{N(z-w)} \ee[ \ea] \notag\\
   + &  \frac{1}{N(z-w)} \ee[ \ea (G_{22} (z) - \msc(z) + \msc(w) - G_{22} (w) ) ] .
\end{align}
We now claim that,
\begin{align} \label{eqn:g12g21-dif-1}
& \left|   \frac{1}{N(z-w)} \ee[ \ea (G_{22} (z) - \msc(z) + \msc(w) - G_{22} (w) ) ]  \right| \notag\\
\sd & \frac{1}{N ( y + \eta )} \left( N^{-1} ( y^{-1} + \eta^{-1} ) + |\lambda| [\Phi^2(1 + |\lambda| ) ( \Psi(z) + \Psi (w) ) + \min\{ N^{-1} ( y^{-1} + \eta^{-1}  ) , \Phi^2 \} ]\right)
\end{align}
which will complete the proof. In order to prove \eqref{eqn:g12g21-dif-1}, WLOG assume $\eta > y$. If $|z-w| > \eta/2$ then we can estimate the prefactor $|N(z-w)|^{-1} \prec (N \eta)^{-1}$ and the expectation in \eqref{eqn:g12g21-dif-1} directly by Proposition \ref{prop:g11-est}, yielding the claim.  Otherwise, assume that $|z-w| < \eta/2$. Then $y \asymp \eta$ and $\Psi(z) \asymp \Psi (w)$. We then write,
\begin{align}
 &\frac{1}{N(z-w)} \ee[ \ea (G_{22} (z) - \msc(z) + \msc(w) - G_{22} (w) ) ]  \notag\\
 = & N^{-1} \int_{0}^1 \del_s \ee[ \ea (  G_{22} (w + s(z-w)) - \msc(w + s(z-w) ) ) ] \d s.
\end{align} 
We conclude \eqref{eqn:g12g21-dif-1} by estimating the integrand using Proposition \ref{prop:g11-est} and the Cauchy Integral formula. \qed

Now for every $k$ let us define the vector
\beq
v_j^{(k)} := \ee[ \ea G_{jk} (z) G_{kj} (w) ] \1_{ \{ j \neq k \} }. 
\eeq
 Define $M = \msc(z) \msc (w)$ as well as,
\beq
F :=M \frac{ \msc(z) - \msc(w)}{N(z-w)} - \frac{ M^2}{N}.
\eeq
 Define the matrix $B^{(k)}$ by $B^{(k)}_{ab} = \frac{s_{ab}}{N} \1_{ \{ a = k \} }$ and the vector,
 \beq
 Y_j^{(k)} := \frac{ s_{jk}}{N} M^2 \1_{ \{ j \neq k \} } 
 \eeq

We have the following equation for the $v_j^{(k)}$. 
\bec
We have that,
\begin{align} \label{eqn:g12g21-self-2}
v^{(k)}_j &= \ee[\ea]  \sum_i (1 - \msc (z) \msc(w) (A-B^{(k)}) )^{-1}_{ji} \left\{ Y^{(k)}_i + F   \right\} \notag\\
+  & \Osd ( (N y)^{-1} (N \eta)^{-1} ) \notag\\
+ & \frac{|\lambda|}{N ( y+ \eta ) } \Osd (  \Phi^2 (1 + |\lambda| ) ( \Psi (z) + \Psi (w) ) + \min \{ N^{-1} (y^{-1} + \eta^{-1} ), \Phi^2 \}  ) \notag \\
+ & \Osd ( |\lambda| \Psi_1 (z) \Psi_1 (w) \Phi \min \{ \Psi_1 (z), \Psi_1 (w), \Phi \} )
\end{align}

\eec
\proof This follows by replacing the index $2$ by $k$ in \eqref{eqn:g12g21-self-1} and Lemma \ref{lem:g12g21-dif}, after noticing that Lemma \ref{lem:inverse} applies to bound the matrix elements of $(1 - M (A - B^{(k)} ) )^{-1}$. \qed

For future use, we now compute the deterministic term that arises when we apply \eqref{eqn:g12g21-self-2} to compute $\sum_{ij} \frac{s_{ij}}{N} v_j^{(i)}$. 
\bel
We have that,
\begin{align} \label{eqn:g12g21-determ}
& \frac{1}{N} \sum_{ij} s_{ij} [ (1 - M(A-B^{(i)}) )^{-1} ( Y^{(i)} + \be F ) ]_{j} = - M \tr S+  M \tr ( S (1- M S)^{-1} )\notag\\
 +&  \Osd ( N^{-1} (y + \eta)^{-1} )
\end{align}

\eel
\proof We first show that we can remove the operator $B^{(i)}$ up to a negligible error.  We have,
\begin{align}
& \sum_j s_{ij} [ (1- M (A-B^{(i)} )^{-1} (Y^{(i)} + \be F) ]_j - \sum_j s_{ij} [ (1- M A )^{-1} (Y^{(i)} + \be F) ]_j \notag \\
= & \left( (1- M(A - B^{(i)} ))^{-1} s^{(i)} \right) \cdot \left( B^{(i)} (1- MA)^{-1} (Y^{(i)} + \be F ) \right)
\end{align}
where $s^{(i)}$ is the vector with entries $s_{ij}$. 
Now, since $(1- M (A - B^{(i)} ))^{-1}$ and $(1- MA )^{-1}$ are bounded on $\ell^\infty$, the entries of the vector $ (1- M(A - B^{(i)} ))^{-1} s^{(i)}$ are bounded by a constant, and the entries $(1- MA)^{-1} (Y^{(i)} + \be F )$ are  $\O ( N^{-1} ( \eta + y )^{-1})$ (see \eqref{eqn:msc-difference}), and so
\begin{align}
 & \left|  \left( (1- M(A - B^{(i)} ))^{-1} s^{(i)} \right) \cdot \left( B^{(i)} (1- MA)^{-1} (Y^{(i)} + \be F ) \right)  \right|  \notag\\
 \sd & \frac{1}{N(y + \eta ) } \sum_{ab} |B_{ab}^{(i)} | \sd \frac{1}{N^2(y + \eta ) } \sum_{ab} \1_{ \{ a = i \} } \prec \frac{1}{ N ( y + \eta ) }. 
\end{align}
Therefore, up to an error of $\Osd ( N^{-1} (y + \eta)^{-1} )$ it suffices to compute,
\beq
\frac{1}{N} \sum_{ij} s_{ij} [(1 - M A )^{-1} ( Y^{(i)} + \be F ) ]_j .
\eeq
The difference between the vector $Y^{(i)}$ and the vector $M^2 s^{(i)} N^{-1}$ is the vector with a single non-zero entry $\1_{ \{ j  =i \} } M^2 s_{ii} / N$. Therefore,
\begin{align}
& \left| \frac{1}{N} \sum_{ij} s_{ij} [( 1 - M A)^{-1} ( Y^{(i)} - M^2 s^{(i)} N^{-1} ) ]_j \right| \notag\\
\sd & \frac{1}{N^2} \sum_{ijk}  |(1- M A)^{-1}_{jk}| \1_{ \{ k = i \} } = \frac{1}{N^2} \sum_{j,k} |(1 - MA )^{-1}_{jk} | \sd \frac{1}{N}
\end{align}
with the last estimate following from the fact that $(1- MA )^{-1}$ is bounded from $\ell^\infty$ to $\ell^\infty$. We have therefore reduced the proof to computing,
\begin{align} \label{eqn:deterministic-1}
\frac{1}{N} \sum_{ij} s_{ij} (1 - M A)^{-1}_{jk} s_{ik} M^2 N^{-1} + \frac{1}{N} \sum_{ij} s_{ij} ( (1 - MA )^{-1} \be F )_j
\end{align}
Since $(1 - M A)^{-1} \be = \be$ by the definition of $A$ we have,
\beq \label{eqn:deterministic-2}
 \frac{1}{N} \sum_{ij} s_{ij} ( (1 - MA )^{-1} \be F )_j = \frac{F}{N} \sum_{ij} s_{ij} = M \frac{ \msc(z) - \msc(w)}{z-w} - M^2 = M \frac{ M}{1 - M} - M^2
\eeq
with the last equality using the identity \eqref{eqn:msc-dif-identity}.  For the other term in \eqref{eqn:deterministic-1} we have,
\begin{align}
\frac{1}{N} \sum_{ij} s_{ij} (1 - M A)^{-1}_{jk} s_{ik} M^2 N^{-1} = M^2 \tr ( S(1- M A)^{-1} S).
\end{align}
On the other hand via the Sherman Morrison formula,
\beq
\frac{1}{ 1 - MS} = \frac{1}{ 1- M A } + M N^{-1} \frac{ \be \be^*}{ 1- \msc(z) \msc (w) }
\eeq
so
\begin{align}
& \tr ( S ( 1- M A)^{-1} M^2 S ) = \tr (S (1- M S )^{-1} M^2 S ) - M^3 \frac{1}{ 1 - M} \notag \\
= & -M \tr S + M \tr (S (1- M S)^{-1} ) - M^3 / (1- M)
\end{align}
Therefore, using the above and \eqref{eqn:deterministic-2} we find that the expression in \eqref{eqn:deterministic-1} equals,
\begin{align}
& \frac{1}{N} \sum_{ij} s_{ij} (1 - M A)^{-1}_{jk} s_{ik} M^2 N^{-1} + \frac{1}{N} \sum_{ij} s_{ij} ( (1 - MA )^{-1} \be F )_j \notag\\
= & - M \tr S + M \tr (S (1- MS)^{-1} ) - M^3 / (1- M) + M \frac{ M}{1-M} - M^2 \notag\\
= & - M \tr S+  M \tr ( S (1- M S)^{-1} )
\end{align}
which yields the claim. \qed

\bec \label{cor:gij-main-est}
For $| \lambda \Phi | \leq 1$ we have,
\begin{align}
& \frac{1}{N} \sum_{i \neq j } s_{ij} \ee[ \ea G_{ij} (z) G_{ji} (w) ] = \ee[\ea] \left( - M \tr S + M \tr (S (1- MS)^{-1} ) \right)\notag\\
+ &   \Osd ( (N y)^{-1} ( \eta)^{-1} ) \notag\\
+ &  \frac{|\lambda|}{y + \eta} \Osd (  \Phi^2 (1 + |\lambda| ) ( \Psi (z) + \Psi (w) ) + \min \{ N^{-1} (y^{-1} + \eta^{-1} ), \Phi^2 \}  ) \notag \\
+ & N \Osd ( |\lambda| \Psi_1 (z) \Psi_1 (w) \Phi \min \{ \Psi_1 (z), \Psi_1 (w), \Phi \} )
\end{align}

\eec
\proof We have,
\beq
 \frac{1}{N} \sum_{i \neq j } s_{ij} \ee[ \ea G_{ij} (z) G_{ji} (w) ] =  \frac{1}{N} \sum_{i \neq j } s_{ij} v^{(i)}_j =  \frac{1}{N} \sum_{i , j } s_{ij} v^{(i)}_j .
\eeq
We then apply \eqref{eqn:g12g21-self-2} and \eqref{eqn:g12g21-determ}. \qed

\section{Estimates for $G_{11} G_{22}$} \label{sec:g11g22}

Let $b_a \in \cc$ for $1 \leq a \leq N$ satisfy $|b_a| \leq C$ some fixed $C>0$ and define,
\beq \label{eqn:Xi-est}
X_i := \sum_{a \neq i } \frac{b_a}{N} ( G_{aa} (w) - \msc (w) ) = \Osd ( (N y)^{-1} )
\eeq
where the estimate follows by fluctuation averaging \eqref{eqn:fluct-av}, where we denote $w = x \pm \i y$, with $y>0$. In this section we are going to find a good estimate for the quantity $\ee[ \ea G_{ii} (z) X_i ]$. We again denote $z = E \pm \i \eta$ with $\eta >0$. We will assume throughout this section that $z, w \in \Oma$. The following performs a cumulant expansion for this quantity.

\bel \label{lem:giigjj-fourth} For any $j$ and $k \geq 0$,
\beq \label{eqn:giigjj-fourth-1}
\del_{j1}^k X_1 = \Osd ( \Psi(w)^2 + \1_{\{ k = 0\}} (N y)^{-1})
\eeq
For $|\lambda \Phi | \leq 1$ we have,
\beq \label{eqn:giigjj-fourth-2}
N^{-2} \del_{11}^3 ( G_{11} (z) \ea X_1 ) = \Osd ((1 + |\lambda|)^3 N^{-2} ( N y)^{-1} )
\eeq
and for $j \neq 1$,
\beq \label{eqn:giigjj-fourth-3}
N^{-2} | \del_{j1}^3 ( \ea G_{j1}(z) X_1 ) | \sd N^{-2} ( 1 + |\lambda| ) (N y)^{-1} ) .
\eeq
For any $\delta >0$, the above estimates also hold with an additional $N^\delta$ factor on the RHS for any $\eps$-regular matrix, if $\eps >0$ is sufficiently small depending on $\delta >0$. 

We have for $|\lambda \Phi| \leq 1$,
\begin{align} \label{eqn:giigjj-fourth-4}
z\ee[ \ea G_{11}(z) X_1 ] + \ee[ \ea X_1 ]&= \sum_{j} \frac{ s_{1j}}{N} \ee[ \del_{j1} (G_{j1} (z) \ea X_1 ) ]  + \sum_{j} \frac{ s_{1j}^{(3)}}{2N^{3/2}} \ee[ \del_{j1}^2 (G_{j1} (z) \ea X_1 ) ] \notag\\
&+ \Osd (  N^{-1} ( 1 + |\lambda| ) (N y)^{-1} )
\end{align}
\eel
\proof The estimate \eqref{eqn:giigjj-fourth-1} follows from direct computation and \eqref{eqn:entry-wise}. For \eqref{eqn:giigjj-fourth-2}, if all derivatives hit $\ea$, we can estimate this by $(1 + |\lambda| )^3 N^{-2} (N y)^{-1}$ using \eqref{eqn:Xi-est} and \eqref{eqn:del-gen-2}. Otherwise using \eqref{eqn:giigjj-fourth-1} or \eqref{eqn:Xi-est} we find an estimate of $\Osd ( (1+ | \lambda| )^2 N^{-2} (N y)^{-1} )$, completing the proof of \eqref{eqn:giigjj-fourth-2}. For \eqref{eqn:giigjj-fourth-3} we find an estimate of $\Osd ( (1 + |\lambda| ) N^{-2} \Psi(z) (N y)^{-1} )$ if no derivatives hit $G_{1j}$ using \eqref{eqn:del-gen-2}, \eqref{eqn:Xi-est} and \eqref{eqn:entry-wise}. Otherwise, if a derivative hits $G_{1j}$ we find an estimate of $\Osd ( N^{-2} (1 + |\lambda| ) (N y)^{-1} )$. The statement about $\eps$-regular matrices is clear. The expansion \eqref{eqn:giigjj-fourth-4} now follows from applying the cumulant expansion Lemma \ref{lem:cumu-exp}, similar to, e.g., Lemma \ref{lem:loops-expand-1}. \qed

The following lemma takes care of the third order terms in the cumulant expansion above.
\bel \label{lem:giigjj-third}
Assume $|\lambda \Phi | \leq 1$. We have,
\beq \label{eqn:giigjj-third-1}
N^{-3/2} |\ee[ \del_{11}^2 ( \ea G_{11} (z) X_1 ) ] | \sd N^{-3/2} ( 1 + | \lambda | )^2 (N y)^{-1} 
\eeq
and
\beq \label{eqn:giigjj-third-2}
\sum_{j \neq 1 } \frac{s_{1j}^{(3)}}{N^{3/2}} \ee[ G_{1j}(z) X_1 \del_{j1}^2 \ea ] = \Osd ( |\lambda \Phi| N^{-1/2} \Psi(z) (N y)^{-1} )
\eeq
and
\beq \label{eqn:giigjj-third-3} 
\sum_{j \neq 1 } \frac{ s_{1j}^{(3)}}{N^{3/2}} \ee[ ( \del_{j1} \ea ) \del_{j1} ( G_{j1} X_1 ) ] = \Osd( \Psi(z) (N y)^{-1} N^{-1/2} |\lambda| \Phi )
\eeq
and
\beq \label{eqn:giigjj-third-4}
\sum_{j \neq 1} \frac{ s_{1j}^{(3)}}{N^{3/2}} \ee[ \ea \del_{j1}^2 (G_{j1} X_1 ) ] = \Osd( N^{-1/2} (N y)^{-1} \Psi (z) )
\eeq
\eel
 \proof The estimate \eqref{eqn:giigjj-third-1} is straightforward from \eqref{eqn:giigjj-fourth-1} and \eqref{eqn:del-gen-2}. The estimate \eqref{eqn:giigjj-third-2} follows from the fact that $\del_{1j}^2 \ea = \i \lambda c_f + \Osd ( | \lambda \Phi | )$ (see \eqref{eqn:deleaij2}), the isotropic local law \eqref{eqn:iso} and \eqref{eqn:Xi-est}. For \eqref{eqn:giigjj-third-3}, the term where $\del_{j1}$ hits $X_1$ is easily bounded by $\Osd ( N^{-1/2} |\lambda \Phi| \Psi(z) (N y)^{-1} )$ using \eqref{eqn:giigjj-fourth-1} and \eqref{eqn:deleaij}. For the other term we have since $\del_{j1} G_{j1}(z) = - \msc(z)^2 + \Osd ( \Psi(z) )$,
 \begin{align}
 &\sum_{j \neq 1 } \frac{ s_{1j}^{(3)}}{N^{3/2}} \ee[ ( \del_{1j} \ea ) (\del_{j1}  G_{j1} (z)) X_1  ] = - \msc(z)^2 \sum_{j \neq 1 } \frac{ s_{1j}^{(3)}}{N^{3/2}} \ee[ X_1 ( \del_{j1} \ea ) ] \notag\\
  + & \Osd( \Psi(z) (N y)^{-1} N^{-1/2} |\lambda| \Phi ) =   \Osd( \Psi(z) (N y)^{-1} N^{-1/2} |\lambda| \Phi )
 \end{align}
 with the second estimate following from the isotropic local law \eqref{eqn:iso} applied to the equality in \eqref{eqn:deleaij}. This completes \eqref{eqn:giigjj-third-3}.  For \eqref{eqn:giigjj-third-4}, since $\del_{1j}^2 X_1 = \Osd ( (N y)^{-1})$ and $\del_{j1}^2 G_{1j} = \Osd ( \Psi (z))$ and $\del_{j1} G_{1j} = - \msc(z)^2 + \Osd ( \Psi (z) )$ we have,
 \begin{align} \label{eqn:giigjj-third-4a}
 &  \sum_{j \neq 1} \frac{ s_{1j}^{(3)}}{N^{3/2}} \ee[ \ea \del_{j1}^2 (G_{j1} X_1 ) ] =- 2\msc(z)^2 \sum_{j\neq 1} \frac{ s_{1j}^{(3)}}{N^{3/2}} \ee[ \ea ( \del_{j1} X_1 ) ] 
 +  \Osd ( N^{-1/2} (N y)^{-1} \Psi(z) ).
 \end{align}
Then,
\begin{align} \label{eqn:giigjj-third-4b}
& -\sum_{j \neq 1} \frac{ s_{1j}^{(3)}}{N^{3/2}} \del_{j1} X_1   = 2 \sum_{j \neq 1 } \sum_{a \neq 1} \frac{ b_a  s_{1j}^{(3)} }{N^{5/2}}  G_{aj} (w) G_{a1}(w) = 2\sum_{a \neq 1 } G_{a1} \frac{ b_a}{N^2} \sum_{j \neq 1} \frac{ s_{1j}^{(3)}}{N^{1/2}} G_{aj}   \notag\\
= & 2\sum_{a \neq 1 } G_{a1} \frac{ b_a}{N^{5/2}} s_{1a}^{(3)} \msc(w) + \Osd ( N^{-1} \Psi(w)^2 ) = \Osd( N^{-1} \Psi(w)^2 )  .
\end{align}
Here, the third equality follows from the isotropic law \eqref{eqn:iso}. This completes the proof of \eqref{eqn:giigjj-third-4}. \qed

We now compute the second order terms on the RHS of \eqref{eqn:giigjj-fourth-4}. 
\bel \label{lem:giigjj-second} Assume that there is a $c>0$ so that $|\lambda \Phi| \leq N^{-c}$. Then we have
\beq \label{eqn:giigjj-second-1}
N^{-1} | \ee[ \del_{11} ( \ea G_{11} (z) X_1 ) ] | \sd N^{-1} ( 1 + |\lambda| ) (N y)^{-1} 
\eeq
and
\beq \label{eqn:giigjj-second-2}
\sum_{j \neq 1 } \frac{ s_{1j}}{N} \ee[ G_{j1} (z) X_1 \del_{j1} \ea ] = \Osd ( \Psi(z) |\lambda \Phi | ( N y )^{-1} )
\eeq
and
\beq \label{eqn:giigjj-second-3}
\sum_{j \neq 1 } \frac{ s_{1j}}{N} \ee[ \ea G_{j1}  (z) \del_{j1} X_1 ] = \Osd ( (N y)^{-1} (N \eta)^{-1} )
\eeq
and 
\beq
\label{eqn:giigjj-second-4}
\sum_{j \neq 1 } \frac{ s_{1j}}{N} \ee[ \ea X_1 \del_{j1} G_{j1} ] = - \msc(z) \ee[ G_{11} X_1 \ea] + \Osd ( (N y)^{-1} (N \eta)^{-1} )
\eeq
\eel
\proof The estimates \eqref{eqn:giigjj-second-1} and \eqref{eqn:giigjj-second-2} are straightforward. For \eqref{eqn:giigjj-second-3} we write,
\begin{align}
-&\sum_{j \neq 1 } \frac{ s_{1j}}{N} \ee[ \ea G_{j1} (z)\del_{j1} X_1 ] =2  \sum_{j \neq 1 } \sum_{a \neq 1} \frac{ s_{1j} b_a}{N^2} \ee[ \ea G_{j1}(z) G_{1a}(w) G_{aj}(w) ] \notag\\
= & 2 \sum_{j\neq 1} \sum_{a \neq 1, j} \frac{ s_{1j} b_a}{N^2} \ee[ \ea G_{j1}(z) G_{1a}(w) G_{aj}(w) ] + \Osd( N^{-1} \Psi(z) \Psi (w) )
\end{align}
Now \eqref{eqn:giigjj-second-3} follows since $\ee[ \ea G_{j1}(z) G_{1a} (w) G_{aj}(w) ] = \Osd ( (N y)^{-1} (N \eta )^{-1} )$ for indices $1, a, j$ all distinct, by Proposition \ref{prop:main-loop-estimate}. 

Finally, for \eqref{eqn:giigjj-second-4} we have,
\begin{align}
& \sum_{j \neq 1 } \frac{ s_{1j}}{N} \ee[ \ea X_1 \del_{j1} G_{j1} (z)]  = - \msc(z) \sum_{j \neq 1 } \frac{ s_{1j}}{N} \ee[ \ea X_1 G_{11} (z)] \notag\\
- & \sum_{j \neq 1} \frac{s_{1j}}{N} \ee[ \ea X_1 G_{11} (G_{jj}(z) - \msc(z) ) ] + \Osd (\Psi(z)^2 (N y)^{-1} )
\end{align}
Now the first term on the second line is $\Osd ( (N \eta)^{-1} (N y )^{-1} )$ by fluctuation averaging \eqref{eqn:fluct-av}. \qed

\bep \label{prop:giigaa-est}
Assume that $|\lambda \Phi| \leq N^{-c}$ for some $c>0$. For any $i$ we have,
\begin{align}
\sum_{a} \frac{b_a}{N} \ee[ \ea (G_{ii} (z) - \msc(z) ) (G_{aa} (w) - \msc (w) ) ] = \Osd ( (N y)^{-1} (N \eta )^{-1} + |\lambda \Phi| \Psi(z) (N y)^{-1} )
\end{align}
\eep
\proof WLOG we can assume $i=1$. Then the term in the sum with $a=1$ contributes $\Osd (N^{-1} \Psi(z) \Psi (w) )$ and so can be ignored. The equation \eqref{eqn:giigjj-fourth-4} combined with Lemmas \ref{lem:giigjj-third} and  \ref{lem:giigjj-second}  imply,
\begin{align}
(z + \msc(z) ) \ee[ (G_{11} (z) - \msc (z) ) X_1 ] = \Osd ( (N y)^{-1} (N \eta )^{-1} + |\lambda \Phi| \Psi(z) (N y)^{-1} )
\end{align}
and so the claim follows since $ | z + \msc(z) | \asymp 1$. \qed

\subsection{Expansion for expectation} \label{sec:expectation}

In this section we derive an expansion for the expectation,
\beq
\frac{1}{\pi} \int_{\Oma} ( \del_{\bar{z}} \tilf (z) ) \sum_a \ee[ G_{aa} (z) ] \d z \d \bar{z} .
\eeq
We will use the expansion \eqref{eqn:intermediate-expectation}, and begin by computing the second order terms on its RHS. 
\bel \label{lem:exp-second}
We have
\begin{align} \label{eqn:exp-second-1}
&\sum_{ja} \frac{ s_{ja}}{N} \ee[ \del_{ja} G_{ja} (z) ] = - N \msc(z)^2 - 2 \msc(z) \sum_a \ee[ G_{aa}(z) - \msc(z) ] \notag\\
& + \msc(z)^2 \tr S - \msc(z)^2 \tr (S (1- \msc(z)^2S)^{-1} )+ \Osd (N^{-1} \eta^{-2} ).
\end{align}
and so
\begin{align} \label{eqn:exp-second-2}
& \sum_a \ee[ G_{aa}(z) - \msc(z) ] =  \sum_{j \neq a } \frac{ s_{aj}^{(4)}}{N^2} \msc(z)^3 \msc'(z) - \sum_a \frac{s_{aa}^{(3)}}{2 N^{3/2}} \msc(z)^2 \msc'(z) \notag\\
& - \msc(z) \msc'(z) \tr S + \msc(z) \msc'(z) \tr (S (1- \msc(z)^2 S)^{-1} ) + \Osd (| \msc' (z) | N^{-1} \eta^{-2} ).
\end{align}
\eel
\proof We first have,
\beq \label{eqn:exp-second-1a}
\sum_{ja} \frac{ s_{ja}}{N} \ee[ \del_{ja} G_{ja} (z) ] = -\sum_{ja} \frac{ s_{ja}}{N} \ee[ G_{aa} (z) G_{jj} (z) ] - \sum_{j\neq a} \frac{ s_{ja}}{N} \ee[ G_{ja}^2 (z) ]
\eeq
We then write,
\begin{align}
& \sum_{ja} \frac{ s_{ja}}{N} \ee[ G_{aa} G_{jj} (z) ] = N \msc(z)^2 + 2 \msc \sum_a \ee[ G_{aa}(z) - \msc(z) ] \notag\\
+ & \sum_{ja} \frac{s_{ja}}{N} \ee[ (G_{aa} (z) - \msc(z) )( G_{jj} (z) - \msc (z) ) ] 
\end{align}
Now by Proposition \ref{prop:giigaa-est} the second line is $\Osd ( N^{-1} \eta^{-2} )$. By Corollary \ref{cor:gij-main-est} we see that
\begin{align}
\sum_{j\neq a} \frac{ s_{ja}}{N} \ee[ G_{ja}^2 (z) ] = - \msc(z)^2 \tr S + \msc(z)^2 \tr (S (1- \msc(z)^2S)^{-1} ) + \Osd (N^{-1} \eta^{-2} ).
\end{align}
This completes the proof of \eqref{eqn:exp-second-1}. The equation \eqref{eqn:exp-second-2} follows from substituting the result of \eqref{eqn:exp-second-1} into \eqref{eqn:intermediate-expectation} and using that $\msc'(z) = - \frac{ \msc(z) }{ 2 \msc(z) +z }$. \qed

\bel \label{lem:exp-determ}
Define,
\begin{align}  \label{eqn:exp-determ-0}
 &E(z) :=\sum_{j \neq a } \frac{ s_{aj}^{(4)}}{N^2} \msc(z)^3 \msc'(z) - \sum_a \frac{s_{aa}^{(3)}}{2N^{3/2}} \msc(z)^2 \msc'(z) \notag\\
& - \msc(z) \msc'(z) \tr S + \msc(z) \msc'(z) \tr (S (1- \msc(z)^2 S)^{-1} )
\end{align}
We have,
\beq \label{eqn:exp-determ-1} 
\int_{\Oma} ( \del_{\bar{z}} \tilf (z)) E(z) \d z \d \bar{z} = \int_{\cc} ( \del_{\bar{z}} \tilf (z)) E(z) \d z \d \bar{z} + \Osd ( N^{\mfa-1} \| f''  \|_{1,w} )
\eeq
and
\begin{align} \label{eqn:exp-determ-2}
& \frac{1}{\pi} \int_{\cc} ( \del_{\bar{z}} \tilf (z) ) E(z) \d z \d \bar{z} = \frac{ f(2) + f (-2)}{4}  + \left( \sum_{j \neq a } \frac{ s_{aj}^{(4)}}{N^2} \right) \frac{1}{ 2 \pi} \int_{-2}^2 f(x) \frac{ x^4 - 4 x^2 + 2}{\sqrt{ 4 -x^2}} \d x \notag\\
+ & \tr (S) \frac{1}{ 2 \pi} \int_{-2}^2 f(x) \frac{ 2 -x^2}{ \sqrt{4 -x^2}} \d x + \left( \sum_a \frac{ s_{aa}^{(3)}}{N^{3/2}} \right) \frac{1}{ 8 \pi} \int_{-2}^2 f(x) \frac{ x^3-x^2-2x+4}{\sqrt{4-x^2}} \d x \notag\\
+ &  \frac{1}{ \pi} \int_{-2}^2 \frac{ f(x)}{ \sqrt{4-x^2}} \left( \Re[ \msc(x+ \i 0)^2 \tr( S (1- \msc(x + \i 0 )^2 S )^{-1} )] \right) \d x =: E_S (f) .
\end{align}
Furthermore,
\beq \label{eqn:exp-determ-4}
\left| \Re[ \msc(x+ \i 0)^2 \tr( S (1- \msc(x + \i 0 )^2 S )^{-1} )] \right| \lesssim 1
\eeq
for $ x \in (-2, 2)$. 
\eel
\proof By the Sherman-Morrison formula we have that,
\beq \label{eqn:sherm-1}
\tr (S(1 - \msc^2(z) S)^{-1} ) = \tr (S (1- \msc^2(z) A)^{-1} ) + \frac{ \msc(z)^2}{ 1 - \msc(z)^2},
\eeq
with $A$ defined as in Lemma \ref{lem:spectral-gap}. 
By Lemma \ref{lem:inverse} the first term on the RHS is bounded. Since $| \msc(z) | \leq 1 - c | \Im[z]|$ for some $c >0$ we have that,
\beq
|E(z) | \lesssim \frac{| \msc'(z) | }{ |\Im[z] |}.
\eeq
From this, the estimate $|\msc'(z)| \lesssim |4-\Re[z]^2|^{-1/2}$, and the explicit form of $\del_{\bar{z}} \tilf (z)$ (see \eqref{eqn:quasi-derivative}), the estimate \eqref{eqn:exp-determ-1} follows. 

We now turn to the computation of the LHS of \eqref{eqn:exp-determ-2}. 
By \cite[(4.69)]{meso} we have,
\begin{align}
\frac{1}{\pi} \int_{\cc} ( \del_{\bar{z}} \tilf (z) )\sum_{j \neq a} \frac{ s_{aj}^{(4)}}{N^2} \msc(z)^3 \msc'(z) \d z \d \bar{z} = \left( \sum_{j \neq a } \frac{ s_{aj}^{(4)}}{N^2} \right) \frac{1}{ 2 \pi} \int_{-2}^2 f(x) \frac{ x^4 - 4 x^2 + 2}{\sqrt{ 4 -x^2}} \d x
\end{align}
and by \cite[(4.68)]{meso}
\begin{align}
\frac{1}{ \pi} \int_{\cc} ( \del_{\bar{z}} \tilf (z) ) (- \msc(z) \msc'(z) \tr S ) \d z \d \bar{z} = \tr (S) \frac{1}{ 2 \pi} \int_{-2}^2 f(x) \frac{ 2 -x^2}{ \sqrt{4 -x^2}} \d x
\end{align}
and then by similar arguments as to the proof the above equalities,
\begin{align}
 & \left( - \sum_a \frac{ s_{aa}^{(3)}}{2 N^{3/2}} \right) \frac{1}{ \pi} \int_{\cc} ( \del_{\bar{z}} \tilf (z) ) \msc(z)^2 \msc'(z) \d z \d \bar{z} \notag\\
 =& \left( \sum_a \frac{ s_{aa}^{(3)}}{N^{3/2}} \right) \frac{1}{ 8 \pi} \int_{-2}^2 f(x) \frac{ x^3-x^2-2x+4}{\sqrt{4-x^2}} \d x .
\end{align}
We have only left to compute the contribution of the $\tr (S(1- \msc(z)^2 S)^{-1} )$ term of $E(z)$. We write,
\begin{align}
\msc(z) \msc'(z) \tr (S (1 - \msc(z)^2 S)^{-1} ) = \msc(z) \msc'(z) \left( \tr (A(1 - \msc(z)^2 A )^{-1} ) + \frac{1}{ 1 - \msc(z)^2 } \right).
\end{align}
By \cite[(4.70)]{meso},
\begin{align}
\frac{1}{ \pi} \int_{\cc} ( \del_{\bar{z}} \tilf (z ) ) \frac{ \msc(z) \msc'(z) }{ 1 - \msc(z)^2} \d z \d \bar{z} = - \frac{1}{ 2 \pi} \int_{-2}^2 \frac{ f(x)}{ \sqrt{4-x^2}} \d x + \frac{ f(2) - f(-2)}{4}.
\end{align}
By Green's theorem \eqref{eqn:green-2},
\begin{align}
& \frac{1}{\pi} \int_{\cc} ( \del_{\bar{z}} \tilf (z) ) \msc(z) \msc'(z) \tr (A (1- \msc(z)^2 A)^{-1} ) \d z \d \bar{z} \notag\\
= & \frac{1}{ \pi} \int_{-2}^2 \frac{ f(x)}{ \sqrt{4-x^2}} \Re\left[ \msc(x+ \i 0 )^2 \tr (A (1 - \msc(x + \i 0 )^2 A )^{-1} ) \right] \d x .
\end{align}
Therefore,
\begin{align}
 & \frac{1}{ \pi} \int_{\cc} ( \del_{\bar{z}} \tilf (z) ) \msc(z) \msc'(z) \tr (S (1- \msc (z)^2 S)^{-1} ) \d z \d \bar{z} \notag\\
= &\frac{f(2) - f(-2)}{4} + \frac{1}{ \pi} \int_{-2}^2 \frac{ f(x) }{ \sqrt{4-x^2}} \left( \Re[ \msc(x+ \i 0)^2 \tr( A (1- \msc(x + \i 0 )^2 A )^{-1} )] - \frac{1}{2} \right) \d x \notag\\
= & \frac{f(2) - f(-2)}{4}+ \frac{1}{ \pi} \int_{-2}^2 \frac{ f(x)}{ \sqrt{4-x^2}} \left( \Re[ \msc(x+ \i 0)^2 \tr( S (1- \msc(x + \i 0 )^2 S )^{-1} )] \right) \d x.
\end{align}
The second equality follows from,
\beq \label{eqn:sherm-3}
\msc(z)^2 \tr (S(1 - \msc(z)^2 S)^{-1} ) = \msc(z)^2 \tr ( A (1- \msc(z)^2 A)^{-1} ) + \frac{ \msc(z)^2}{ 1 -\msc(z)^2}
\eeq
and $\msc'(z) = \frac{ \msc(z)^2}{ 1 - \msc(z)^2}$ and $\Re[ \msc' (x + \i 0 ) ] = - \frac{1}{2}$. This completes the proof of \eqref{eqn:exp-determ-2}. By taking the real part on both sides of \eqref{eqn:sherm-3} we conclude the estimate \eqref{eqn:exp-determ-4}.

\qed

\bep \label{prop:exp-correction}
We have that,
\begin{align} \label{eqn:exp-correction}
\ee[ \tr f(H) ] = N \int f(x) \rhosc (x) \d x + E_S (f) + \Osd( N^{\mfa-1} \| f''\|_{1, w} )
\end{align}
\eep
\proof From \eqref{eqn:exp-determ-1}, \eqref{eqn:exp-determ-2} and \eqref{eqn:HS-est} it suffices to show that
\beq \label{eqn:exp-1a}
\int_{ \Oma} ( \del_{\bar{z}} \tilf (z) ) ( N\ee[ m_N (z) - \msc(z) ] - E(z) ) \d z \d \bar{z} = \Osd ( N^{-1} \| f''\|_{1, w} ).
\eeq
From \eqref{eqn:exp-second-2} we have that,
\beq
|N \ee[ m_N (z) - \msc(z) ] - E(z) | \sd \frac{1}{N \eta^2 \sqrt{ \eta + |4-E^2|}}.
\eeq
It is then straightforward to modify the proof of \eqref{eqn:H-est} to derive \eqref{eqn:exp-1a}. \qed

\section{Characteristic function} \label{sec:char}

In this section we complete our expansion for the characteristic function,
\beq
\ee[ \ea ] = \ee\left[ \exp \left( \i \lambda \int_{\Oma} ( \del_{\bar{u}} \tilf (u) ) ( \sum_a G_{aa} (u) - \ee[ G_{aa} (u) ] ) \d u \d \bar{u} \right) \right].
\eeq
Throughout this section we will assume that $|\lambda \Phi| \leq N^{-c}$ for some $c>0$. 
We will apply Stein's method which involves showing that,
\beq
\frac{ \d}{ \d \lambda} \ee[ \ea (\lambda) ] \approx - \lambda V(f) \ee[ \ea (\lambda)]
\eeq
for a certain function $V(f)$ to be computed later. Differentiating $\ee[ \ea(\lambda)]$ we see that we need to investigate the quantity, 
\beq
\sum_a \ee[ \ea (G_{aa}(z) - \ee[ G_{aa} (z) ] ) ].
\eeq
In order to compute this quantity, we will apply the expansion \eqref{eqn:intermediate-characteristic}, on which the above appears on the LHS. In order to use \eqref{eqn:intermediate-characteristic}, we must compute the first term on its RHS. We write this term as, 
\begin{align} \label{eqn:char-second-basic}
 & \sum_{ja} \frac{ s_{aj}}{N} \left( \ee[ \del_{ja} ( \ea G_{ja} ) ] - \ee[ \ea] \ee[ \del_{ja} G_{ja} ] \right) \notag \\
 = & - \frac{ \i \lambda}{ \pi}  \int_{\Oma} ( \del_{\bar{w}} \tilf (w) ) \del_w \sum_{ja} \frac{ s_{aj}}{N} (2 - \delta_{ja} ) \ee[ \ea G_{aj} (z) G_{ja} (w) ] \d w \d \bar{w} \notag\\
 - & \sum_{ja} \frac{ s_{aj}}{N} \ee[ \ea (G_{aa}(z) G_{jj} (z) - \ee[ G_{aa} (z) G_{jj} (z) ] ) ] \notag\\
 - & \sum_{j \neq a } \frac{ s_{aj}}{N} \ee[ \ea (G_{aj}(z)^2 - \ee[ G_{aj}(z)^2] ) ] 
\end{align}
and begin computing the various terms on the RHS of the above. 
The terms on the last two lines are computed in the following. We will assume that $ z= E \pm \i \eta$ with $ \eta >0$ and that $z \in \Oma$. 
\bel \label{lem:char-second-1}
We have,
\begin{align} \label{eqn:char-second-1}
 & \sum_{ja} \frac{ s_{aj}}{N} \ee[ \ea (G_{aa}(z) G_{jj} (z) - \ee[ G_{aa} (z) G_{jj} (z) ] ) ] = 2 \msc(z) \sum_a \ee[ \ea (G_{aa} (z) - \ee[ G_{aa} (z) ] ) ] \notag\\
& + \Osd ( N^{-1} \eta^{-2} + |\lambda| \Phi \Psi(z) \eta^{-1} )
\end{align}
and
\begin{align} \label{eqn:char-second-2}
& \left| \sum_{j \neq a } \frac{ s_{aj}}{N} \ee[ \ea (G_{aj}(z)^2 - \ee[ G_{aj}(z)^2] ) ] \right|  \sd   N^{-1} \eta^{-2}  \notag\\
+ &  ( \eta )^{-1} |\lambda| [(  \Phi^2 (1 + |\lambda| )  \Psi (z)  + \min \{ N^{-1}  \eta^{-1} , \Phi^2 \}]
+ N |\lambda| \Psi_1 (z)^2 \Phi \min \{ \Psi_1 (z),  \Phi \} 
\end{align}
\eel
\proof The estimate \eqref{eqn:char-second-2} follows immediately from Corollary \ref{cor:gij-main-est}. We can write the LHS of \eqref{eqn:char-second-1} as,
\begin{align} \label{eqn:char-second-1a}
& \sum_{ja} \frac{ s_{aj}}{N} \ee[ \ea (G_{aa}(z) G_{jj} (z) - \ee[ G_{aa} (z) G_{jj} (z) ] ) ] \notag\\
= & 2 \sum_{ja} \frac{s_{aj}}{N} \ee[ G_{jj}(z)] \ee[ \ea (G_{aa} (z) - \ee[ G_{aa}(z)\ ) ] \notag\\
&+ \sum_{ja} \frac{ s_{aj}}{N} \ee[ \ea (G_{aa}(z) - \ee[ G_{aa} (z) ] )(G_{jj}(z) - \ee[G_{jj} (z) ] ) ] \notag\\
&- \sum_{ja} \frac{ s_{aj}}{N} \ee[\ea] \ee[ (G_{aa}(z) - \ee[ G_{aa} (z) ) (G_{jj} (z) - \ee[ G_{jj}(z) ] ) ]
\end{align}
Due to Proposition \ref{prop:g11-est} and \eqref{eqn:fluct-av} we have,
\begin{align}
 & 2 \sum_{ja} \frac{s_{aj}}{N} \ee[ G_{jj}] \ee[ \ea (G_{aa} (z) - \ee[ G_{aa}(z)\ ) ] \notag\\
  =& 2 \msc(z) \sum_a \ee[ \ea (G_{aa} (z) - \ee[ G_{aa} (z) ] ) ] + \Osd (N^{-1} \eta^{-2} )
\end{align}
We now turn to the term on the second last line of \eqref{eqn:char-second-1a}. We have,
\begin{align} \label{eqn:char-second-1b}
& \sum_{ja} \frac{ s_{aj}}{N} \ee[ \ea (G_{aa}(z) - \ee[ G_{aa} (z) ] )(G_{jj}(z) - \ee[G_{jj} (z) ] ) ]  \notag\\
= & \sum_{ja} \frac{s_{aj}}{N} \ee[ \ea (G_{aa}(z) - \msc(z) ) (G_{jj} (z) - \msc(z) ) ] \notag\\
- & 2 \sum_{ja} \frac{ s_{aj}}{N} \ee[ \ea (G_{aa}(z) - \msc(z) )] \ee[ G_{jj}(z) - \msc(z) ]  \notag\\
+ & \sum_{ja} \ee[\ea] \frac{ s_{aj}}{N} \ee[G_{aa}- \msc(z) ] \ee[ G_{jj} - \msc(z) ]
\end{align}
By Proposition \ref{prop:giigaa-est}, the term on the second line is $\Osd ( N^{-1} \eta^{-2}  + |\lambda| \Phi \Psi(z) \eta^{-1})$. By Proposition \ref{prop:g11-est} and \eqref{eqn:fluct-av}, the terms on the third and fourth lines are $\Osd ( N^{-1} \eta^{-2} )$. Therefore,
\beq
\sum_{ja} \frac{ s_{aj}}{N} \ee[ \ea (G_{aa}(z) - \ee[ G_{aa} (z) ] )(G_{jj}(z) - \ee[G_{jj} (z) ] ) ]  = \Osd ( N^{-1} \eta^{-2} + |\lambda| \Phi \Psi(z) \eta^{-1} ).
\eeq
By almost the same argument, the same estimate (with $\lambda=0$) holds for the term on the fourth line of \eqref{eqn:char-second-1a}.  This completes the proof of \eqref{eqn:char-second-1}. \qed

Recall the notation $M := \msc(z) \msc(w)$. The following lemma computes the first term on the RHS of \eqref{eqn:char-second-basic}.
\bel \label{lem:char-second-2}
We have,
\begin{align} \label{eqn:char-var-1}
 & \frac{ \i \lambda}{\pi} \int_{\Oma} ( \del_{\bar{w}} \tilf (w) ) \del_w \sum_{j \neq a} \frac{ s_{aj}}{N} \ee[ \ea G_{aj} (z) G_{ja} (w) ] \d w \d \bar{w} \notag\\
= & \ee[ \ea] \frac{ \i \lambda}{\pi} \int_{\Oma} ( \del_{\bar{w}} \tilf (w) ) \del_w \left( - M \tr S + M \tr (S (1- M S)^{-1} ) \right) \d w \d \bar{w} +  |\lambda| \Osd ( \eta^{-1} \Phi^2(1 + |\lambda| ) )
\end{align}
and
\begin{align} \label{eqn:char-var-2}
 & \frac{ \i \lambda}{\pi} \int_{\Oma} ( \del_{\bar{w}} \tilf (w) ) \del_w \sum_{a} \frac{ s_{aa}}{N} \ee[ \ea G_{aa} (z) G_{aa} (w) ] \d w \d \bar{w} \notag\\
= & \ee[ \ea] \frac{ \i \lambda }{\pi} \int_{\Oma}  ( \del_{\bar{w}} \tilf (w) ) \del_w M \tr S \d w \d \bar{w} 
+  \Osd ( \eta^{-1} \Phi^2 |\lambda| ).
\end{align}
\eel
\proof The estimate \eqref{eqn:char-var-1} is an immediate application of Corollary \ref{cor:gij-main-est}, and using \eqref{eqn:H-est} to integrate the error terms in $w$. The estimate \eqref{eqn:char-var-2} follows from,
\begin{align}
\sum_a \frac{ s_{aa}}{N} G_{aa}(z) G_{aa}(w) = \msc(z) \msc(w) \tr (S) + \Osd ( N^{-1}y^{-1} + N^{-1} \eta^{-1} )
\end{align}
which follows from \eqref{eqn:entry-wise} and \eqref{eqn:fluct-av}. \qed

In the following we summarize our expansion for the characteristic function.  In preparation, define
\beq
\tilde{s}_4 := \sum_{a \neq j} \frac{ s_{aj}^{(4)}}{N^2} , \qquad \tilde{s}_3 := \sum_a \frac{ s_{aa}^{(3)}}{N} 
\eeq
and
\beq \label{eqn:ak-def}
a_{k, \mfa} := \frac{1}{\pi} \int_{\Oma} ( \del_{\bar{u} } \tilf (u ) ) \del_u (\msc(u)^k ) \d u \d \bar{u}.
\eeq
\bep \label{prop:char-stein-1}
We have,
\begin{align}
 & (z + 2 \msc (z) ) \sum_a \ee[ \ea (G_{aa}(z) - \ee[G_{aa} (z) ] ) ] \notag\\
 = & \ee[ \ea ] \frac{ \i \lambda}{\pi} \int_{\Oma} ( \del_{\bar{w}} \tilf (w) ) \del_w ( - 2 M \tr (S (1- MS)^{-1} +  M \tr S ) \d w \d \bar{w} \notag\\
 + & \ee[\ea] \bigg\{ - \tilde{s}_4 \i \lambda a_{2, \mfa} \msc(z)^2 + \frac{\tilde{s}_3 \i \lambda}{2 N^{1/2}}  \msc(z)  \left[ 2 \msc(z)  a_{1, \mfa} +  \i \lambda a_{1, \mfa}^2 +   a_{2, \mfa} \right] \bigg\} \notag\\
 + & \Osd \left( N^{-1} ( 1 + |\lambda| )^3  + N^{-1} \eta^{-2} + |\lambda \Phi \Psi(z) | \eta^{-1} + |\lambda| \eta^{-1} \Phi^2 (1 + |\lambda| ) \right) \label{eqn:good-Gaa-equation}
\end{align}

Therefore,
\begin{align} \label{eqn:stein-1}
 & \frac{ \d }{ \d \lambda} \ee[ \ea ( \lambda) ] = - \lambda V_\mfa (f) \ee[ \ea ( \lambda) ] + \i \lambda^2 B_\mfa (f) \ee[ \ea ( \lambda) ] \notag\\
+ & \Osd ( (1+ |\lambda| )^3 N^{-1} + \Phi_w^2 (1 + |\lambda| )^2  )
\end{align}
where,
\begin{align} \label{eqn:Va-def}
& V_\mfa (f) := \frac{1}{\pi^2} \int_{\Oma^2} ( \del_{\bar{z} } \tilf (z) ) ( \del_{\bar{w}} \tilf (w) ) \del_w \left( \msc'(z) \msc(w)(2 \tr (S (1- M S)^{-1}) - \tr S)\right) \d w \d \bar{w} \d z \d \bar{z} \notag\\
& + \frac{1}{2} \tilde{s}_4 a_{2, \mfa}^2 - \frac{ \tilde{s}_3}{ N^{1/2}}  a_{2, \mfa} a_{1, \mfa} 
\end{align}
and
\beq
B_\mfa (f) := \frac{ \tilde{s}_3}{2 N^{1/2}} ( a_{1, \mfa} )^3.
\eeq
\eep
\proof The equation \eqref{eqn:good-Gaa-equation} follows from \eqref{eqn:intermediate-characteristic}, \eqref{eqn:char-second-basic} (which expands the first term on the RHS of \eqref{eqn:intermediate-characteristic}) and Lemmas \ref{lem:char-second-1} and \ref{lem:char-second-2} (which compute the terms on the RHS of \eqref{eqn:char-second-basic}). The equation \eqref{eqn:stein-1} then follows from using \eqref{eqn:good-Gaa-equation} to solve for the expression $\sum_a \ee[ \ea ( G_{aa} (z) - \ee[ G_{aa} (z) ] ) ]$ in the equality,
\begin{align}
\frac{ \d }{ \d \lambda} \ee[ \ea ( \lambda) ] = \frac{ \i }{ \pi} \int_{\Oma} ( \del_{\bar{z}} \tilf (z) ) \sum_a \ee[ \ea (G_{aa} (z) - \ee[ G_{aa} (z) ] ) ] \d z \d \bar{z},
\end{align}
and using a slight modification of Lemma \ref{lem:H-est} to integrate the error term (i.e., we use the estimate $|z + 2 \msc(z) |^{-1} \leq C / \sqrt{ |E-2| + \eta}$ and arrive at errors with $\| f''\|_{1, w}$ instead of $\| f''\|_1$). \qed

\subsection{Computation of deterministic coefficients}

\subsubsection{Chebyshev polynomials} \label{sec:cheby-1}

We recall now the Chebyshev polynomials of the first kind $T_n(x)$ and coefficients $t_n(f)$ as defined in Definition \ref{def:cheby}. We note the 
the orthogonality relation,
\begin{align}
\int_{-2}^2 T_n(x) T_n (x) \frac{ \d x }{ \sqrt{ 4-x^2}} = \begin{cases} 0, & \mbox{if }n \neq m \\ \pi, & \mbox{if }n=m=0 \\ \frac{\pi}{2} , & \mbox{if }n = m \neq 0 \end{cases} ,
\end{align}
so that we have the expansion, 
\beq
f = T_0 (x) \frac{ t_0(f)}{2} + \sum_{n \geq 1 }T_n (x) t_n (f).
\eeq
as $T_n(x)$ is a complete orthogonal family of  $L^2 ((-2, 2), \frac{ \d x}{ \sqrt{4-x^2}} )$. 
We have for $ x \in (-2, 2)$,
\beq \label{eqn:Chebyshev-explicit}
T_n (x) = \frac{(-1)^n}{2} \left( \msc(x + \i 0 )^n + \msc(x - \i 0)^n \right) = (-1)^n \Re\left[ \msc(x + \i 0)^n \right]
\eeq
which is easily derived from \eqref{eqn:cheby-def}. In particular, 
\beq\label{eqn:tn-msc}
t_n (f) = \frac{2 (-1)^n}{ \pi} \int_{-2}^2 f(x) \Re[ \msc(x + \i 0)^n] \frac{ \d x }{ \sqrt{4-x^2}}.
\eeq

\subsubsection{Easy coefficients} \label{sec:easy}

In this section we compute the coefficients $a_{k, \mfa}$ defined in \eqref{eqn:ak-def} in terms of the $t_n(f)$.

\bel \label{lem:easy}
For each $k \geq 1$,
\begin{align} \label{eqn:easy-coeff-1} 
a_{k, \mfa} = \frac{1}{\pi} \int_{\Oma} ( \del_{\bar{u}} \tilf (u ) ) \del_u ( \msc^k (u) ) \d u \d \bar{u} = \frac{1}{\pi} \int_{\cc} ( \del_{\bar{u}} \tilf (u ) ) \del_u ( \msc^k (u) ) \d u \d \bar{u} + \Osd ( (N^{\mfa-1})^{3/2} \| f''\|_1 ) 
\end{align}
and for $k \geq 1$,
\beq \label{eqn:easy-coeff-2}
\frac{1}{\pi} \int_{\cc} ( \del_{\bar{u}} \tilf (u ) ) \del_u ( \msc^k (u) ) \d u \d \bar{u} = (-1)^k \frac{k}{2} t_k (f).
\eeq
\eel
\proof The first equality in \eqref{eqn:easy-coeff-1} is by definition. The  estimate in \eqref{eqn:easy-coeff-1} follows from the fact that $|\del_z \msc(z) | \lesssim |\Im[z]|^{-1/2}$ and direct computation. 

For \eqref{eqn:easy-coeff-2} we first note that Green's theorem \eqref{eqn:green-2} implies that,
\begin{align}
\int_{\cc} ( \del_{\bar{u}} \tilf (u)) \frac{ \msc(u)^k}{ 2\msc(u) + u } \d u \d \bar{u} =  \int_{\rr} f(x) \Im\left[ \frac{ \msc( x + \i 0)^k}{ x + 2 \msc ( x + \i 0) } \right] \d x .
\end{align}
The equality \eqref{eqn:easy-coeff-2} follows from this, as well as the identities, $\del_u ( \msc(u)^k) = - k \msc(u)^k / (u + 2 \msc (u ) )$ and 
\beq
\Im\left[ \frac{ \msc( x + \i 0)^k}{ x + 2 \msc ( x + \i 0) } \right] = -  \frac{ \Re[ \msc(x + \i 0)^k ] }{ \sqrt{4-x^2}},
\eeq
and \eqref{eqn:tn-msc}.  
\qed

\subsubsection{Variance functional} \label{sec:var-functional}

In this section we compute the variance function $V_\mfa (f)$ as defined in \eqref{eqn:Va-def}. 
Let us define,
\beq \label{eqn:var-F-def}
F(z, w) := \del_w \left( \msc'(z) \msc(w) ( 2 \tr ( S(1- \msc(z) \msc(w) S )^{-1} ) - \tr S ) \right)
\eeq

\bel \label{lem:var-1}
We have,
\begin{align} \label{eqn:var-aa1}
 & \int_{\Oma^2} ( \del_{\bar{z}} \tilf ) ( \del_{\bar{w}} \tilf (w) ) F(z, w)  \d z \d \bar{z} \d w \d \bar{w} \notag\\
= & \int_{\cc^2} ( \del_{\bar{z}} \tilf ) ( \del_{\bar{w}} \tilf (w) )  F(z, w)  \d z \d \bar{z} \d w \d \bar{w} + \Osd ( (N^{\mfa-1}( \| f''  \|_{1,w} + \|f'\|_{1, w} + \|f \|_{1, w}))^2 )
\end{align}
\eel
\proof By the Sherman Morrison formula,
\beq \label{eqn:sherman-1}
\tr (S(1 - MS )^{-1}) = \tr (S(1 - MA)^{-1}) + \frac{ \msc(z) \msc(w)}{ 1 - \msc(z) \msc(w) }.
\eeq
Since $\| M A \|_{\ell^2 \to \ell^2 } \leq 1 -c $ for some $c>0$ we can conclude, using Lemma \ref{lem:inverse} that
\beq
\left| \del_w [ \msc' (z) \msc(w)  \tr ( S(1 - MA )^{-1} )] \right| \lesssim  | \msc'(z) \msc'(w) |.
\eeq
Since $|1 - \msc(z)\msc(w) | \geq c (y + \eta)$ we conclude that,
\beq \label{eqn:var-aa2}
|F(z, w) | \lesssim \frac{ |\msc'(z) \msc'(w)|}{ (y + \eta)^2}.
\eeq
From this estimate we easily conclude, by direct integration,
\begin{align}
\left| \int_{\Oma^c \times \Oma^c}   ( \del_{\bar{z}} \tilf (z))   ( \del_{\bar{w}} \tilf (w) (w)) F(z, w) \d z \d \bar{z} \d w \d \bar{w} \right| \sd (N^{\mfa-1} \| f'' \|_{1, w} )^2.
\end{align}
as well as,
\beq
\left| \int_{\Oma^c \times \Oma}   ( \del_{\bar{z}} \tilf ) ( \del_{\bar{w}} \tilf (w) ) F(z, w) \d z \d \bar{z} \d w \d \bar{w} \right| \sd ( N^{\mfa-1} ( \| f'' \|_{1, w} + \| f' \|_{1, w} + \| f \|_{1, w} ) )^2 ,
\eeq
which yields the claim.
 \qed

%
%

\bep \label{prop:var-1}
We have that,
\begin{align} \label{eqn:determ-var-1}
 \frac{1}{ \pi^2} & \int_{\cc^2} ( \del_{\bar{z}} \tilf (z) ) ( \del_{\bar{w}} \tilf (w) ) \del_w ( \msc'(z) \msc(w) \tr (S (1 - \msc(z) \msc(w) S)^{-1} ) \d z \d \bar{z} \d w \d \bar{w}  \notag\\
= & \frac{1}{4} \sum_{j=1}^\infty j t_j (f)^2 \tr (S^j)
\end{align}
as well as,
\beq \label{eqn:determ-var-2}
\frac{1}{\pi^2} \int_{\cc^2} ( \del_{\bar{z}} \tilf (z) ) ( \del_{\bar{w}} \tilf (w) ) \msc'(z) \msc'(w) \tr S = \tr S \frac{1}{4} t_1(f)^2
\eeq

\eep
\proof For $0 < \delta < 1$ we have $|1 - \delta \msc(z) \msc(w) | \geq 1 - \delta$ and so,
\begin{align} 
 & \left| \frac{1}{ 1 - \delta \msc(z) \msc(w) } - \frac{1}{ 1- \msc(z) \msc(w)} \right| \leq \frac{ 1-\delta}{ |1- \delta \msc(z) \msc(w) | } \frac{1}{ |1 - \msc(z) \msc(w) | } \notag\\
\leq & \frac{1}{ |1 - \msc(z) \msc(w) | } \lesssim \frac{1}{ \eta + y}
\end{align}
Therefore using \eqref{eqn:sherman-1} and dominated convergence we find,
\begin{align}
\int_{\cc^2} ( \del_{\bar{z}} \tilf (z) ) ( \del_{\bar{w}} \tilf (w) ) \del_w ( \msc'(z) \msc(w) \tr (S (1 - \msc(z) \msc(w) S)^{-1} ) \d z \d \bar{z} \d w \d \bar{w}  \notag\\ 
= \lim_{\delta \to 1^-} \int_{\cc^2} ( \del_{\bar{z}} \tilf (z) ) ( \del_{\bar{w}} \tilf (w) ) \del_w ( \msc'(z) \msc(w) \tr (S (1 -  \delta \msc(z) \msc(w) S)^{-1} ) \d z \d \bar{z} \d w \d \bar{w}  .
\end{align}
Now,
\begin{align}
& \del_w  ( \msc'(z) \msc(w) \tr (S (1 -  \delta \msc(z) \msc(w) S)^{-1} )  = \del_w \sum_{j=0}^\infty \msc'(z) \msc(w)^{j+1} \msc(z)^j \delta^j \tr S^{j+1} \notag\\
= & \sum_{j=0}^\infty (j+1) \msc'(z) \msc'(w)( \delta \msc(z) \msc(w) )^j \tr S^{j+1} .
\end{align}
The integration of summation and differentiation is justified because $|\msc(z) | \leq 1$ and
\beq \label{eqn:trS-bd}
\left| \tr S^j \right| \leq N \| S \|_{\ell^2 \to \ell^2}^j \leq N.
\eeq
Therefore,
\begin{align}
 & \frac{1}{ \pi^2} \int_{\cc^2} ( \del_{\bar{z}} \tilf (z) ) ( \del_{\bar{w}} \tilf (w) ) \del_w ( \msc'(z) \msc(w) \tr (S (1 -  \delta \msc(z) \msc(w) S)^{-1} ) \d z \d \bar{z} \d w \d \bar{w} \notag\\
= &\sum_{j=1}^\infty \delta^{j-1} \frac{ j t_j (f)^2}{4} \tr (S^j) ,
\end{align}
using \eqref{eqn:tn-msc} as in the proof of Lemma \ref{lem:easy}. Again, interchanging the order of integration and summation is justified using that $|\msc(z) | \leq 1$ and \eqref{eqn:trS-bd}. Finally, by monotone convergence,
\beq
\lim_{\delta \to 1^-}  \sum_{j=1}^\infty \delta^{j-1} \frac{ j t_j (f)^2}{4} \tr (S^j) = \frac{1}{4} \sum_{j=1}^\infty j t_j (f)^2 \tr (S^j)
\eeq
which completes the proof of \eqref{eqn:determ-var-1}. The equation \eqref{eqn:determ-var-2} follows from \eqref{eqn:easy-coeff-2}. 
\qed

\bel \label{lem:var-2}
We have that,
\beq \label{eqn:Va-V}
V_\mfa (f) = V_1 (f) + \Osd( N^{\mfa-1} \left( \| f''\|_{1,w} + \| f'\|_{1, w} + \| f \|_{1, w} \right)
\eeq
where $V_1 (f)$ is defined in \eqref{eqn:intro-V-def}.  
Furthermore, 
\beq \label{eqn:V-pos}
V_1(f) \geq 0.
\eeq
We have also,
\beq \label{eqn:Ba-s}
B_\mfa (f) = \frac{ \tilde{s}_3}{8 N^{1/2}} t_1(f)^3 + \Osd ( N^{\mfa-3/2} \| f''\|_1 )
\eeq
\eel
\proof First of all since $\hat{s}_{4,1} = \tilde{s}_4 + \O ( N^{-1} )$ and $|t_2 (f) | \lesssim \| f \|_{1, w}$, we have that \eqref{eqn:Va-V} follows from Lemmas \ref{lem:easy} and \ref{lem:var-1} and Proposition \ref{prop:var-1}. Similarly, \eqref{eqn:Ba-s} follows from Lemma \ref{lem:easy}. It  remains to prove \eqref{eqn:V-pos}. 

%
%
%
%
%
%
%
%
%

Let $a = t_1 (f)$ and $b= t_2 (f)$ for notational simplicity. We have,
\begin{align} \label{eqn:var-nonneg-1}
V_1(f) \geq&  \frac{a^2}{4} \tr (S) + b^2 \left(  \tr(S^2) + \frac{\hat{s}_{4, 1}}{2}\right) + \frac{\tilde{s}_3}{2 N^{1/2}} ab \notag\\
= & b^2 \left(\frac{1}{N^2}  \sum_{i \neq j } s_{ij}^2 + \frac{ s_{ij}^{(4)}}{2} \right) \notag\\
+ & \frac{a^2}{4} \tr (S) + b^2 \left( \frac{1}{N^2} \sum_i s_{ii}^2 + \frac{ s_{ii}^{(4)}}{2} \right) + ab \frac{ \tilde{s}_3}{2N^{1/2}}
\end{align}
We recall that the third and fourth cumulants of a centered random variable $X$ are given by,
\beq
\kappa_3 (X) = \ee[ X^3], \qquad \kappa_4 (X) = \ee[ X^4] - 3 \ee[ X^2]^2.
\eeq
The nonnegativity of the RHS of \eqref{eqn:var-nonneg-1} follows from the following two inequalities. First,
\beq \label{eqn:var-nonneg-3}
\frac{1}{N^2} \sum_{ i \neq j } 2 s_{ij}^2 +  s_{ij}^{(4)}= \sum_{i \neq j } \ee[ H_{ij}^4] - \ee[ H_{ij}^2]^2 = \sum_{i \neq j } \ee[ (H_{ij}^2 - \ee[ H_{ij}^2] )^2] \geq 0
\eeq
Second, 
\begin{align} \label{eqn:var-nonneg-2}
& \left| \frac{ \tilde{s}_3}{2 N^{1/2}} ab \right|= \left| \frac{1}{2} \Cov \left( a \sum_i H_{ii} , b \sum_i H_{ii}^2 \right)\right| \leq \frac{1}{4} \Var\left( a \sum_i H_{ii} \right) + \frac{1}{4} \Var \left( b \sum_i H_{ii}^2 \right) \notag\\
= & \frac{a^2 \tr (S) }{4} + \frac{b^2}{4 N^2} \left( \sum_i  2 s_{ii}^2  +  s_{ii}^{(4)} \right) \notag\\
\leq & \frac{a^2 \tr (S) }{4} + \frac{b^2}{2 N^2} \left( \sum_i  2 s_{ii}^2  + s_{ii}^{(4)}  \right)  .
\end{align}
The last inequality follows because $\kappa_4 (X) + 2 \kappa_2 (X)^2 \geq 0$ for any random variable $X$. 
\qed

Define,
\beq
\hat{V}_p (f) := \frac{1}{2} \sum_{j=1}^\infty j t_j(f)^2 \tr S^j.
\eeq
\bel \label{lem:var-3}
We have that,
\begin{align}
\hat{V}_p (f) & = \frac{1}{2 \pi^2} \int_{-2}^2 \int_{-2}^2 \left( \frac{ f(x) - f(y)}{ x -y} \right)^2 \frac{ 4 - xy}{ \sqrt{4-x^2} \sqrt{4-y^2}} \d x \d y  \notag\\
&+ \frac{1}{ \pi^2} \int_{-2}^2 \int_{-2}^2 f(x) f(y) \frac{ g (x, y)}{ \sqrt{4-x^2} \sqrt{4-y^2} } \d x \d y
\end{align}
where,
\beq
g(x, y) := \Re\left[ \tr \frac{ \msc (x + \i 0) \msc(y+ \i 0) A}{(1 - \msc (x + \i 0) \msc(y+ \i 0) A)^2} + \tr \frac{ \msc (x + \i 0) \msc(y- \i 0) A}{(1 - \msc (x + \i 0) \msc(y- \i 0) A)^2}\right] 
\eeq
The function $g$ obeys $|g(x, y) | \lesssim 1$ for all $x, y \in (-2, 2)$.
\eel
\proof We need to re-compute the integral on the LHS of \eqref{eqn:determ-var-1}. By \eqref{eqn:sherman-1} we have that
\begin{align}
 & \del_w ( \msc'(z) \msc(w) \tr ( S(1- \msc(z) \msc(w) S)^{-1} ) ) = \frac{ \msc'(z) \msc'(w)}{ (1 - \msc(z) \msc(w) )^2} \notag\\
 + & \msc'(z) \msc'(w) \tr \frac{A}{ (1 - M A)^2}.
\end{align}
By \cite[(4.60)]{meso} and \eqref{eqn:msc-dif-identity} we have,
\begin{align}
& \frac{2}{ \pi^2} \int_{\cc^2} ( \del_{\bar{z}} \tilf (z) ) ( \del_{\bar{w}} \tilf (w) ) \frac{ \msc'(z) \msc'(w)}{ (1 - \msc(z) \msc(w) )^2} \d z \d \bar{z} \d w \d \bar{w}  \notag\\
= &
\frac{1}{2 \pi^2} \int_{-2}^2 \int_{-2}^2 \left( \frac{ f(x) - f(y)}{ x -y} \right)^2 \frac{ 4 - xy}{ \sqrt{4-x^2} \sqrt{4-y^2}} \d x \d y 
\end{align} 
On the other hand, by applying Green's theorem \eqref{eqn:green-2} we see that,
\begin{align}
 & \frac{2}{ \pi^2} \int_{\cc^2} ( \del_{\bar{z}} \tilf (z) ) ( \del_{\bar{w}} \tilf (w) ) \msc'(z) \msc'(w) \tr \frac{A}{ (1- MA)^2} \d z \d \bar{z} \d w \d \bar{w} \notag\\
= & \frac{1}{ 2 \pi^2} \int_{-2}^2 \int_{-2}^2 \frac{ f(x) f(y)}{ \sqrt{4-x^2} \sqrt{4-y^2}}  \bigg\{  h (x + \i 0, y + \i 0) + h ( x - \i 0, y - \i 0) \notag\\
 & \qquad \qquad \qquad + h(x - \i 0, y - \i 0) + h (x + \i 0, y - \i 0)  \bigg\}  \d x \d y
\end{align}
where $h(z, w) = \msc(z) \msc(w) \tr (A (1- MA)^{-2})$. The claim now follows, with the bound on $g$ following 
 Lemma \ref{lem:inverse}. \qed

\subsection{Stein's method}

We finally apply Stein's method to compute the characteristic function $\ee[ \ea ( \lambda)]$.

\bel \label{lem:stein}
Suppose there is a $c>0$ so that $| \lambda \Phi | \leq N^{-c}$. We have,
\begin{align} \label{eqn:stein-2}
&\frac{ \d }{ \d \lambda} \ee[ \ea ( \lambda) ] = - \lambda V_1 (f) \ee[ \ea ( \lambda) ] + \i \lambda^2 B (f) \ee[ \ea ( \lambda) ] \notag\\
+ & \Osd ( (1 + | \lambda| )^3 N^{-1} + \Phi_w^2 ( 1 + |\lambda| )^2 ) + \Osd( |\lambda| N^{-1+\mfa} ( \| f''\|_{1, w} + \| f'\|_{1, w} + \| f\|_{1, w} )
\end{align}
Consequently, if $V_1(f) \geq c_1$ for some $c_1>0$ then, for $|\lambda \Phi| \leq N^{-c}$ we have,
\begin{align} \label{eqn:result-1}
&\ee[ \ea  ( \lambda) ] = \e^{ - \lambda^2 V_1(f)/2 + \i \lambda^3 B(f)/3 }  + \Osd ( \Phi_w^2 (1 + |\lambda|)  + N^{-1} ( 1 + |\lambda| )^2 ) \notag\\
+ &\Osd(  N^{-1+\mfa} ( \| f''\|_{1, w} + \| f'\|_{1, w} + \| f\|_{1, w} ) )
\end{align}
Otherwise for $|\lambda \Phi| \leq N^{-c}$ we have,
\begin{align} \label{eqn:result-2}
&\ee[ \ea ( \lambda) ] = \e^{ - \lambda^2 V_1(f)/2 + \i \lambda^3 B(f)/3 } + \Osd ( (1 + |\lambda| )^4 N^{-1} + \Phi_w^2 (1 + |\lambda| )^3 ) \notag\\
 +& \Osd( |\lambda|^2 N^{-1+\mfa} ( \| f''\|_{1, w} + \| f'\|_{1, w} + \| f\|_{1, w} )
\end{align}
\eel
\proof  The estimate \eqref{eqn:stein-2} follows from \eqref{eqn:stein-1} and \eqref{eqn:Va-V} and \eqref{eqn:Ba-s}. Letting now $\psi (\lambda) := \ee[ \ea ( \lambda)]$ and $\Phi (\lambda) := \e^{ - \lambda^2V_1(f)/2 + \i \lambda^3 B(f)/3}$ we have by differentiating $\psi(\lambda)\Phi(\lambda)^{-1}$ that,
\begin{align}
| \psi ( \lambda) - \Phi ( \lambda) | \sd \int_0^{|\lambda|} \e^{ (s^2-\lambda^2) V_1 (f)/2} \eps (s) \d s
\end{align}
with
\beq
\eps(s) := (1 + |s| )^3 N^{-1} + \Phi_w^2 ( 1 + |s| )^2 + |s| N^{-1+ \mfa} ( \| f''\|_{1, w} + \| f'\|_{1, w} + \| f\|_{1, w} ).
\eeq
The estimate \eqref{eqn:result-2} follows from the fact that $| \e^{ (s^2-\lambda^2) V_1 (f)/2}  | \leq 1 $ for $|s| \leq |\lambda|$ since $V_1(f) \geq 0$ from \eqref{eqn:var-nonneg-1}. The estimate \eqref{eqn:result-1} follows from using that
\beq
\int_0^{|\lambda|} \e^{ s^2 V_1(f) /2 } |s|^k \d s \lesssim (1+ |\lambda|)^{k-1} \e^{ \lambda^2 V(f)/2} .
\eeq
for $k \geq 0$ if $V_1(f) \geq c_1$. This completes the proof. \qed

\subsection{Proof of main theorems in the real symmetric case} \label{sec:real-proofs}

In this section we prove Theorem \ref{thm:main} and Proposition \ref{prop:funct} in the real symmetric case.

\subsubsection{Proof of Theorem \ref{thm:main}} \label{sec:main-proof}
 The estimate for the expectation follows from Proposition \ref{prop:exp-correction}. 
The two estimates for the characteristic function follow from Lemma \ref{lem:stein}, \eqref{eqn:HS-est} and the just proved estimate for the expectation. \qed

\subsubsection{Proof of Proposition \ref{prop:funct}} \label{sec:funct-proof}

The lower bound for $V_1 (f)$ follows from \eqref{eqn:V-pos}. The upper bound for the part of $V_1 (f)$ that we denoted by $\hat{V}_p (f)$ can be deduced from bounding the integral on the LHS of \eqref{eqn:determ-var-1} using \eqref{eqn:var-aa2} and Lemma \ref{lem:H-est}. The remaining components can be bounded by \eqref{eqn:determ-var-2} and the fact that $|t_n(f) | \leq C_n \| f \|_{1, w}$ for any $n>0$. The estimate for $B(f)$ also follows from this. Finally, the estimate for $E_1 (f)$ follows from its definition as well as \eqref{eqn:exp-determ-4}. \qed

%
%
%
%
%

\section{Application: max of log-characteristic polynomial} \label{sec:max}

In this section we show how to extend the results of \cite{BLZ} to generalized Wigner matrices. We will discuss only the real symmetric case, the complex Hermitian case being similar. In this section we define the complex logarithm by $\log (r \e^{ \i \theta} ) := \log (r) + \i \theta $ for $\theta \in (- \pi, \pi]$. Introduce,
\beq
L_N (z) := \sum_{j=1}^N \log (z- \lambda_j ) - N \int_{\rr} \log (z - x) \rhosc (x) \d x
\eeq
where $\lambda_j$ are the eigenvalues of a generalized Wigner matrix $H$. We will denote $z = E + \i \eta$ and assume $\eta \geq 0$.  We will make the following assumption on the tail of the distribution of entries of $H$.
\bas \label{ass:tail}
We assume that there is a $c>0$ so that 
\beq
\pp\left[ |H_{ij} | > \sqrt{N} x \right] \leq c^{-1} \e^{ - x^c}
\eeq
for all $x >0$, uniformly in $i$ and $j$. 
\eas
We will give a proof of the following. It extends Theorem 1.2 and Theorem 1.8(i) from the Wigner to generalized Wigner setting.
\bet \label{thm:max}
Let $H$ be a generalized Wigner matrix satisfying Assumption \ref{ass:tail}. Then, for any $\eps$ and $\kappa >0$,
\beq \label{eqn:max-main-1}
\pp\left[ \sup_{ |E| < 2 - \kappa } \frac{ \Re[ L_N(E)]}{ \sqrt{2} \log N } \in (1- \eps, 1+ \eps ) \right] = 1 + o (1),
\eeq
as well as, (for any choice of $\pm$)
\beq \label{eqn:max-main-2}
\pp\left[ \sup_{ |E| < 2 - \kappa } \frac{ \pm \Im[ L_N(E)]}{ \sqrt{2} \log N } \in (1- \eps, 1+ \eps ) \right] = 1 + o (1),
\eeq
and
\beq \label{eqn:max-main-3}
\pp\left[ \sup_{ \kappa N \leq k \leq ( 1- \kappa ) N } \frac{ \pi}{ \sqrt{2}} \rhosc ( \gamma_k) \frac{ \pm N( \lambda_k - \gamma_k ) }{\log N} \in (1-\eps, 1 + \eps) \right] = 1 + o (1)
\eeq
for any choice of $\pm$.
\eet

The main input from our work will be the following. 
\bep \label{prop:max-lss}
Fix $K>0, \frac{1}{2} > \eps >0, \kappa >0$. We have that,
\beq
\pp\left[ \sup_{z : |E| \leq 2 - \kappa, \e^{- K ( \log \log N)^2} \leq \eta \leq 1 } |L_N (z) | > ( \log N)^\eps \right] \leq C \e^{ - ( \log N)^{\eps} }
\eeq
for some $C>0$.
\eep
\proof We show how to bound the real part of $L_N (z)$, with the imaginary part being similar.  We start by proving the bound for fixed $z$.  Define,
\beq
X_z := \frac{1}{2} \sum_i \log ( (\lambda_i-E)^2 + \eta^2) - \frac{N}{2} \int \log ((x-E)^2 + \eta^2) \rhosc (x) \d x = \Re[ L_N(z)]
\eeq
Let us denote by $V_z, \E_z$ and $B_z$ the coefficients in Theorem \ref{thm:main} with the function $\Re[ \log (z - \cdot ) ]$. By \eqref{eqn:max-var-2} below we have, for $z$ as in the statement of the proposition, that $V_z = - \log ( \eta) + \O (1)$, and that $|B_z| \leq C N^{-1/2}$ and $|\E_z| \leq C$. 
It therefore follows from Theorem \ref{thm:main} that
\beq \label{eqn:max-char}
\ee[ \e^{ \i \lambda X_z} ] = \e^{ - \lambda^2 V_z/2  + \i \lambda \E_z } + \Osd ( N^{-1/2} ( 1 + |\lambda| )^4) , 
\eeq
for all $|\lambda| \leq N^{1/2-1/10}$.
Let $u \geq 2$ and $0.01> \delta >0$. A straightforward argument using the rigidity estimates  \cite[Theorem 2.2]{erdHos2012rigidity} shows that,
\begin{align}
\pp\left[ X_z-\E_z > u \right] \leq \pp\left[ u  < X_z - \E_z < N^{\delta} \right] + N^{-D}.
\end{align}
for any $D>0$. 
Let $\chi$ be a smooth function s.t. $\chi (x) = 1$ for $x \in (u , N^{\delta})$ and $\chi (x) = 0$ for $x \notin (u-1, N^{\delta} + 1 )$. We have that for any $M>0$,
\beq \label{eqn:max-chi}
| \hat{\chi} ( \lambda) | \leq N^{\delta} C_M (1 + |\lambda| )^{-M}.
\eeq
 Then,
\begin{align}
 \pp\left[ u < X_z - \E_z < N^{\delta} \right] &\leq \ee[ \chi (X_z-\E_z) ] \notag \\
 = & \int_{\rr} \hat{ \chi} ( \lambda) \ee[ \e^{ \i \lambda X_z } ] \d \lambda  \notag \\
 = & \int_{ |\lambda| \leq N^{\delta}} \hat{ \chi} ( \lambda) \ee[ \e^{ \i \lambda (X_E -\E_z)} ] \d \lambda + \O ( N^{-1} )  \notag \\
 = & \int_{ | \lambda| \leq N^{\delta} } \hat{ \chi ( \lambda) } \e^{ - \lambda^2 V_z /2 } \d \lambda+ \O ( N^{6 \delta}N^{-1/2}  ). \notag\\
 = & \int_{ \rr } \hat{ \chi ( \lambda) } \e^{ - \lambda^2 V_z /2  } \d \lambda + \O ( N^{6\delta - 1/2} )
\end{align}
Above, the third and fifth lines follow from \eqref{eqn:max-chi} and the fact that $V_z \geq 0$. In the fourth line we used \eqref{eqn:max-char}. Now let $Z$ be a standard normal random variable. Then,
\begin{align}
& \int_{ \rr } \hat{ \chi ( \lambda) } \e^{ - \lambda^2 V_z /2  } \d \lambda = \ee[ \chi ( V_z^{1/2} Z )] \leq \pp\left[ V_z^{1/2} Z > (u-1) \right] \leq C \e^{- (u-1)^2 / (2 V_z )}
\end{align}
From all of this we conclude that for all $0 < \eps < 1/2$ that
\beq \label{eqn:max-est-1}
\pp\left[ |X_z| > ( \log N)^{\eps} ]\right] \leq C \e^{ - (\log N)^{\eps} }
\eeq
for some $C>0$.

We have that
\beq
\del_E X_z = N \Re[ m_N (z) - \msc (z) ], \qquad \del_\eta X_z = N \Im[ m_N (z) - \msc (z) ].
\eeq
Let $D_{r} (z)$ be the disc of radius $r$ centered at $z= E + \i \eta$. Assume that $r < \eta/2$. From the above estimates and \cite[Theorem 2.1]{erdHos2012rigidity} we have that for some $C_1>0$,
\beq \label{eqn:max-est-2}
\pp\left[  \sup_{w \in D_{r}} |X_z-X_w | >   \e^{ C_1 ( \log \log N)^2 } \frac{r}{\eta} \right] \leq N^{-D} ,
\eeq
We conclude the proof by taking a union bound of the estimate \eqref{eqn:max-est-1} over $\e^{ C_2 ( \log \log N)^2}$ points, where $C_2$ is chosen to depend on $C_1$ and $K$. The claim follows.  \qed

\subsection{Proof of Theorem \ref{thm:max}}

In this section we detail how the proof of \cite{BLZ} can be modified to obtain the result for generalized Wigner matrices. The main point is that Proposition \ref{prop:max-lss} is a sufficient substitute for \cite[Corollary A.2]{BLZ} in order to obtain the main results. We first note that the main inputs of \cite{BLZ} from other works are: the rigidity estimates, from \cite{erdHos2012rigidity}, listed as \cite[Theorem 2.2]{BLZ}; the relaxation of eigenvalues under DBM, listed as \cite[Proposition 4.1]{BLZ}, using \cite{bourgade2018extreme} as input; and the comparision techniques of \cite{landon2020comparison}, developed further in \cite[Section 5]{BLZ}. All of these three input works prove their results for generalized Wigner matrices, and so these aspects of \cite{BLZ} hold for the generalized Wigner class without change.

 We now give further details, beginning with the proof of \eqref{eqn:max-main-1}. This theorem is proven in Section 6 of \cite{BLZ}, building on the treatment of Gaussian divisible ensembles in Section 4, and the comparison arguments of Section 5. As indicated above, Section 5 goes through for generalized Wigner matrices without change, as it is based on \cite{landon2020comparison} where the arguments were developed for generalized Wigner matrices. Section 4 requires some changes. Reading through line-by-line, one finds that the first place requiring changes is the proof of Proposition 4.3, which must be proven with the $\lambda_k (t)$ (defined near \cite[(4.1)]{BLZ}) starting from a generalized Wigner matrix instead of only a Wigner matrix. Reading through line-by-line, one can use Proposition \ref{prop:max-lss} in place of \cite[Corollary A.2]{BLZ} to obtain \cite[(4.16)]{BLZ} and the estimate $\pp[ \mathcal{G}_{\eta'} ] = 1 - o (1)$, with $\mathcal{G}_{\eta'}$ defined in \cite[(4.23)]{BLZ}, in the same manner as in the proof of \cite[Proposition 4.3]{BLZ}. Everything else is identical. The remainder of Section 4 is not needed for the proof of \eqref{eqn:max-main-1}. 
 
 With \cite[Proposition 4.3]{BLZ} extended to the case of generalized Wigner matrices, the proof of \eqref{eqn:max-main-1} given in \cite[Section 6]{BLZ} goes through without change. We only remark that first: \cite[Lemma 6.1]{BLZ} holds without change for generalized Wigner matrices.  Second, the proof of the upper bound given in \cite[Section 6.1]{BLZ} does not appear to be complete as of the writing of this paper, but the gap is extremely minor. In particular, we could not find the proof of the deterministic upper bound \cite[(6.3)]{BLZ} within the paper. Instead, one can easily check that such an estimate holds with overwhelming probability,  as long as one allows the set $J$ to have cardinality $N \e^{ C ( \log \log N)^2}$ for some large $C>0$. Indeed, the estimate \cite[(3.1)]{BLZ} allows one to add $\i /N$ to $E$. Then, estimates similar to \eqref{eqn:max-est-2} allow one to replace the supremum over a line with a supremum over the discrete set $J$. The arguments in \cite[Section 5]{BLZ} are unchanged even though the size of $J$ has increased slightly.
 
 With the upper bound proven, we turn to the proof of the lower bound in Section 6.2 of \cite{BLZ}. The proof goes through without change, except that we could not quite follow all of the arguments in \cite{BLZ}. In particular, the proof of \cite[(6.7)]{BLZ} relies on the proof of \cite[(3.21)]{BLZ} the latter, as written, relies on \cite[Lemma 2.5]{BLZ} which is proven only for $\beta$-ensembles. However one can see from the proof  of \cite[(3.21)]{BLZ} that the rigidity estimates of \cite[Theorem 2.2]{BLZ} are sufficient to obtain \cite[(6.7)]{BLZ}.
 
 For the remaining two results of Theorem \ref{thm:max} we first note that \eqref{eqn:max-main-2} and \eqref{eqn:max-main-3} are equivalent due to the fact that
 \beq
 | \{ i : \lambda_i \leq E \} | \geq k \iff \{ \lambda_k \leq E \} ,
 \eeq
 and so it suffices to prove \eqref{eqn:max-main-3}. For this, we follow directly the proof given in \cite[Section 6.3]{BLZ}. Again, \cite[Section 5]{BLZ} carries over without change to the generalized Wigner setting. The proof of \cite[(6.14)]{BLZ} follows from \cite[Proposition 4.1]{BLZ} and Proposition \ref{prop:max-lss} above. The rest of the proof is identical.  \qed

\subsection{Estimates for the CLT functionals for the logarithm}

Define,
\beq
V_{G\beta E} (f) := \frac{1}{2 \pi^2 \beta} \int_{-2}^2 \int_{-2}^2 \left( \frac{ f(x) - f(y)}{ x -y} \right)^2 \frac{ 4 - xy}{ \sqrt{4-x^2} \sqrt{4-y^2}} \d x \d y .
\eeq
\bep
Let $\kappa >0$ and let $|E| \leq 2 - \kappa$ and $0 < \eta < 1$. Define $f_\Re (x)$ and $f_\Im (x)$ by, 
\beq
f_\Re (x) := \Re[ \log (E + \i \eta - x ) ], \qquad f_\Im (x) := \Im[ \log (E + \i \eta  -x ) ]
\eeq
Then,
\beq \label{eqn:max-var-1}
V_{G\beta E} ( f_\Re  ) = \beta^{-1} | \log ( \eta) | + \O (1) = V_{G\beta E} ( f_\Im )
\eeq
For any generalized Wigner matrix we have that for $g = f_\Re$ or $g=f_\Im$ that
\beq \label{eqn:max-var-2}
V_\beta (g ) = \beta^{-1} | \log ( \eta ) | + \O (1), \qquad \E_\beta (g) = \O (1), \qquad |B(g) | \leq C N^{-1/2}. 
\eeq
\eep
\proof By \cite[Lemma A.1]{baik2016fluctuations} if we set,
\beq
L(z, w) := \int_{-2}^2 \int_{-2}^2 ( \log (z-x) - \log (z-y) )  ( \log ( w -x ) - \log ( w - y) ) Q(x, y) \d x \d y
\eeq
where,
\beq
Q(x, y) := \frac{4-xy}{(x-y)^2 \sqrt{4-x^2} \sqrt{4-y^2}}
\eeq
then
\beq
L(z, w) = 2 \pi^2 \log \left[ \frac{ (z+R(z) ) (w+ R(w) )}{2 (zw-4+R(z) R(w) ) } \right], \qquad R(z) := \sqrt{ z^2-4}
\eeq
for $z, w \in \cc \backslash (-\infty, 2]$ and $R(z)$ defined with branch cut in $[-2, 2]$. Note that,
\beq
2 \beta \pi^2 V_{G\beta E} ( \Re[ \log ( z - \cdot ) ] )  = \frac{1}{4} ( L(z, z)  + 2L(z, \bar{z} ) + L ( \bar{z}, \bar{z} ) )
\eeq
and 
\beq
2 \beta \pi^2 V_{G\beta E} ( \Im[ \log (z - \cdot ) ) = \frac{1}{4} (2 L(z, \bar{z} ) - L ( z, z) - L(\bar{z}, \bar{z} ) ).
\eeq
It is clear that $L(z, z) = \O ( 1) = L( \bar{z}, \bar{z} )$ for $|E| \leq 2 - \kappa$ and $0 < \eta < 10$. Note that,
\begin{align}
|z|^2 -4 + R(z) R(\bar{z} )  &= |z|^2 -4 + |R(z)|^2= E^2 + \eta^2 -4 + |E^2 -4 - \eta^2 + 2 \i E \eta | \notag\\
&= \frac{ 8 \eta^2}{4 - E^2} + \O ( \eta^4).
\end{align}
Hence,
\beq
\frac{1}{2} L(z, \bar{z} ) = 2 \pi^2 | \log \eta | + \O ( 1).
\eeq
The claim \eqref{eqn:max-var-1} follows. The first estimate of \eqref{eqn:max-var-2} follows from the definition of $V$ in \eqref{eqn:intro-V-def} and Lemma \ref{lem:var-3}. The estimate for $\E_\beta (f)$ follows from its definition as well as \eqref{eqn:exp-determ-4}. The estimate for $B(f)$ follows directly from its definition. \qed

\appendix

\section{Changes for complex Hermitian Wigner matrices} \label{a:cplx}

In the case of complex Hermitian Wigner matrices, we introduce the Wirtinger derivatives,for $i \neq j$ by
\beq
\del_{H_{ij}} = \frac{1}{2} \left( \del_{ \Re[ H_{ij} ] } - \i \del_{ \Im[ H_{ij} ] } \right), \qquad \del_{\bar{H}_{ij}} = \frac{1}{2} \left( \del_{ \Re[ H_{ij} ] } + \i \del_{ \Im[ H_{ij} ] } \right)
\eeq
We will denote $\del_{ij} = \del_{H_{ij}}$ so that $\del_{ji} = \del_{\bar{H}_{ij}}$. We still denote $\del_{ii} = \del_{H_{ii}}$. In this section let us denote by $\delrr_{ij}$ the derivative of a function on the space of real symmetric matrices with respect to the $(i, j)$th element.  The main algebraic difference between the complex Hermitian and real symmetric cases is that,
\beq \label{eqn:derivative-complex}
\del_{ij} G_{ab} (z) = - G_{ai} (z) G_{j b } (z)
\eeq
instead of
\beq \label{eqn:derivative-real}
(1+ \delta_{ij} ) \delrr_{ij} G_{ab} (z) = - G_{ai}(z) G_{jb}(z) - G_{aj}(z) G_{ib}(z) .
\eeq
Some general observations are then that any kind of monomial of Green's function elements that arises in the complex Hermitian case already has arisen in the real symmetric case. For example, in estimating third order contributions to cumulant expansions in the real case, we had to account for every term in an expression like $(\delrr_{ij})^2 G_{ab} (z)$. In the complex case, the third derivatives will be combinations of $\del_{ij} \del_{ji} G_{ab}$, $\del_{ij}^2 G_{ab}$ and $\del_{ji}^2 G_{ab}$. These three expressions all contribute monomials in Green's function elements that already arose in the real case $(\delrr_{ij})^2 G_{ab} (z)$, and so the work we have already done will apply without much change to the complex Hermitian case.

The main changes to the complex Hermitian case will then be that not every second order term that arose in the real symmetric case will arise in the complex Hermitian case. For example, when computing $\ee[ G_{11}(z)]$ via the cumulant expansion, the second order terms in the real symmetric and complex Hermitian cases will be,
\beq  \label{eqn:derivative-compare}
\delrr_{1j} G_{j1}(z) = - G_{jj} (z) G_{11}(z) - G_{1j}(z)^2, \qquad \del_{j1} G_{j1} (z) = - G_{jj}(z) G_{11}(z),
\eeq
respectively. So some second order terms present in the real symmetric case will not be there in the complex Hermitian case, simplifying the argument somewhat.

The only thing that one needs to  be careful about is that in the real symmetric case one has $G_{ij} (z) = G_{ji} (z)$ but in the complex Hermitian case, $G_{ij}(\bar{z}) = \widebar{ G_{ji} (z ) }$. We have already accounted for this in our notion of hypergraphs, in that we assume that the graphs are directed. In this way there is no ambiguity in the resolvent entry associated to an edge in our hypergraph (i.e., if $e = (a, b)$ then the associated resolvent entry is $G_{i_a i_b} (z)$ and not $G_{i_b i_a} (z)$). Of course the estimates \eqref{eqn:entry-wise} and \eqref{eqn:iso} are insensitive to whether we consider $G_{ij} (z)$ or $G_{ji} (z)$. 

However, there is one crucial aspect in a few of our proofs (in particular, Proposition \ref{prop:main-loop-estimate}) for which the order of indices matters. We sometimes use the fact that $\sum_{k} G_{ik}(z) G_{ka}(w) = \langle \delta_a , (H-z)^{-1} (H-w)^{-1} \delta_b \rangle$ crucially. We need to check that the terms that arise in our expansions are of the form $G_{ik} (z) G_{ka} (w)$ and not $G_{ik} (z) G_{ak} (w)$ as the same trick would not work for the latter monomials.

Let us now explain the underlying reason why the monomials we consider always have their indices appearing in the correct order. These terms all arise from second order terms in our cumulant expansions. The operation that generates the second order terms, differentiation by $\del_{ij}$, has good algebraic properties that will not change the order that the indices appear in, and so our method will nonetheless work. I.e., if $F$ is some monomial in resolvent entries, then the resolvent and cumulant expansion is schematically (keeping only the dominant second order terms)
\beq
G_{ab} F \to H_{aj} G_{jb} F \to \del_{ja} (G_{jb} F ) .
\eeq
Then, due to \eqref{eqn:derivative-complex}, if $\del_{ja}$ hits any resolvent entry, then it adds $j$ as a right index and $a$ as a left index. Therefore, this expansion does not change whether $a$ or $b$ appear as right or left indices and introduces a new index $j$ that appears once as a left index and once as a right index. Due to \eqref{eqn:delea-cplx} below, a similar observation applies to when the derivative hits $\ea$.

In what follows we will go through the arguments of the main body of the paper, detailing the necessary changes. As is clear from the above discussion, for the most part the arguments needed to handle the complex Hermitian case are the same as the real symmetric cases and so little comment will be necessary. We mainly include this appendix for completeness, as there is little to no novelty introduced.

In the complex Hermitian case we fix the following notation,
\beq
\frac{ s_{ij}}{N} := \ee[ | H_{ij}|^2], \qquad  s_{ij}^{(k)} := \kappa_k ( \sqrt{N} \Re[ H_{ij} ] ), \qquad  t_{ij}^{(k)} := \kappa_k ( \sqrt{N} \Im[ H_{ij} ] ) ,
\eeq
for $k \geq 3$. With this notation, by appling the cumulant expansion  of Lemma \ref{lem:cumu-exp} twice (once to the real part of $H_{ij}$ and once to the imaginary part) we find that,
\begin{align} \label{eqn:cumu-complex}
\ee[ H_{ij} F(H)] =&  \frac{ s_{ij}}{N} \ee[ \del_{ji} F(H) ] \notag\\
+&  \sum_{k=3}^K \frac{ s_{ij}^{(k)}}{(k-1)! N^{k/2}(1+ \delta_{ij} )^{k-1}} \ee[( \del_{ij} + \del_{ji} )^{k-1} F(H) ]  \notag\\
+& \i \sum_{k=3}^K \frac{ t_{ij}^{(k)}}{(k-1)! N^{k/2}} \ee[ ( \i ( \del_{ij}  - \del_{ji} ) )^{k-1} F(H) ] + \mbox{error}
\end{align}
with an error term having a similar form as in Lemma \ref{lem:cumu-exp}.

We first note that in Lemma \ref{lem:delea}, the equation \eqref{eqn:deleaij} is replaced by,
\beq \label{eqn:delea-cplx} 
\del_{ij} \ea = - \frac{ \i \lambda \ea}{ \pi} \int_{\Oma} ( \del_{\bar{z}} \tilf (z) ) \del_z G_{ji} (z) \d z \d \bar{z} = \Osd ( | \lambda \Phi | )
\eeq
and \eqref{eqn:deleaij2} is replaced by,
\beq \label{eqn:delea-cplx-2}
\del_{ij}^{2-n} \del_{ji}^n \ea = \Osd ( |\lambda \Phi| + |\lambda|^2 \Phi^2 ) + \1_{ \{ n=1\}} \lambda c_f .
\eeq
The proof is very similar.  Moreover,  the estimate in \eqref{eqn:del-gen-1} holds with the LHS  replaced by $| \del_{ij}^{n-m} \del_{ji}^m \ea| $ for any $0 \leq m \leq n$. Again, the proof is similar.

\subsection{Loops estimate}

We now detail the necessary changes in the complex Hermitian case in the arguments of Section \ref{sec:loops}. Definition \ref{def:loop-admissible} and the main result, Proposition \ref{prop:main-loop-estimate} are unchanged. The statement of Lemma \ref{lem:del-loops-1} changes slightly as the estimates for the quantity \eqref{eqn:del-loops-1} hold for any term of the form,
\beq \label{eqn:del-loops-1-cplx}
\del_{1j}^{k-n} \del_{j1}^n ( G_{j2} (z) M_1 G_{a1} (w) M_2 )
\eeq
for any $0 \leq n \leq k$. The proof is almost identical (the order of derivative referred to in \ref{it:gen-der-4} is in this case the power $k$ in the expression \eqref{eqn:del-loops-1-cplx}).  Similarly, the statement of Lemma \ref{lem:loops-fourth} is unchanged except that the  $\del_{1j}^3$ in \eqref{eqn:loops-fourth} is replaced by $\del_{1j}^{3-n} \del_{j1}^n$ for $0 \leq n \leq 3$ and the estimates hold for any such $n$.

In \eqref{eqn:loops-expand-1} of Lemma \ref{lem:loops-expand-1}, the third order term on the second line is replaced by,
\begin{align}  \label{eqn:loops-expand-cplx}
& \sum_j \frac{ s_{1j}^{(3)}}{2 N^{3/2}(1+ \delta_{1j} )^2} \ee[( \del_{1j} + \del_{j1} )^{2} (\ea G_{j2} (z) M_1 G_{a1} (w) M_2) ] \notag\\
- & \i \sum_j \frac{ t_{1j}^{(3)}}{2 N^{3/2}} \ee[( \del_{1j} - \del_{j1} )^{2} (\ea G_{j2} (z) M_1 G_{a1} (w) M_2) ] ,
\end{align}
due to \eqref{eqn:cumu-complex}.

The estimates in Lemma \ref{lem:loops-third} for the third order terms must then be modified as follows. First, the estimates \eqref{eqn:loops-third-1} and \eqref{eqn:loops-third-3} hold with the $\del_{1j}^2$ on the LHS replaced by any of $\del_{1j}^2, \del_{1j} \del_{j1}$ or $\del_{j1}^2$. Similarly, the estimate \eqref{eqn:loops-third-2} holds with any number of the $\del_{1j}$ appearing on the LHS replaced by $\del_{j1}$. Again, the proof is almost identical.

%
%
%
%
%

After having addressed the third order terms, we now turn to the computation of the second order terms in the complex analog of \eqref{eqn:loops-expand-1}. First, Lemma \ref{lem:loops-second-1} is unchanged. For Lemma \ref{lem:loops-second-2}, we have that the term with $D_2$ on the last line of \eqref{eqn:loops-second-2} is not present, but the estimate \eqref{eqn:loops-second-2} is otherwise unchanged. The difference is due to the fact that the second term on the RHS of \eqref{eqn:loops-second-2b} does not appear in the complex Hermitian case. In Lemma \ref{lem:loops-second-3}, we have that the term with $D_4$ is not present in \eqref{eqn:loops-second-3}. Again, this is due to the action of operator $\del_{j1}$ in the complex case versus the real case (i.e., \eqref{eqn:derivative-complex} vs \eqref{eqn:derivative-real}). Finally, the estimate \eqref{eqn:loops-second-4} of Lemma \ref{lem:loops-second-4} is unchanged except the prefactor $\frac{2 \i \lambda}{ \pi}$ is replaced by $\frac{ \i \lambda}{ \pi}$ due to \eqref{eqn:delea-cplx}.  This completes the discussion of the second order terms in \eqref{eqn:loops-expand-1}.

The outcome of all the above discussion is then that in \eqref{eqn:loops-expand-2} of Proposition \ref{prop:loops-expand}, the terms with $D_2$ and $D_4$ do not appear, and that the $2 \i \lambda$ is replaced by $\i \lambda$ in the last line. The rest of Section \ref{sec:loop-self} is then unchanged, after accounting for this difference. The arguments of Section \ref{sec:loop-iterative} and the proof of Proposition \ref{prop:main-loop-estimate} given in Section \ref{sec:main-loops-proof} then go through without change.

We now turn to detailing how Section \ref{sec:short-loops} changes in the complex case. First, the statement of Lemma \ref{lem:weakgij-expand} holds without change, given the changes to \eqref{eqn:loops-expand-1} and Lemma \ref{lem:loops-third} we have outlined above. Next, the estimates \eqref{eqn:weakgij-second-1}, \eqref{eqn:weakgij-second-2} and \eqref{eqn:weakgij-second-3} hold without change, via similar proofs. The estimate \eqref{eqn:weakgij-second-4} also holds without change except that in the proof, some terms estimated using Proposition \ref{prop:main-loop-estimate} in the real symmetric case are simply not present in the complex case (and so the proof is easier). The remainder of the arguments of Section \ref{sec:short-loops} go through without change and in particular the estimate \eqref{eqn:weak-g12g21g33} holds without change in the complex case.

\subsection{Line estimates}

We now detail the changes to Section \ref{sec:line} in the complex case. First, Definition \ref{def:line-admissible} and the main result of this section, Proposition \ref{prop:main-line-estimate} are unchanged.  For Lemma \ref{lem:line-der}, the estimates of the term \eqref{eqn:line-der} hold with $\del_{j1}^k$ replaced by $\del_{1j}^{k-n} \del_{j1}^n$ for any $0 \leq n \leq k$. The proof is similar.  The estimates of Lemma \ref{lem:line-fourth} hold for the quantity \eqref{eqn:line-fourth} with the $\del_{j1}^3$ on the LHS replaced by $\del_{1j}^{3-n} \del_{j1}^n$ for any $0 \leq n \leq 3$, and the proof is similar.

Now, in the expansion \eqref{eqn:line-expansion-1} of Lemma \ref{lem:line-expand-1}, the third order terms on the second line are replaced by,
\begin{align}
& \sum_j \frac{ s_{1j}^{(3)}}{2 N^{3/2}(1+ \delta_{1j} )^2} \ee[( \del_{1j} + \del_{j1} )^{2} (\ea G_{j2} (z) M_1  M_2) ] \notag\\
- & \i \sum_j \frac{ t_{1j}^{(3)}}{2 N^{3/2}} \ee[( \del_{1j} - \del_{j1} )^{2} (\ea G_{j2} (z) M_1  M_2) ] .
\end{align}
For the third order terms, the estimates \eqref{eqn:line-third-1},\eqref{eqn:line-third-2} \eqref{eqn:line-third-3}, and \eqref{eqn:line-third-4} of Lemma \ref{lem:line-third} all hold, with any of the $\del_{j1}$ replaced by $\del_{1j}$. The proofs are almost identical. 

For the second order terms, the estimates of Lemma \ref{lem:line-second-1} hold without change. The estimate \eqref{eqn:line-second-2} of Lemma \ref{lem:line-second-2} holds except that the term with $D_2$ on the second line is not present in the complex case. This is due to the fact that the second term on the last line of \eqref{eqn:line-second-2a} is not present in the complex case.  In Lemma \ref{lem:line-second-3}, the equation \eqref{eqn:line-second-3} changes as the term with $D_4$ in the second line is not present. In \eqref{eqn:line-second-4}, the prefactor $\frac{ 2 \i \lambda}{ \pi}$ is replaced by $\frac{ \i \lambda}{ \pi}$.  Finally, the estimate \eqref{eqn:line-exx} of Proposition \ref{prop:line-expand} is unchanged except that the terms with $D_2$ and $D_4$ are not present, and the $ 2 \i \lambda$ in the last line is replaced by $ \i \lambda$. 

With the above, the arguments of Section \ref{sec:line-iterative} and the proof of Proposition \ref{prop:main-line-estimate} given in Section \ref{sec:main-line-proof} are then unchanged.

We now turn to outlining how Section \ref{sec:gij-better} changes. First, the result of Lemma \ref{lem:gij-expand} is unchanged, but the handling of the third order terms (i.e., the analog of the first term on the second line of \eqref{eqn:gij-expand-a}) changes somewhat. Specifically, one uses $\del_{1j}^{2-n} \del_{j1}^n G_{j2} (z) = \1_{ \{ n=1\}} \msc(z)^2 G_{j2} (z) + \Osd ( \Psi(z)^2)$ for $0 \leq n \leq 2$, and then an argument similar to \eqref{eqn:gij-expand-b} yields \eqref{eqn:gij-expand}. 

All of the three estimates of Lemma \ref{lem:gij-second} are unchanged. Note that in the proof of \eqref{eqn:g12-second-3}, the second term on the second line of \eqref{eqn:gij-second-3a} is not present in the complex case. Finally, Proposition \ref{prop:gij-est} follows in the complex case as in the real case.

\subsection{Preliminary expansion for $\ee[ \ea G_{11} ]$}

We now turn to the necessary changes to the arguments of Section \ref{sec:prelim} in the complex case. First, the estimates of Lemma \ref{lem:G11-fifth} hold, with \eqref{eqn:G11-fifth-der-2} also holding when $\del_{j1}^4$ on the LHS is replaced by $ \del_{j1}^{4-n} \del_{1j}^n$ for any $0 \leq n \leq 4$.

For the estimate \eqref{eqn:prelim-expand-1}, the third order terms are replaced with,
\begin{align} \label{eqn:cplx-c2}
&\sum_j \frac{ s_{1j}^{(3)}}{2 N^{3/2}(1+ \delta_{1j} )^2} \ee[( \del_{1j} + \del_{j1} )^{2} (\ea G_{j1} (z)) ] 
-  \i \sum_j \frac{ t_{1j}^{(3)}}{2 N^{3/2}} \ee[( \del_{1j} - \del_{j1} )^{2} (\ea G_{j1} (z) ) ] 
\end{align}
and the fourth order terms are replaced by,
\begin{align} \label{eqn:cplx-c1}
&\sum_j \frac{ s_{1j}^{(4)}}{6 N^{2}(1+ \delta_{1j} )^3} \ee[( \del_{1j} + \del_{j1} )^{3} (\ea G_{j1} (z)) ] +   \sum_j \frac{ t_{1j}^{(4)}}{6 N^{2}} \ee[( \del_{1j} - \del_{j1} )^{3} (\ea G_{j1} (z) ) ] .
\end{align}
We now turn to Section \ref{sec:prelim-fourth}, which deals with the fourth order terms in \eqref{eqn:cplx-c1}. First, the estimate \eqref{eqn:prelim-fourth-1} is replaced by, 
\begin{align} 
& \sum_j \frac{ s_{1j}^{(4)}}{N^2(1+ \delta_{1j})^3} \ee[ ( \del_{1j} + \del_{j1} )^3 ( \ea G_{j1} ) ] +\sum_j \frac{ t_{1j}^{(4)}}{ N^{2}} \ee[( \del_{1j} - \del_{j1} )^{3} (\ea G_{j1} (z) ) ]   \notag\\
= & - \msc(z)^4 \sum_{j \neq 1} \frac{s_{1j}^{(4)}+ t_{1j}^{(4)}}{N^2} 6 \ee[ \ea ] - 12 \msc(z)^3 \sum_{j \neq 1 } \frac{ s_{1j}^{(4)} + t_{1j}^{(4)}}{N^2} \ee[ \ea ( G_{11} (z)  - \msc(z) ) ]  \notag\\
+ & 3 \sum_{j \neq 1} \frac{ s_{1j}^{(4)}}{N^2} \ee[( ( \del_{1j} + \del_{j1} )^2 \ea) (\del_{1j} + \del_{j1} ) G_{j1} ]+3 \sum_{j \neq 1} \frac{ t_{1j}^{(4)}}{N^2} \ee[( ( \del_{1j} - \del_{j1} )^2 \ea) (\del_{1j} - \del_{j1} ) G_{j1} ] \notag\\
+ & \Osd ( N^{-2} \eta^{-1}  + ( 1 + |\lambda| )^3 N^{-2} + |\lambda \Phi| \Psi(z) (1 + |\lambda| ) .
\end{align}
The proof is the same, except the identity 
\beq
\del_{j1}^{3-n} \del_{1j}^n G_{j1} (z) = - 2 G_{jj}(z)^2 G_{11}(z)^2 \1_{ \{ n=1 \} } + \Osd ( \Psi(z)^2 )\eeq
 replaces \eqref{eqn:prelim-fourth-1a}. Instead of \eqref{eqn:prelim-fourth-2} of  Lemma \ref{lem:prelim-fourth-2} we derive the estimate,
 \begin{align}
 & \sum_{j \neq 1} \frac{ s_{1j}^{(4)}}{N^2} \ee[( ( \del_{1j} + \del_{j1} )^2 \ea) (\del_{1j} + \del_{j1} ) G_{j1} ]+ \sum_{j \neq 1} \frac{ t_{1j}^{(4)}}{N^2} \ee[( ( \del_{1j} - \del_{j1} )^2 \ea) (\del_{1j} - \del_{j1} ) G_{j1} ] \notag\\
 = & - \ee[\ea] \msc(z)^2 \sum_{j \neq 1 } \frac{  s_{1j}^{(4)}- t_{1j}^{(4)}}{N^2} \frac{ 2 \i \lambda}{ \pi} \int_{\Oma} ( \del_{\bar{u}} \tilf (u ) ) \del_u \msc(u)^2 \d u \d \bar{u} \notag\\
 &-  \sum_{j \neq 1} \frac{ s_{1j}^{(4)} -t_{ij}^{(4)}}{N^2} \frac{2 \i \lambda}{\pi} \int_{\Oma} ( \del_{\bar{u}} \tilf (u) ) \del_u \msc(z)^2 \msc(u) \ee[ \ea(  G_{11} (u) - \msc(u) ) ] \notag\\ 
  &- \sum_{j \neq 1} \frac{ s_{1j}^{(4)}- t_{ij}^{(4)}}{N^2}\frac{2 \i \lambda}{\pi} \int_{\Oma} ( \del_{\bar{u}} \tilf (u) ) \del_u \msc(u)^2 \msc(z) \ee[ \ea(  G_{11} (z)  -\msc(z) ) ] \notag\\ 
  + &  \Osd ( |\lambda|N^{-2} \eta^{-1} + |\lambda|(1+ |\lambda| ) \Phi^2 N^{-1} )  .
 \end{align}
 The proof is the same, except that we use the estimates
 for $n=0, 1$,
\beq \label{eqn:prelim-cplx-1}
\del_{j1}^{1-n} \del_{1j}^n G_{j1} (z) = - \1_{ \{ n=1 \} } G_{jj} G_{11} + \Osd ( \Psi(z)^2 )
\eeq
and for $m=0, 1, 2$,
\beq \label{eqn:prelim-cplx-2}
\del_{j1}^{2-m} \del_{1j}^m \ea = \1_{ \{ m =1 \} } \frac{ \i \lambda}{ \pi} \int_{\Oma} ( \del_{\bar{u}} \tilf (u) ) \del_u (G_{11} (u) G_{jj} (u) ) + \Osd ( | \lambda \Phi^2 | ( 1+ | \lambda| ) ) 
\eeq
instead of \eqref{eqn:prelim-fourth-2a} and \eqref{eqn:prelim-fourth-2b}.  This completes our outline of the changes to the fourth order terms in \eqref{eqn:cplx-c1}

We now turn to outlining the changes to the handling of the third order terms \eqref{eqn:cplx-c2} in Section \ref{sec:prelim-third}. First, estimate \eqref{eqn:prelim-third-4} of Lemma \ref{lem:prelim-third-1} holds without change. Next, estimate \eqref{eqn:prelim-third-1} of Lemma \ref{lem:prelim-third-2} changes to,
\begin{align}
&\sum_{j\neq 1} \frac{ s_{1j}^{(3)}}{ N^{3/2}} \ee[\ea ( \del_{1j} + \del_{j1} )^{2} G_{j1} (z) ] 
-  \i \sum_{j \neq 1} \frac{ t_{1j}^{(3)}}{ N^{3/2}} \ee[\ea ( \del_{1j} - \del_{j1} )^{2}  G_{j1} (z) ]  \notag\\
= &  \msc (z)\sum_{j \neq 1 } \frac{ 1}{N^{3/2}} \ee[ \ea ( G_{jj} -\msc(z) )(2(s_{ij}^{(3)} - \i t_{ij}^{(3)} ) G_{1j}+4 (s_{ij}^{(3)} + \i t_{ij}^{(3)} )G_{j1} ) ] \notag\\
+ & \Osd ( N^{-1/2} \Psi(z)^3 + |\lambda \Phi | N^{-1} \Psi(z) ).
\end{align}
The proof is the same except that we use
\beq
\del_{1j}^{2-n} \del_{j1}^n G_{j1} (z) = \1_{ \{ n=2\}} 2 G_{jj} G_{11} G_{1j} + \1_{ \{ n =1 \} } 2 G_{jj} G_{11} G_{j1} + \Osd ( \Psi(z)^3)
\eeq
instead of \eqref{eqn:prelim-third-1a} and we also use Proposition \ref{prop:gij-est} to estimate both $\ee[ \ea G_{1j} (z) ] $ and $\ee[ \ea G_{j1} (z) ]$. 

Instead of \eqref{eqn:prelim-third-2} of Lemma \ref{lem:prelim-third-3}  we have,
\begin{align}
&\sum_{j\neq 1} \frac{ s_{1j}^{(3)}}{ N^{3/2}} \ee[(( \del_{1j} + \del_{j1} )\ea) ( \del_{1j} + \del_{j1} ) G_{j1} (z) ] 
-  \i \sum_{j \neq 1} \frac{ t_{1j}^{(3)}}{ N^{3/2}} \ee[(( \del_{1j} - \del_{j1} )\ea) ( \del_{1j} - \del_{j1} ) G_{j1} (z) ]  \notag\\
&= \sum_{ j \neq 1 } \frac{ \msc(z)}{N^{3/2}} \ee[  ( (s_{1j}^{(3)}(\del_{1j}+ \del_{j1} ) + \i t_{1j}^{(3)} ( \del_{1j} - \del_{j1} )) \ea ) ( G_{jj} (z)-  \msc (z) ) ] \notag\\
+ & \Osd ( N^{-1/2} |\lambda \Phi | \Psi(z)^2 + N^{-1} |\lambda| ( 1 + |\lambda| ) \Phi^2 )
\end{align}
via almost the same proof, using the identity \eqref{eqn:prelim-cplx-1}.

Instead of \eqref{eqn:prelim-third-3} of Lemma \ref{lem:prelim-third-4} we have,
\begin{align}
& \sum_{j\neq 1} \frac{ s_{1j}^{(3)}}{ N^{3/2}} \ee[(( \del_{1j} + \del_{j1} )^2\ea)  G_{j1} (z) ] 
-  \i \sum_{j \neq 1} \frac{ t_{1j}^{(3)}}{ N^{3/2}} \ee[(( \del_{1j} - \del_{j1} )^2\ea)  G_{j1} (z) ]  \notag\\ 
= & \frac{2 \i \lambda}{\pi} \sum_{ j \neq 1 } \frac{ s_{1j}^{(3)}+ \i t_{1j}^{(3)}}{N^{3/2}} \int_{\Oma} \d u \d \bar{u} ( \del_{ \bar{u}} \tilf (u) ) \del_u \msc (u) \ee[ \ea G_{j1} (z) ( G_{jj} (u) -  \msc (u) ) ] \notag\\
+ & \Osd ( N^{-1/2} |\lambda \Phi^2| ( 1+ |\lambda|) \Psi (z) + |\lambda| N^{-1} \Psi (z)^2 )
\end{align}
where we use \eqref{eqn:prelim-cplx-2} instead of \eqref{eqn:prelim-fourth-2b}. This completes the changes to the third order terms in \eqref{eqn:cplx-c1}.

We now turn to detailing the changes in Section \ref{sec:intermediate}. First, the estimate \eqref{eqn:intermediate-expansion} of Proposition \ref{prop:intermediate-expansion} becomes,
\begin{align}
& z  \sum_a \ee[ \ea G_{aa} (z) + N \ee[ \ea ] = \sum_{ja} \frac{ s_{1j}}{N} \ee[ \del_{ja} ( \ea G_{ja} ) ] \notag \\
&  -  \sum_{j \neq a } \frac{ s_{aj}^{(4)} + t_{aj}^{(4)}}{N^2} \msc(z)^4 \ee[ \ea ] - \ee[ \ea] \msc(z)^2 \sum_{j \neq a } \frac{ s_{aj}^{(4)} - t_{aj}^{(4)}}{N^2} \frac{\i \lambda}{\pi} \int_{\Oma} ( \del_{\bar{u}} \tilf (u) ) \del_u \msc(u)^2 \d u \d \bar{u}  \notag\\
&+ \sum_a  \frac{ s_{aa}^{(3)}}{2 N^{3/2}} \ee[ \ea] \bigg\{ \msc(z)^3 +  2 \msc(z)^2 \left( \frac{ \i \lambda}{\pi} \ea \int_{\Oma} \d u \d \bar{u} \del_{\bar{u}} \tilf (u) \del_u \msc(u) \right)  \notag\\
+ & \msc(z) \left( - \frac{ \i \lambda}{\pi}  \int_{\Oma} \d u \d \bar{u} \del_{\bar{u}} \tilf (u) \del_u \msc(u) \right)^2 
+ \msc(z)  \left( \frac{ \i \lambda}{\pi}  \int_{\Oma} \d u \d \bar{u} \del_{\bar{u}} \tilf (u) \del_u \msc(u)^2 \right) \bigg\} \notag\\
+ &\Osd ( N^{-1} \eta^{-1} (1 + |\lambda| ) + ( 1 + |\lambda| )^3 N^{-1} +N^{1/2} \Psi(z)^3 + N^{1/2} |\lambda \Phi | \Psi(z)^2) \notag \\
+ &  \Osd ( N^{1/2} |\lambda \Phi^2| ( 1+ |\lambda|) \Psi (z)  )  
\end{align}
by using the above changes to the various lemmas of Sections \ref{sec:prelim-fourth} and \ref{sec:prelim-third}. As in the proof of Proposition \ref{prop:intermediate-expansion} in the real case, there are some extra terms that we can handle by \eqref{eqn:fluct-av} in the exact same manner.

  We then find that \eqref{eqn:intermediate-characteristic} is replaced by,
\begin{align} \label{eqn:intermediate-characteristic-cplx}
&z  \sum_a \ee[ \ea (G_{aa} (z) - \ee[ G_{aa} (z) ] ) ] = \sum_{ja} \frac{ s_{aj}}{N} \left( \ee[ \del_{ja} ( \ea G_{ja} (z) ) ] - \ee[ \ea] \ee[ \del_{ja} G_{ja} (z) ] \right) \notag\\
& - \ee[ \ea] \msc(z)^2 \sum_{j \neq a } \frac{ s_{aj}^{(4)}+ t_{aj}^{(4)}}{N^2} \frac{ \i \lambda}{\pi} \int_{\Oma} ( \del_{\bar{u}} \tilf (u) ) \del_u \msc(u)^2 \d u \d \bar{u}  \notag\\
&+ \sum_a  \frac{ s_{aa}^{(3)}}{2 N^{3/2}} \ee[ \ea] \bigg\{   2 \msc(z)^2 \left( \frac{ \i \lambda}{\pi} \ea \int_{\Oma} \d u \d \bar{u} \del_{\bar{u}} \tilf (u) \del_u \msc(u) \right)  \notag\\
+ &  \msc(z) \left( - \frac{ \i \lambda}{\pi}  \int_{\Oma} \d u \d \bar{u} \del_{\bar{u}} \tilf (u) \del_u \msc(u) \right)^2 
+  \msc(z)  \left( \frac{ \i \lambda}{\pi}  \int_{\Oma} \d u \d \bar{u} \del_{\bar{u}} \tilf (u) \del_u \msc(u)^2 \right) \bigg\}  \notag\\
+ &\Osd ( N^{-1} \eta^{-1} (1 + |\lambda| ) + ( 1 + |\lambda| )^3 N^{-1} +N^{1/2} \Psi(z)^3 + N^{1/2} |\lambda \Phi | \Psi(z)^2) \notag \\
+ &  \Osd ( N^{1/2} |\lambda \Phi^2| ( 1+ |\lambda|) \Psi (z)  )
\end{align}
and that \eqref{eqn:intermediate-expectation} is replaced by,
\begin{align} \label{eqn:intermediate-expectation-cplx}
& z   \sum_a \ee[ G_{aa} (z) ]  +N= \sum_{ja} \frac{ s_{aj}}{N} \ee[ \del_{ja} G_{ja} (z) ]  \notag\\
-&   \sum_{j \neq a } \frac{ s_{aj}^{(4)} + t_{aj}^{(4)}}{N^2} \msc(z)^2 + \sum_{a} \frac{ s_{aa}^{(3)}}{2 N^{3/2}} \msc(z)^3 + \Osd ( N^{-1} \eta^{-1} + N^{1/2} \Psi(z)^3 ) .
\end{align}

We now turn to outlining how Section \ref{sec:g11-better} changes in the complex case. First, the estimate \eqref{eqn:g11-expand-1} of Proposition \ref{prop:g11-expand-1} holds without change. Indeed, the third and fourth order terms in the cumulant expansion are estimated using the estimates listed above, and the remainder of the terms are estimated in the same fashion as the proof of Proposition \ref{prop:g11-expand-1}. The estimates \eqref{eqn:g11-expand-2} and \eqref{eqn:g11-expand-3} of Lemma \ref{lem:g11-second-1} hold without change. Given these inputs, the rest of Section \ref{sec:g11-better} then holds without change.

\subsection{Estimates for $G_{12} G_{21}$}

We now outline the changes to Section \ref{sec:g12g21} in the complex case. First, the statement of Lemma \ref{lem:g12g21-expand} holds without change, but the proof changes somewhat. First, the third order term on the second line of \eqref{eqn:g12g21-third-1b} becomes
\beq
\sum_{j \neq 1 } \frac{ s_{1j}^{(3)}}{2N^{3/2}} \ee[ ( \del_{1j} + \del_{j1} )^2 ( \ea G_{j2} (z) G_{21} (w) ) ] - \i \sum_{j \neq 1 } \frac{ t_{1j}^{(3)}}{2N^{3/2}} \ee[ ( \del_{1j} - \del_{j1} )^2 ( \ea G_{j2} (z) G_{21} (w) ) ]
\eeq
In order to estimate this term, we then use, instead of \eqref{eqn:g12g21-third-1a}, the estimate
\begin{align}
\del_{1j}^{2-n} \del_{j1}^n G_{j2} (z) &= \1_{ \{ n=1\}} \msc(z)^2 G_{j2} (z) + \Osd ( \Psi(z)^2)   \notag\\
\del_{1j}^{2-n} \del_{j1}^n G_{21} (w) &= \1_{ \{ n=1\}} \msc(w)^2 G_{21} (w) + \Osd ( \Psi(w)^2) \notag\\
\del_{1j}^{1-m} \del_{j1}^m G_{j2} (z) &= -\1_{ \{ m =1 \} } \msc(z) G_{12} (z) + \Osd ( \Psi(z)^2 )\notag\\
\del_{1j}^{1-m} \del_{j1}^m G_{21} (w) &= -\1_{ \{ m =1 \} } \msc(w) G_{2j} (w) + \Osd ( \Psi(w)^2)
\end{align} 
The rest of the proof is similar (using at one point \eqref{eqn:delea-cplx-2} instead of \eqref{eqn:deleaij2}).  Next, the estimates of Proposition \ref{prop:g12g21-second} all hold without change. In fact, the proofs are slightly easier as some terms we estimated in the real case are not present in the complex case (e.g., compare the two derivatives in \eqref{eqn:derivative-compare}). After this, the remainder of the arguments of Section \ref{sec:g12g21} do not involve any differentiation with respect to matrix entries and so the arguments are identical. 

\subsection{Estimates for $G_{11}G_{22}$}

We now outline the changes to Section \ref{sec:g11g22} in the complex case. First, for Lemma \ref{lem:giigjj-fourth}, the estimate \eqref{eqn:giigjj-fourth-2} is unchanged. The estimates \eqref{eqn:giigjj-fourth-1} and \eqref{eqn:giigjj-fourth-3} hold for any of the powers of $\del_{j1}$ on the LHS replaced by $\del_{1j}$; the proofs are similar. In \eqref{eqn:giigjj-fourth-4}, the third order term (the last term on the RHS of the first line) is replaced by,
\begin{align}
&\sum_j \frac{ s_{1j}^{(3)}}{2 N^{3/2}(1+ \delta_{1j} )^2} \ee[( \del_{1j} + \del_{j1} )^{2} (\ea G_{j1} (z) X_1) ] 
-  \i \sum_j \frac{ t_{1j}^{(3)}}{2 N^{3/2}} \ee[( \del_{1j} - \del_{j1} )^{2} (\ea G_{j1} (z) X_1 ) ].
\end{align} 
For the estimates of the third order terms carried out in Lemma \ref{lem:giigjj-third}, we have first that \eqref{eqn:giigjj-third-1} is unchanged. The three remaining estimates hold also for any of the $\del_{j1}$ on the LHS replaced by $\del_{1j}$, and the proofs are modified as follows. For \eqref{eqn:giigjj-third-2}, the proof uses \eqref{eqn:delea-cplx-2}. The proof of \eqref{eqn:giigjj-third-3} uses 
\beq
\del_{j1}^{1-n} \del_{1j}^n G_{j1} (z) = - \1_{ \{ n=1 \} } \msc(z) + \Osd ( \Psi(z) )
\eeq
for $n=0, 1$. The proof of \eqref{eqn:giigjj-third-4} also uses this identity to reduce the estimate to estimating the first term on the RHS of \eqref{eqn:giigjj-third-4a} (note that this term is only present if a factor of $\del_{j1}$ hits $G_{j1} (z)$). In this case, we need to estimate the first term on the RHS of \eqref{eqn:giigjj-third-4a} also with $\del_{j1} X_1$ replaced by $\del_{1j} X_1$. The argument in \eqref{eqn:giigjj-third-4b} works for both of these cases, yielding the estimate.

All of the estimates of Lemma \ref{lem:giigjj-second} are unchanged. The proofs are similar, except that some terms that were estimated in the real case are simply not present in the complex case. Finally, Proposition \ref{prop:giigaa-est} and its proof are unchanged.

We now outline the changes to Section \ref{sec:expectation}, in which we compute the correction to the expectation of a linear spectral statistic. First,  recall that in the complex case, the expansion \eqref{eqn:intermediate-expectation} is replaced by \eqref{eqn:intermediate-expectation-cplx} above.  Next, the estimates of Lemma \ref{lem:exp-second} change somewhat. First, the two terms on the second line of \eqref{eqn:exp-second-1} involving $S$ are not present in the complex case. This is due to the fact that the second term on the RHS of \eqref{eqn:exp-second-1a} are not present in the complex case; the rest of the proof is the same. Consequently, \eqref{eqn:exp-second-2} is replaced by,
\begin{align}
&\sum_a \ee[ G_{aa}(z) - \msc(z) ] =  \sum_{j \neq a } \frac{ s_{aj}^{(4)} + t_{aj}^{(4)}}{N^2} \msc(z)^3 \msc'(z) - \sum_a \frac{s_{aa}^{(3)}}{2 N^{3/2}} \msc(z)^2 \msc'(z) \notag\\
&  + \Osd (| \msc' (z) | N^{-1} \eta^{-2} ),
\end{align}
where we note that the fourth order term also changed because \eqref{eqn:intermediate-expectation}  changed to \eqref{eqn:intermediate-expectation-cplx}.  Consequently, the equations of Lemma \ref{lem:exp-determ} change. First, one should define $E(z)$ in \eqref{eqn:exp-determ-0} as
\beq
E_\cc (z) := \sum_{j \neq a } \frac{ s_{aj}^{(4)} + t_{aj}^{(4)}}{N^2} \msc(z)^3 \msc' (z) - \sum_a \frac{ s_{aa}^{(3)}}{2 N^{3/2}} \msc(z)^2\msc' (z).
\eeq
Then, \eqref{eqn:exp-determ-1} holds with $E(z)$ replaced by $E_\cc (z)$. Equation \eqref{eqn:exp-determ-2} should be replaced with,
\begin{align}
&\frac{1}{ \pi} \int_\cc ( \del_{\bar{z}} \tilf (z) ) E(z) \d z \d \bar{z} = \left( \sum_{j \neq a } \frac{ s_{aj}^{(4)} + t_{aj}^{(4)}}{N^2} \right) \frac{1}{ 2 \pi} \int_{-2}^2 f(x) \frac{ x^4 - 4x^2 + 2 }{ \sqrt{4-x^2}} \d x \notag\\
+& \left( \sum_a \frac{ s_{aa}^{(3)}}{N^{3/2}} \right) \frac{1}{ 8 \pi} \int_{-2}^2 f(x) \frac{ x^3 - x^2 - 2 x + 4 }{ \sqrt{ 4 -x^2}} \d x =: E_{S, \cc} (f)
\end{align}
The estimate \eqref{eqn:exp-determ-4} is not needed in the complex case. Note that the majority of the proof of Lemma \ref{lem:exp-determ} is devoted to dealing with the $\tr (S (1- \msc(z)^2 S)^{-1} )$ term which is not present in the complex case. Finally, in \eqref{eqn:exp-correction}, the term $E_S (f)$ should be replaced with $E_{S, \cc} (f)$, but Proposition \ref{prop:exp-correction} and its proof is otherwise unchanged.

\subsection{Characteristic function}

We now outline how the arguments of Section \ref{sec:char} change in the complex case. First note that in the real symmetric case we rely on \eqref{eqn:intermediate-characteristic} whereas in the complex case we rely on \eqref{eqn:intermediate-characteristic-cplx}. As in the real case, we must compute the second order terms in \eqref{eqn:intermediate-characteristic-cplx}. In this case, we replace \eqref{eqn:char-second-basic} with,
\begin{align}  \label{eqn:char-second-basic-cplx}
 & \sum_{ja} \frac{s_{aj}}{N} \left( \ee[ \del_{ja} ( \ea G_{ja} (z) ) ] - \ee[ \ea ] \ee[ \del_{ja} G_{ja} (z) ] \right) \notag\\
= & - \frac{ \i \lambda}{ \pi} \int_{\Oma} ( \del_{\bar{w}} \tilf (w) ) \del_w \sum_{ja} \frac{ s_{aj}}{N} \ee[ \ea G_{ja} (z) G_{aj} (w) ] \d w \d \bar{w} \notag\\
- & \sum_{ja} \frac{s_{aj}}{N} \ee[ \ea (G_{aa} (z) G_{jj} (z) - \ee[ G_{aa} (z) G_{jj} (z)]) ] .
\end{align}
For Lemma \ref{lem:char-second-1}, the estimate \eqref{eqn:char-second-1} is unchanged and \eqref{eqn:char-second-2} is unneeded. The proof is identical. The estimates and proofs of Lemma \ref{lem:char-second-2} are also unchanged.  As a consequence of replacing \eqref{eqn:intermediate-characteristic} by \eqref{eqn:intermediate-characteristic-cplx} and replacing \eqref{eqn:char-second-basic} by \eqref{eqn:char-second-basic-cplx}, we have that \eqref{eqn:good-Gaa-equation} of Proposition \ref{prop:char-stein-1} is replaced by,
\begin{align}
 & (z + 2 \msc (z) ) \sum_a \ee[ \ea (G_{aa}(z) - \ee[G_{aa} (z) ] ) ] \notag\\
 = & \ee[ \ea ] \frac{ \i \lambda}{\pi} \int_{\Oma} ( \del_{\bar{w}} \tilf (w) ) \del_w ( -  M \tr (S (1- MS)^{-1}  ) \d w \d \bar{w} \notag\\
 + & \ee[\ea] \bigg\{ - \tilde{s}_{4, \cc} \i \lambda a_{2, \mfa} \msc(z)^2 + \frac{\tilde{s}_3 \i \lambda}{2 N^{1/2}}  \msc(z)  \left[ 2 \msc(z)  a_{1, \mfa} +  \i \lambda a_{1, \mfa}^2 +   a_{2, \mfa} \right] \bigg\} \notag\\
 + & \Osd \left( N^{-1} ( 1 + |\lambda| )^3  + N^{-1} \eta^{-2} + |\lambda \Phi \Psi(z) | \eta^{-1} + |\lambda| \eta^{-1} \Phi^2 (1 + |\lambda| ) \right) 
 \end{align}
where,
\beq
\tilde{s}_{4, \cc} := \sum_{a \neq j } \frac{ s_{aj}^{(4)} + t_{aj}^{(4)}}{N^2}.
\eeq
Consequently, \eqref{eqn:stein-1} still holds after we replace $V_\mfa (f)$ by,
\begin{align}
V_{\mfa, \cc} (f) &:= \frac{1}{ \pi^2} \int_{\Oma^2} ( \del_{\bar{z}} \tilf (z) ) ( \del_{\bar{w}} \tilf (w) ) \del_w ( \msc' (z) \msc (w) \tr (S (1- M S)^{-1} ) ) \d w \d \bar{w} \d z \d \bar{z} \notag\\
&+ \frac{1}{2} \tilde{s}_{4, \cc} a_{2, \mfa}^2 - \frac{ \tilde{s}_3}{N^{1/2}} a_{2, \mfa} a_{1, \mfa} .
\end{align}
We can leave $B_\mfa (f)$ unchanged. 

Next, Sections \ref{sec:cheby-1} and \ref{sec:easy} are unchanged. Some parts of Section \ref{sec:var-functional} change. First, the relevant function replacing $F$ in \eqref{eqn:var-F-def} is,
\beq
F_\cc (z, w):= \del_w ( \msc' (z) \msc(w) \tr ( S(1- \msc(z) \msc(w) S)^{-1} ) ).
\eeq
With this replacement, the statement and proof of Lemma \ref{lem:var-1} are otherwise unchanged. Proposition \ref{prop:var-1} is unchanged. For Lemma \ref{lem:var-2}, we replace $V_1(f)$ with $V_2 (f)$ as defined in \eqref{eqn:intro-V-def}. 
With this definition the analog of \eqref{eqn:Va-V} holds. For the remainder of Lemma \ref{lem:var-2} we need only verify that $V_\cc (f) \geq 0$. For this we write with $a= t_1 (f)$ and $b = t_2 (f)$,
\begin{align}
V_\cc (f)  \geq & \frac{b^2}{2} \left( \sum_{ i \neq j } \frac{ s_{ij}^{(4)} + t_{ij}^{(4)}}{N^2} + \sum_{i \neq j } \frac{ s_{ij}^2}{N^2} \right) \notag\\
+ & \frac{1}{4} a^2 \sum_i \frac{s_{ii}}{N} + \frac{b^2}{4} \left( \sum_i \frac{ s_{ii}^{(4)}}{N^2} + 2 \sum_i \frac{ s_{ii}^2}{N^2} \right) + \frac{ ab}{2 N^{1/2}} \left( \sum_i \frac{ s_{ii}^{(3)}}{N^{3/2}} \right) 
\end{align}
Non-negativity of the second line follows from the first two lines of \eqref{eqn:var-nonneg-2}. The non-negativity of the first line is proven in the same manner as \eqref{eqn:var-nonneg-3} once one realizes that $s_{ij}^2/N^2 = 2 \ee[  \Re[H_{ij} ] ^2]^2 + 2 \ee[  \Im[ H_{ij} ]^2 ]^2$. Finally, Lemma \ref{lem:var-3} is unchanged.

The three estimates of Lemma \ref{lem:stein} are unchanged after replacing $V_1(f)$ with $V_2 (f)$.  Finally, the proofs in Section \ref{sec:real-proofs} are carried out in the same manner.

\section{Linear algebra results}

The following allows us to estimate matrix elements of operators similar to $(1 -\msc(z) \msc (w) S)^{-1}$. 
\bel \label{lem:inverse}
Let $M$ be an $N \times N$ matrix with complex entries such that,
\beq
|M_{ij} | \leq \frac{C}{N}, \qquad \| M \|_{\ell^2 \to \ell^2} \leq 1 - c
\eeq
for some $C, c>0$. Then,
\beq
| (1- M)^{-1}_{ij} | \leq C' ( \delta_{ij} + N^{-1} ).
\eeq
\eel
\proof We have the identity,
\beq
\frac{1}{1-M} = 1 + M + M \frac{1}{1-M} M.
\eeq
To prove the estimate it suffices to show that the matrix entries $M(1-M)^{-1} M$ are bounded by $C'/N$. But this follows since,
\begin{align}
| (M(1-M)^{-1} M)_{ab} | \leq \frac{C^2}{N^2} \sum_{ij} | (1-M)^{-1}_{ij} | \leq C' /N
\end{align}
where the last inequality uses the fact that $\| (1-M)^{-1} \|_{\ell^2 \to \ell^2} \leq C''$ for some $C'' >0$. \qed

\subsection{Derivatives of the resolvent wrt matrix entries}

If we view the resolvent $G = (H-z)$ as a function on the space of real symmetric $N \times N$ matrices, then we have,
\beq \label{eqn:delabgij}
(1 + \delta_{ab} ) \del_{ab} G_{ij} = - G_{ia} G_{bj} - G_{ib} G_{aj} .
\eeq

\section{Various proofs}

\bel \label{lem:G-dif}
We have,
\beq
\frac{1}{N} \sum_j G_{aj}(z) G_{j2}(w) = \frac{1}{N} \frac{ G_{a2} (z) - G_{a2} (w) }{z-w} = \frac{\delta_{a2}}{N} \frac{ \msc(z) - \msc(w)}{z-w} + \Osd \left( \frac{ \Psi (z) + \Psi (w) }{ N( |\Im[z] | + | \Im[w] | ) } \right)
\eeq
\eel
\proof The first equality follows from $(H-z)^{-1} (H-w)^{-1} (z-w) = (H-z)^{-1} - (H-w)^{-1}$. For the estimate, let us assume that $| \Im[z] | \geq | \Im[w] |$, and for simplicity that $\Im[z] >0$, other cases being similar. If $|z-w| > \Im[z]/2$, then the claim follows immediately from \eqref{eqn:entry-wise}, as $|z-w| \geq c \Im[z]$. On the other hand, if $|z-w| < \Im[z]/2$, then the estimate follows from,
\begin{align}
&\frac{G_{a2} (z) - G_{a2} (w)}{z-w} = \int_{0}^1 (\del_u G_{a2} ) (w + s (z-w) ) \d s \notag\\
&= \delta_{a2} \int_0^1 \msc' (w + s (z-w) ) \d s + \Osd ( \Psi(z) / \Im[z] ) 
\end{align}
\qed

\subsection{Proof of Lemma \ref{lem:eps-regular}} \label{a:eps-regular}

$H$ and $M^{(\theta)}$ be as in the statement of the lemma. Let $\Delta_{ab}$ be the matrix that is $1$ only for entries $(a,b)$ and $(b, a)$ and $0$ otherwise. Assuming that $\delta < \frac{1}{4}$ and using that $\| (M^{(\theta)} -z )^{-1} \|_{\ell^2 \to \ell^2} \leq N$ we have by the resolvent expansion and \eqref{eqn:entry-wise} for $H$,
\begin{align}
(M^{(\theta)} -z )^{-1}_{ij} = (H-z )^{-1}_{ij}  + \sum_{m=1}^{20} \theta^m [ (H-z)^{-1} ( \Delta_{ab} (H-z)^{-1} )^m ] _{ij} + \Osd ( N^{-2} )
\end{align}
As long as, say, $\delta < \eps/2$ then every term in the sum is $\Osd ( N^{-1/2+\eps} )$ by applying the entry-wise local law \eqref{eqn:entry-wise}, and so \eqref{eqn:eps-regular-def-2} follows for $M^{(\theta)}$. Moreover, the terms for $m \geq 2$ are $\Osd ( N^{-1+\eps} )$ so to establish \eqref{eqn:eps-regular-def-1} it suffices to to estimate the normalized trace of the term with $m=1$. This term equals,
\beq
\frac{ \theta}{N} \sum_a G_{ia} (z) G_{bi} (z) \sd \frac{ \theta}{N \eta} ( \Im[ G_{aa}(z) ]  + \Im[ G_{bb} (z) ] ) \sd \frac{1}{ N \eta}
\eeq
by the Ward identity. This completes the proof. \qed

\subsection{Proof of \eqref{eqn:HS-est}} \label{a:HS-est}

The difference between the LHS and the first term on the RHS of \eqref{eqn:HS-est} is bounded by,
\beq
\int_{ 0 < y< N^{\mfa-1} } y |f''(x)| N| \Im[ m_N(z) ] - \ee[ \Im[m_N (z) ] ] | \d x \d y
\eeq
Fixing an $\eps >0$ we bound the integral over $ N^{\eps-1} < y < N^{\mfa-1}$ by $\Osd ( N^{\mfa-1} \| f''\|_1 )$ using \eqref{eqn:local-law}. For the integral over $y < N^{\eps-1}$ we can use the fact that $y \to y \Im[ m_N (x + \i y ) ]$ is increasing to bound,
\beq
y \Im[ m_N (x + \i y ) ] \leq N^{\eps-1} \Im[m_N ( x + \i N^{\eps-1} )] \sd N^{\eps-1} 
\eeq
with the second estimate following from \eqref{eqn:local-law}. The claim now follows. \qed

\subsection{Proof of Lemma \ref{lem:H-est}} \label{a:H-est}

From the Cauchy integral formula,
\beq
\del_z H(x + \i y ) = \frac{1}{ 2 \pi \i } \int_{ |x+ \i y - w | = ry } \frac{ H(w)}{ (w- (x+ \i y ) )^2} \d w
\eeq
for any $0 < r < 1$, and so 
\beq \label{eqn:a-H-est}
| \del_z H (x + \i y ) | \leq \frac{C}{ry} \sup_{ |w - (x + \i y ) | = r y } |H(w) | .
\eeq
The estimate \eqref{eqn:del-z-H} follows by this and induction. The estimate \eqref{eqn:H-est} is a special case of \eqref{eqn:general-H-est} and so it remains to prove the latter.  

The terms on the LHS of \eqref{eqn:general-H-est} involving $\i ( f(x) + \i y f'(x) ) \chi'(y)$ are clearly bounded by the first term on the RHS of \eqref{eqn:general-H-est}, and so it suffices to bound the integral involving only $f''(x) y \chi(y)$.  Let $ T = \| f''\|_1^{-1}$. We first estimate,
\beq
\left| \int_{N^{\mfa -1} < |y| < T } y \chi (y) f''(x) H(x + \i y ) \d x \d y \right| \leq C  \| f''\|_1 \sup_{ |x| \leq 10 } \int_{ N^{\mfa-1} < |y| < T }  |y H(x + \i y ) | \d y
\eeq
and so this is bounded by the second term on the RHS of \eqref{eqn:general-H-est}. By integration by parts and the fact that $\del_x H = \del_z H$ by the Cauchy-Riemann equations we have,
\begin{align}
 & \left| \int_{ T < |y| < 2 } \chi (y) y f''(x) H(x) \d x \d  y \right| =  \left| \int_{ T < |y| < 2 } \chi (y) y f'(x)\del_z H(x) \d x \d  y \right| \notag\\
\leq & \| f' \|_1 \sup_{ |x| \leq 10 } \int_{ T < |y| < 2 } |S(x + \i y ) |
\end{align}
where we used \eqref{eqn:a-H-est} in the last inequality (with $r=1/2$). This completes the proof. \qed

\subsection{Semicircle calculations} \label{a:semi-calc}

\subsubsection{Proof of Lemma \ref{lem:msc}}

The estimate \eqref{eqn:msc-upper} follows from \cite[(4.2)]{erdHos2013local}. The estimate \eqref{eqn:msc-difference} is a consequence of \eqref{eqn:msc-dif-identity} and \cite[(4.2)]{erdHos2013local}. 
The identity  \eqref{eqn:msc-dif-identity} follows from \cite[Lemma B.1]{landon2022almost}. \qed


\bibliography{mybib}{}
\bibliographystyle{abbrv}

\end{document}